\documentclass[11pt]{article}
\usepackage{amssymb}
\usepackage{amsmath}
\usepackage{amsthm}
\usepackage{hyperref}
\usepackage{amsfonts}
\usepackage{color}
\usepackage{epsfig}
\usepackage{graphics}                 
\usepackage{graphicx}
\usepackage{float}
\usepackage{epsf}
\usepackage{mathrsfs}
\usepackage{subfigure}

\newtheoremstyle{thmm}{1.5ex plus 1ex minus .2ex}{1.5ex plus 1ex minus
.2ex}{\it}{}{\bfseries}{}{1em}{}
\theoremstyle{thmm}
\newtheorem{theorem}{Theorem}[section]
\newtheorem{algorithm}{Algorithm}[section]
\newtheorem{lemma}{Lemma}[section]

\newtheorem{remark}{Remark}[section]

\allowdisplaybreaks \textwidth  6.0in \textheight 9in
\newcommand{\nn}{\nonumber}

\def\refe#1{(\ref{#1})}

\definecolor{myGreen}{rgb}{0.9, 0.99, 0.9}

\def\e{\epsilon}

\def\d{\delta}

\def\d{\,{\rm d}}
\begin{document}

\date{\today}
\allowdisplaybreaks{}
\title{\bf  A mathematical model and its FEMs for the electro-osmotic flow in micro-channels\thanks{Thi work is supported by the Natural Science Foundation of Henan Province (No. 252300421989) and Innovative Research Team of Henan Polytechnic University, China (No. T2024-4).}}

\author{
	Yunxia Wang\footnotemark[2]\ \ \& Zhiyong Si\footnote{School of Mathematics and Information Science,
		Henan Polytechnic University, 454003, Jiaozuo, P.R. China. {\tt wangyunxia@hpu.edu.cn} (Y. Wang), {\tt sizhiyong@hpu.edu.cn} (Z. Si).}
}

\maketitle

\begin{abstract}
	In this paper, we will provide the finite element method for the electro-osmotic flow in micro-channels, in which a convection-diffusion type equation is given for the charge density $\rho^e$. A time-discrete method based on the backward Euler method is designed. The theoretical analysis shows that the numerical algorithm is unconditionally stable and has optimal convergence rates. To show the effectiveness of the proposed model, some numerical results for the electro-osmotic flow in the T-junction micro-channels and  in rough micro-channels are provided. Numerical results indicate that the proposed numerical method is suitable for simulating electro-osmotic flows. \\
	
\maketitle


\end{abstract}

\section{Introduction}
\setcounter{equation}{0}
With the development of bio-MEMS and bio-NEMS, micro-fluids have been widely studied. They are widely used in modern biological, chemical, and medical analysis, and have brought new revolutionary capabilities to biology, chemistry, and medicine.
The flows in macroscopic scale and microsystems are different from flows in normal scale. There is a great demand for developing numerical methods for micro-fluids, and there are many great works in this area. The simulation of an electroosmotic flow in rectangular micro-scale channel networks was presented by Bianchi et al~\cite{BFG} by using the finite element method. In~\cite{JL}, a micro-pump in which the pumping mechanism is based on MHD principles was derived by Jang and Lee. In~\cite{LP}, a full 3-dimensional conjugated heat transfer model is used to simulate the heat transfer performance of silicon-based, parallel microchannel heat sinks.
A micro-fluidic flow model where the movement of several charged species is coupled with the electric field and the motion of ambient fluid was presented by B. Mohammadi and J. Tuomela~\cite{MT}. Hong and Cheng~\cite{HC} presented a numerical study on laminar forced convection of water in offset strip-fin micro-channels
network heat sinks for microelectronic cooling was presented. In~\cite{MS}, a numerical method that combines a state of control and shape design for the optimization of microfluidic channels was used for sample extraction and separation of chemical species existing in a buffer solution.  An adaptive finite element method with a large aspect ratio for mass transport in electro-osmosis and pressure-driven micro-flows was shown by Prachittham et al.~\cite{PPG}.
Singh and Agrawal~\cite{SA} studied the Burnett equation in cylindrical coordinates and gave the solution in a micro-tube.
 The effect of geometric parameters on water flow and heat transfer characteristics in a microchannel heat sink with triangular reentrant cavities using a numerical method was studied by Xia et al.~\cite{XCWZC}. The heat transfer characteristics of a double-layered micro-channel heat sink was studied by Huang et al.~\cite{HYL} through finite volume method.
By using the FLUENT software, Sui et al. \cite{STL} studied the fully developed flow and heat transfer in periodic wavy channels with rectangular cross sections.
In \cite{HH13}, Ho and Hung presented hybrid finite element and particle-in-cell simulation of the effect for different Debye lengths on charged ion migration in capillary zone electrophoresis. In~\cite{DDES}, the best possible design for a flow distributor at the entrance of a flat channel was established by Davydova et al.
Lockerby and Collyer~\cite{LC} derived the fundamental solutions (Green's functions) to Grad's steady-state linearized 13-moment equations for non-equilibrium gas flows.
Xia et al. \cite{XCMZY} studied micro-PIV visualization and numerical simulation of flow and heat transfer in three micro pin-fin heat sinks.  The performance of ions in electrophoresis microchips with
different crosses comparatively by simulation and experimental methods was given by Yang et al.~\cite{YHCLY}. Pezeshkpour~\cite{PSR} presented a shape factor model for injection analysis of microchip sample electrophoresis.
Abdollahi et al. \cite{ANS} studied the fluid flow and heat transfer of liquid-liquid Taylor flow in square microchannels.
Li et al.~\cite{LRSSAC} gave the numerical approach for nano-fluid transportation due to electric force in a porous enclosure. Wang~\cite{W} presented the numerical investigation for the heat transfer of a droplet-laden flow in a microfluidic system based on the volume of fluid method. In~\cite{LS}, Lin et al. studied the instability of electrokinetic flows with conductivity gradients. Roughness and cavitations effects on electro-osmotic flows in rough microchannels was studied by Wang et al.~\cite{WWC} and Kamali et al.~\cite{KSH}. 
In \cite{SDW}, Si et al. gave a Modified characteristics finite element method for the electro-neutral micro-fluids.

In this paper, we will give a mathematical model and its numerical method for the electro-osmotic flow in microchannels. In this model, a convection-diffusion type equation is given for the charge density $\rho^e$. A time-discrete method based on the backward Euler method is designed. Then, a finite element algorithm for the model is provided and analyzed. Theoretical analysis shows that the numerical algorithm is unconditionally stable and has optimal convergence orders. Numerical results for the electro-osmotic flow in rough microchannels  are given. The effect of fluid viscosity and roughness of the microchannel has been studied in detail. Numerical results confirm the previous numerical results very well.  It indicates that the proposed numerical method is effective for capturing the dynamics of the electro-osmotic flow in microchannels.

\section{The mathematical model and its numerical method}

In this paper, we consider the electroneutral micro-fluids, which are governed by the following systems \cite{LS}.  
\begin{align}\label{equation}
\left\{\begin {array}{rll}
u_t-\mu \triangle u+(u\cdot\nabla)u+\nabla p-\rho^e\nabla \phi=0,&\\
\nabla\cdot u=0,&\\
c^i_t-d_i\triangle c^i -\nu_i z_i \nabla\cdot(c^i\nabla\phi )+u\cdot\nabla c^i=0,& i=1,\ldots,M,\\
-\varepsilon \triangle \phi-\rho^e=0,&\\
\rho^e-\sum_{i=1}^{m} z_ic^i=0,&
\end{array}\right.
\end{align}
where $u$, $p$, $\phi$, are the velocity, pressure, electric potential, $c^i, i=1,2,\ldots, M$ are the molar concentration of $i^{th}$  ionic species in the electrolytes, $\rho^e$ is the charge density, and $D_0$ is the diffusive coefficient of the charge density $\rho^e$, $d_i, z_i, \nu_i$ are the diffusive coefficient, valence number, and charge number of $i^{th}$ ionic species, respectively.

Times the third equation by $z_i, i= 1,
\ldots, M $  and summing them up, using the last equation of (\ref{equation}), we can get
\begin{align*}
	\rho^e_t-\sum_{i=1}^M d_iz_i\triangle  c^i -\sum_{i=1}^M \nu_iz_i^2\nabla\cdot(c^i\nabla \phi)+u\cdot\nabla \rho^e=0.
\end{align*}
Then, equation (\ref{equation}) can be rewritten as follows.

\begin{align}\label{requ}
\left\{\begin {array}{rll}
\rho^e_t-\sum_{i=1}^M  d_iz_i\triangle c^i -\sum_{i=1}^M \nu_iz_i^2\nabla\cdot(c^i\nabla \phi)+u\cdot\nabla \rho^e=0,&\\
u_t-\mu \triangle u+(u\cdot\nabla)u+\nabla p-\rho^e\nabla \phi=0,&\\
\nabla\cdot u=0,&\\
c^i_t-d_i\triangle c^i -\nu_i z_i \nabla\cdot(c^i\nabla\phi )+u\cdot\nabla c^i=0,& i=1,\ldots,M,\\
-\varepsilon \triangle \phi-\sum_{i=1}^{m} z_ic^i=0.&
\end{array}\right.
\end{align}

\begin{remark}
In \cite{LS}, Lin et. al. gave a model by omitting the difference of $d_i, i=1, \ldots, M$. The advantage is that the first equation in (\ref{requ}) is change as a parabolic equation, but there is a disadvantage that they omit the different diffusion between the concentrations.
\end{remark}

\begin{lemma}\cite{G,H03}\label{Lem2.1}
	{\it
	Assume that $\Omega$ is bounded and $\partial \Omega$ is in $C^2$. For any prescribed $g\in L^p(\Omega) (1 < p < 6)$, the steady Stokes problem
	\begin{eqnarray}
	\begin{array}{ll}
	  -\triangle v+\nabla q=g, & \mbox{ in }\Omega, \\
	  \nabla\cdot v=0, & \mbox{ in }\Omega ,\\
	  v=0, & \mbox{ on }\partial\Omega ,
	\end{array}
	\label{stokes}
	\end{eqnarray}
	admits a unique solution satisfying
	\begin{align}
	\|v\|_{W^{2,p}}+\|q\|_{W^{1,p}}\leq C\|g\|_{L^p}.
	\label{stokes-1}
	\end{align}
	}
	\end{lemma}

The initial and boundary conditions are given as follows. The boundary condition for the fluid is slip boundary condition given by~\cite{MT}
\begin{align*}
u=-\xi \nabla \phi.
\end{align*}
Here $\xi >0$ is a constant which depends on the material of the channels and the fluid permittivity and dynamic viscosity. It is the homogeneous Neumann boundary condition on the outlet, and the Dirichlet condition on the inlet. The boundary conditions for the concentration are Dirichlet on the inlet boundary, homogeneous Neumann boundary conditions along the walls, and the outlet boundary. We impose homogeneous Neumann conditions on boundaries for the electric potential.  It can be given as follows \cite{MT}
\begin{align*}
	u|_{\partial\Omega_{in}}&=u_{in},\quad \frac{\partial u}{\partial n}|_{\partial\Omega_{out}}=0,\quad c|_{\partial\Omega_{in}}=c_{in},\quad \frac{\partial c}{\partial n}|_{\partial \Omega/\partial\Omega_{in}}=0,\\
	\phi|_{\partial\Omega_{in}}&=\phi_{in},\quad
	\phi|_{\partial\Omega_{out}}=\phi_{out},\quad
	\phi|_{\partial \Omega/(\Gamma_D)}=0,
\end{align*}
where $\Gamma_D=\partial\Omega_{in}\cup \partial\Omega_{out}$.

For any integer $m\geq 0$ and $1\leq p\leq\infty$, let $W^{m,p}(\Omega)$
be the usual Sobolev space of functions defined in $\Omega$
equipped with the norm  \cite{Adams}
\begin{align*}
	\|f\|_{m,p}
	=\left\{
	\begin{array}{ll}
	\left(\sum_{|\beta|\leq m}\int_\Omega|D^\beta f|^p\d
	x\right)^\frac{1}{p},
	&1\leq p<\infty
	,\\[10pt]
	\displaystyle
	\sum_{|\beta|\leq m}{\rm ess}\,\sup_{x\in
		\Omega}\,|D^\beta
	f(x)|, & p=\infty,
	\end{array}
	\right.
\end{align*}
where
\begin{align*}
D^\beta=\frac{\partial^{|\beta|}}{\partial x_1^{\beta_1}\cdots \partial x_d^{\beta_d}}
\end{align*}
for the multi-index $\beta=(\beta_1,\cdots,\beta_d), \beta_i\geq 0, i=1,\cdots, d$
and $|\beta|=\beta_1+\cdots+\beta_d$.
We denote
\begin{align*}
X&:=\left\{u\in H^1(\Omega)^d; ~  u|_{\partial\Omega_{in}}=u_{in},\quad \frac{\partial u}{\partial n}|_{\partial\Omega_{out}}=0 \right\}\\
M&:=L_0^2(\Omega)=\left\{\varphi\in
L^2(\Omega);\int_{\Omega}\varphi
dx=0\right\} \, ,\\
W&=\left\{ H^1(\Omega); ~ \phi|_{\partial\Omega_{in}} =\phi_{in},\quad
	\phi|_{\partial\Omega_{out}}=\phi_{out},\quad
	\phi|_{\partial \Omega/(\Gamma_D)}=0\right\},
\end{align*}
and
$$
\|\cdot\|_k: = \| \cdot \|_{H^k(\Omega)}, \qquad \|\cdot\|_0: = \| \cdot \|_{L^2(\Omega)}
\, .
$$

Let $0=t_{0}<t_{1}<\cdots<t_{N}=T$ be a uniform partition of the time interval $[0,T]$ with $t_{n}=n\tau$,
and $N$ being a positive integer. In this paper, we set $\phi^{n}=\phi(x,t_{n})$, for any function $\phi(x,t)$. For any sequence of functions $\{f^{n}\}_{n=0}^{N}$, we define
\begin{align*}
D_{\tau}f^{n}=\frac{f^{n}-f^{n-1}}{\tau}.
\end{align*}

Then, a time-discrete method for the system (\ref{requ}) is given as follows:

	\begin{algorithm}[Time-discrete method] \label{AL2.1}
			\textit{Step 1.} Find $\rho^{e,n} \in W $ by
		\begin{align}\label{9}
		&D_\tau \rho^{e,n}-D_0 \triangle\rho^{e,n} -\sum_{i=1}^M \nu_iz_i^2\nabla \cdot (c^{i,n-1}\nabla \phi^{n-1})+ u^{n-1}\cdot\nabla \rho^{e,n}=0,
		\end{align}
		where $D_0=\frac{\sum_{i=1}^M d_i}{M}$.

		\textit{Step 2.} Compute $\phi^{n}\in W$ as the solution of
		\begin{align}\label{10}
		-\varepsilon\triangle \phi^{n} -\rho^{e,n}=0.
		\end{align}
		
		\textit{Step 3.} Find $u^n\in X$, such
		that
		\begin{align}\label{step21}
		&\frac{u^n-u^{n-1}}{\tau}
	-\mu \triangle u^n+u^{n-1}\cdot\nabla u^n+\nabla p^{n}-\rho^{e,n} \nabla \phi^{n}=0,\\
	&\nabla \cdot u^n=0.\label{step22}
		\end{align}
        
		\textit{Step 4.} Find $c^{i,n}\in W$ such that
		\begin{align}\label{12}
		D_\tau c^{i,n}-d_i\triangle c^{i,n} +u^{n}\cdot\nabla c^{i,n}&+ \nu_i z_i \cdot\nabla (c^{i,n}\nabla\phi^{n}) =0,&i=1,\ldots,M.
		\end{align}		
	\end{algorithm}

We denote $T_{h}$ be a regular and quasi-uniform partition of the domain $\Omega$ into the triangles for $d=2$ or tetrahedra for $d=3$
with diameters by a real positive parameter $h (h\rightarrow 0)$. The finite element pair $(X_{h},M_{h},W_{h})$ is constructed based on
$T_{h}$. We assume that $(X_h,M_h)$ satisfies the discrete
LBB condition
\begin{eqnarray}\label{LBB}
\sup_{v_h\in  X_h}\frac{(\varphi_h,\nabla\cdot v_h)}{\|\nabla v_h\|_0}\geq
\beta \|\varphi_h\|_0, \quad \forall \varphi_h\in M_h,
\end{eqnarray}
for some constant $\beta >0$.

We define the Stokes projection $(R_{h}(u,p),Q_{h}(u,P)):(X,M)\rightarrow (X_{h},M_{h})$ by
\begin{align}\label{3.45}
&\nu(R_{h}(u,p)-u,\nabla v_{h})-(Q_{h}(u,p)-p,\nabla \cdot v_{h})=0, \forall v_h\in X_h,\\
&(\nabla\cdot (R_{h}(u,p)-u),\phi_{h})=0,\forall \phi_h\in M_h.
\end{align}
By the classical FEM theory (\cite{AG,BNV1}), we have the following results.
\begin{lemma}\label{3.2}
Assume that $u \in H_{0}^{1}(\Omega)^{d}\cap H^{r+1}(\Omega)^{d}$ and $p \in L_{0}^{2}(\Omega)\cap H^{r}(\Omega)$. Then,
\begin{align}\label{3.46}
\|R_{h}(u,p)-u\|_{0}+h(\|\nabla(R_{h}(u,p)-u)\|_{0}+\|Q_{h}(u,p)-p\|_{0})\leq Ch^{r+1}(\|u\|_{r+1}+\|p\|_{r}),
\end{align}
and
\begin{align}\label{3.47}
\|R_{h}(u,p)\|_{L^{\infty}}\leq C(\|u\|_{2}+\|p\|_{1}).
\end{align}
\end{lemma}
\begin{lemma}\label{3.b}
If $(u,p)\in W^{2,k}(\Omega)^{d}\times W^{1,k}(\Omega)$ for $k>d$,
\begin{align}\label{3.d}
\|\nabla R_{h}(u,p)\|_{L^{\infty}}\leq C(\|u\|_{W^{1,\infty}}+\|p\|_{L^{\infty}}).
\end{align}

\end{lemma}

Then, we define the Ritz projection \cite{GHL} $R_{ch}(c): H_0^1(\Omega)\rightarrow W_{h}$ by
\begin{align}\label{Ritz}
& (\nabla (R_{h,c}(c)-c),\nabla \zeta_h) =0,\forall\ \zeta_h\in W_h,
\end{align}
and it holds that
\begin{align}\label{er-Ritz}
\|R_{h,c}(c)-c\|_0+h\| \nabla (R_{h,c}(c)-c)\|_0\leq Ch^{r+1}\|c\|_{r+1},
\end{align}
for any function $c\in H^{r+1}(\Omega)$.

With the notations above, we provide the finite element method as follows

\begin{algorithm}[Fully-discrete finite element method]\label{FEM}

	\textit{Step 1.} Find $\rho^{e,n}_h\in W_h $ by
	\begin{align}\label{22}
		&\left(D_\tau \rho^{e,n}_h,\psi_h\right)+D_0 (\nabla \rho_h^{e,n},\nabla \psi_h) +\sum_{i=1}^M \nu_iz_i^2(c^{i,n-1}_h\nabla \phi_h^{n-1},\nabla \psi_h)\nn\\
		&+(u_h^{n-1}\cdot\nabla \rho^{e,n}_h, \psi_h)=0, \forall \psi_h\in W_h.
	\end{align}
	
	\textit{Step 2.} Compute $\phi_h^{n}\in W_h$ as the solution of
	\begin{align}\label{23}
		(\varepsilon\nabla \phi_h^{n},\nabla \theta_h) -( \rho_h^{e,n}, \theta_h )=0, \forall \theta_h\in W_h.
	\end{align}
	
	\textit{Step 3.} Find $u_h^n\in X_h$ and $p_h^n\in M_h$, such
	that
	\begin{align}\label{24}
	&\left(\frac{u_h^n-u_h^{n-1}}{\tau
	},v_h\right)+\mu(\nabla u_h^n,\nabla v_h)+(u_h^{n-1}\cdot\nabla u_h^n,v_h)-(\nabla\cdot v_h,p_h^{n})\nonumber\\
	&-(\rho_h^{e,n}\nabla \phi_h^n,v_h)=0, \forall v_h\in X_h,\\
	&(\nabla \cdot u_h^n,q_h)=0,\forall q_h\in M_h.\label{24b}
	\end{align}
		
	\textit{Step 4.} Find $c_h^{i,n}\in W_h$ such that
	\begin{align}\label{25}
		\left(D_\tau c_h^{i,n},\zeta_h \right)+d_i(\nabla c_h^{i,n+1},\nabla \zeta_h) +(u_h^{n}\cdot\nabla c_h^{i,n},\zeta_h)&+ \nu_i z_i ( c_h^{i,n}\nabla\phi_h^{n-1} ,\nabla \zeta_h)=0,\nonumber\\ & i=1, 2, \cdots, M, \forall \zeta_h\in W_h.
	\end{align}

\end{algorithm}

\section{Theoretical analysis}
\subsection{Theoretical analysis for the time-semi discrete algorithm}

In order to give the error estimation, we define the error as follows
\begin{align*}
	e_\rho^n=\rho^{e,n}-\rho(t_n), & ~~ e_\phi^n=\phi^n-\phi(t_n), &  e_c^{i,n}= c^{i,n}-c^i(t_n),\\
	e_n^n=u^n-u(t_n),&~~ e_p^n=p^n-p(t_n), & n =  0, 1, 2, \cdots, N.
\end{align*}
\begin{theorem}\label{Th3.2}
	Suppose $c^{i,n}\in W$, $\phi^n\in W$, $u\in X$ and $p^n\in M$ be the solutions of (\ref{AL2.1}), when $\tau$ is sufficient small, we get the error estimates as follows
	\begin{align*}
		&\|e_\rho^n\|_0^2+\tau D_0 \sum_{n=1}^N \|\nabla e_\rho^n\|_0^2 + \|e_c^{i,n}\|_0^2 +\tau\sum_{n=1}^N \sum_{i=1}^M d_i\|\nabla e_c^{i,n}\|_0^2+\|e_u^n\|_0^2
		+ \tau \mu \sum_{n=1}^N \|\nabla e_u^n\|_0^2
	\leq  C\tau^2,\\
	&\tau\sum_{m=0}^n\|e_u^{m}\|_2^2
		+\tau\sum_{m=0}^n\| e_p^m \|_1^2\leq C\tau^2.
	\end{align*}
	Furthermore, there holds that
	\begin{align*}
		\tau\sum_{n=0}^N \left (
\| D_{\tau} u^{n} \|_2^2 + \| D_{\tau} p^{n} \|_1^2
+\| u^{n} \|_{W^{2,d^*}}^2+\| p^{n} \|_{W^{1,d^*}}^2
\right )
+\max_{0\leq n\leq N-1}(\|u^{n}\|_2  + \| p^{n} \|_{1})
\leq C_0^*\, ,
	\end{align*}
and	
	\begin{align*}
		&\| e_\rho^{e,n} \|_2 + \tau^{3/4} \| \rho^{e,n} \|_{W^{2,d^*}} \le 1 \, ,\\
		&\| e_c^{i,n} \|_2 + \tau^{3/4} \| c^{i,n} \|_{W^{2,d^*}} \le 1, i=1,\ldots, M,
	\end{align*}
	where $C_0^*$ is a positive constant independent of $h, \tau$ and  $d^*>d$.
\end{theorem}

\begin{proof}
	We prove it by induction method. It is easy to see that
	\begin{align*}
		&\| e_\rho^{e,0} \|_2 + \tau^{3/4} \| \rho^{e,0} \|_{W^{2,d^*}} \le 1 \, ,\\
		&\| e_c^{i,0} \|_2 + \tau^{3/4} \| c^{i,0} \|_{W^{2,d^*}} \le 1, i=1,\ldots, M.	
	\end{align*}
	hold at the initial time step.
We assume that there holds for $0\leq n\leq k$
for some integer $n\ge 0$
\begin{align*}
	&\| e_\rho^{e,n-1} \|_2 + \tau^{3/4} \| \rho^{e,n-1} \|_{W^{2,d^*}} \le 1 \, ,\\
	&\| e_c^{i,n-1} \|_2 + \tau^{3/4} \| c^{i,n-1} \|_{W^{2,d^*}} \le 1, i=1,\ldots, M.	
\end{align*}

	Subtracting the first equation of (\ref{requ}) from (\ref{9}), we get the error equation of $\rho^e$ in the weak form as follows
	\begin{align}
		&\left(D_\tau e_\rho^{n},\psi\right)+D_0 (\nabla \rho^{e,n},\nabla \psi)+\sum_{i=1}^M  d_iz_i (\nabla (c^{i}(t_{n-1})-c^i(t_n)),\nabla \psi) \nonumber\\
		&+\sum_{i=1}^M \nu_iz_i^2(e_c^{i,n-1}\nabla \phi^{n-1},\nabla \psi)
		+\sum_{i=1}^M \nu_iz_i^2(c^{i}(t_{n-1})\nabla e_\phi^{n-1},\nabla \psi)\nn\\
		&+ \sum_{i=1}^M \nu_iz_i^2(c^{i}(t_{n-1})\nabla (\phi(t_{n-1})-\phi(t_n)),\nabla \psi)+ (e_u^n\cdot\nabla \rho^{e,n})+(u(t_{n-1})\cdot\nabla e_\rho^{n},\psi)\nonumber\\
		&+((u(t_{n-1})-u(t_n))\cdot\nabla \rho^{e}(t_n), \psi)=0, \forall \psi\in W.
	\end{align}
	Letting $\psi=2 \tau e_\rho^n$ and using $2(a-b,b)=\|a\|_0^2-\|b\|_0^2+\|a-b\|_0^2$, we deduce that
	\begin{align}
		&\|e_\rho^n\|_0^2-\|e_\rho^{n-1}\|_0^2+\|e_\rho^n-e_\rho^{n-1}\|_0^2+2\tau D_0 \|\nabla e_\rho^n\|_0^2+2\tau \sum_{i=1}^M \nu_iz_i^2(e_c^{i,n-1}\nabla \phi^{n-1},\nabla e_\rho^n)\nonumber\\
		&
		+2\tau \sum_{i=1}^M \nu_iz_i^2(c^{i}(t_{n-1})\nabla e_\phi^{n-1}, \nabla e_\rho^n)+2\tau  \sum_{i=1}^M \nu_iz_i^2(c^{i}(t_{n-1})\nabla (\phi(t_{n-1})-\phi(t_n)),\nabla e_\rho^n)\nn\\
		&+2\tau  (e_u^n\cdot\nabla \rho^{e,n},e_\rho^n)+2\tau (u(t_{n-1})\cdot\nabla e_\rho^{n},e_\rho^n)+2\tau ((u(t_{n-1})-u(t_n))\cdot\nabla \rho^{e}(t_n), e_\rho^n)=Tr_\rho,
	\end{align}
	where $Tr_\rho= \frac{\partial \rho^{e,n}}{\partial t}(t_n)- \frac{\rho^{e,n}(t_n)-\rho^{e,n-1} (t_{n-1})} {\tau } $. By (\ref{10}), we have $\|\phi^{n-1}\|_{H^2}\leq C\|\rho^{e,n-1}\|_0^2$, using Cauchy-Schwarz and Young's inequality, we deduce
	\begin{align*}
		|2\tau \sum_{i=1}^M \nu_iz_i^2(e_c^{i,n-1}\nabla \phi^{n-1},\nabla e_\rho^n)|\leq &C\tau \| \phi^{n-1}\|_{H^2}^2\sum_{i=1}^M \|e_c^{i,n-1}\|_0^2+ \frac{D_0\tau}{8}\|\nabla e_\rho^n\|_0^2\\
		\leq& C\tau \| \rho^{e,n-1}\|_0^2\sum_{i=1}^M \|e_c^{i,n-1}\|_0^2+ \frac{D_0\tau}{8}\|\nabla e_\rho^n\|_0^2.
	\end{align*}
	Subtracting (\ref{requ}) from (\ref{10}), we can get
	\begin{align*}
		-\varepsilon\triangle e_\phi^{n} -e_\rho^{n}=0.
	\end{align*}
	Taking inner product of it with $e_\phi^n$, we have
	\begin{align*}
		\|\nabla e_\phi^{n}\|_0 \leq \|e_\rho^{n}\|_0.
	\end{align*}
	Using Cauchy-Schwarz and Young's inequality, we derive that
	\begin{align*}
		|2\tau \sum_{i=1}^M \nu_iz_i^2(c^{i}(t_{n-1})\nabla e_\phi^{n-1}, \nabla e_\rho^n)|\leq& C\tau \|\nabla e_\phi^{n-1}\|_0^2+\frac{D_0\tau}{8}\|\nabla e_\rho^n\|_0^2\\
		\leq& C\tau \| e_\rho^{n-1}\|_0^2+\frac{D_0\tau}{8}\|\nabla e_\rho^n\|_0^2.
	\end{align*}
	Using Cauchy-Schwarz inequality, Young's inequality and Taylor's formula, we have
	\begin{align*}
		|2\tau  \sum_{i=1}^M \nu_iz_i^2(c^{i}(t_{n-1})\nabla (\phi(t_{n-1})-\phi(t_n)),\nabla e_\rho^n)|\leq&  C\tau\|\nabla (\phi(t_{n-1})-\phi(t_n))\|_0^2+\frac{D_0\tau}{4}\|\nabla e_\rho^n\|_0^2\\
		\leq&  C\tau^3+\frac{D_0\tau}{8}\|\nabla e_\rho^n\|_0^2.
	\end{align*}
	Using the Cauchy-Schwarz and Young's inequality, we arrive at
	\begin{align*}
		|2\tau  (e_u^n\cdot\nabla \rho^{e,n},e_\rho^n)|\leq C\tau \|e_u^n\|_0^2+\frac{D_0\tau}{4}\|\nabla e_\rho^n\|_0^2.
	\end{align*}
	Noting $\nabla \cdot u(t_{n-1})=0$, we derive that
	\begin{align*} 
		2\tau (u(t_{n-1})\cdot\nabla e_\rho^{n}, e_\rho^n) \leq \frac{D_0\tau}{8} \|\nabla e_\rho^n\|_0^2 +C\tau \|e_\rho^n\|_0^2.
	\end{align*}
	Using Cauchy-Schwarz inequality, Young's inequality and Taylor's formula, we can deduce
	\begin{align*}
		|2\tau ((u(t_{n-1})-u(t_n))\cdot\nabla \rho^{e}(t_n), e_\rho^n) |\leq C\tau^3+\frac{D_0\tau}{8}\|\nabla e_\rho^n\|_0^2.
	\end{align*}
	Then, we arrive at
	\begin{align}\label{26}
		&\|e_\rho^n\|_0^2-\|e_\rho^{n-1}\|_0^2+\|e_\rho^n-e_\rho^{n-1}\|_0^2+\tau D_0 \|\nabla e_\rho^n\|_0^2 \nonumber\\
	\leq &  C\tau \| \rho^{e,n-1}\|_0^2 \sum_{i=1}^M \|e_c^{i,n-1}\|_0^2+C\tau \| e_\rho^{n-1}\|_0^2+C\tau \|e_u^n\|_0^2+C\tau^3.
	\end{align}

	Subtracting second equation of (\ref{requ}) with $t=t_n$ from (\ref{step21}), we get the error equation of $u$ as follows
	\begin{align}\label{eu1}
		&D_\tau e_u^n
	-\mu \triangle e_u^n+u^{n-1}\cdot\nabla u^n-u(t_n)\cdot\nabla u(t_n)+\nabla e_p^{n}-\rho^{e,n} \nabla \phi^{n}+\rho^{e}(t_n) \nabla \phi(t_n)=Tr_u,\\
	&\nabla \cdot e_u^n=0,\label{eu2}
		\end{align}
	where $Tr_u= \frac{\partial u}{\partial t}(t_n)- \frac{u(t_n)-u(t_{n-1})}{\tau}$. Taking inner product of it with $2\tau e_u^n$, we have
	\begin{align*}
		& \|e_u^n\|_0^2-\|e_u^{n-1}\|_0^2+\|e_u^n-e_u^{n-1}\|_0^2
	+ 2\tau \mu \|\nabla e_u^n\|_0^2 +2\tau (u^{n-1}\cdot\nabla e_u^n,e_u^n) +2\tau (e_u^{n-1}\cdot\nabla u^n ,e_u^n) \\
	&+2\tau ((u(t_{n-1})-u(t_n))\cdot\nabla u^n,e_u^n)-2\tau (e_\rho^{n} \nabla \phi^{n},e_u^n) -2\tau (\rho^{e}(t_n)\nabla e_{\phi}^n,e_u^n) =2\tau (Tr_u,e_u^n).
\end{align*}
Using Cauchy-Schwarz and Young's inequality, there holds that
\begin{align*}
	2\tau |(e_u^{n-1}\cdot\nabla u^n,e_u^n)|\leq C\tau \|e_u^{n-1}\|_0^2+ \frac{\nu\tau}{10}\|\nabla e_u^n\|_0^2.
\end{align*}
Noting $\nabla \cdot u(t_{n-1})=0$, it yields that
\begin{align*} 
	2\tau (u(t_{n-1})\cdot\nabla e_u^n,e_u^n)\leq \frac{\nu \tau}{10} \|\nabla e_u^n\|_0^2 +C\tau \|e_u^n\|_0^2.
\end{align*}
By Taylor's formula and Cauchy-Schwarz inequality, we deduce that
\begin{align*}
	2\tau |((u(t_{n-1})-u(t_n))\cdot\nabla u(t_n),e_u^n)|\leq C\tau^3+ \frac{\nu\tau}{5}\|\nabla e_u^n\|_0^2.
\end{align*}
Using Cauchy-Schwarz and Young's inequality, we have
\begin{align*}
	2\tau |(e_\rho^{n} \nabla \phi^{n},e_u^n)|\leq C\tau \|e_\rho^n\|_0^2+ \frac{\nu\tau}{5}\|\nabla e_u^n\|_0^2.
\end{align*}
By Cauchy-Schwarz and Young's inequality, it yields that
\begin{align*}
	2\tau |(\rho^{e}(t_n)\nabla e_{\phi}^n,e_u^n) |\leq  C\tau \|e_\phi^n\|_0^2+ \frac{\nu\tau}{5}\|\nabla e_u^n\|_0^2.
\end{align*}
Using Cauchy-Schwarz and Young's inequality, there holds that
\begin{align*}
	2\tau |(Tr_u,e_u^n)|\leq C\tau^3+ \frac{\nu\tau}{5}\|\nabla e_u^n\|_0^2.
\end{align*}
Then, we arrive at
\begin{align}\label{27}
	& \|e_u^n\|_0^2-\|e_u^{n-1}\|_0^2+\|e_u^n-e_u^{n-1}\|_0^2
+ \tau \mu \|\nabla e_u^n\|_0^2 \leq  C\tau \|e_u^{n-1}\|_0^2+C\tau \|e_\rho^n\|_0^2+C\tau \|e_\phi^n\|_0^2 +C\tau^3.
\end{align}

Subtracting the fourth equation of (\ref{requ}) from (\ref{12}), we get the error equation of $c^i$ as follows
\begin{align*}
	D_\tau e_c^{i,n}-d_i\triangle e_c^{i,n} +u^{n}\cdot\nabla c^{i,n}-u(t_n)\cdot\nabla c^i(t_n) + \nu_i z_i (\nabla\cdot( c^{i,n}\nabla\phi^{n})-\nabla\cdot (c^{i}(t_n)\nabla\phi(t_n))) =Tr_c^i,
\end{align*}
where $Tr_c^i=D_\tau c^{i,n}-\frac{\partial c^{i,n}}{\partial t}, ~i=1,\ldots,M.$ Taking inner product of it with $2\tau e_c^{i,n}$, we deduce that
\begin{align*}
	&\|e_c^{i,n}\|_0^2-\|e_c^{i,n-1}\|_0^2 +2\tau d_i\|\nabla e_c^{i,n}\|_0^2 +2\tau (e_u^{n}\cdot\nabla c^{i,n}, e_c^{i,n})+2\tau (u(t_n)\cdot\nabla e_c^{i,n},e_c^{i,n})\\
	& - 2\tau \nu_i z_i ( e_c^{i,n}\nabla\phi^{n},\nabla e_c^{i,n})+2\tau (c^{i}(t_n)\nabla e_\phi^n,\nabla e_c^{i,n}) =2\tau (Tr_c^i,e_c^{i,n}).
\end{align*}
Using Cauchy-Schwarz inequality and Young's inequality, we have
\begin{align*}
	2\tau |(e_u^{n}\cdot\nabla c^{i,n}, e_c^{i,n})|\leq c\tau\|e_n^n\|_0^2+\frac{\tau d_i}{8}\|\nabla e_c^{i,n}\|_0^2.
\end{align*}
Noting $\nabla \cdot u(t_n)=0$, we deduce
\begin{align*} 
	2\tau (u(t_n)\cdot\nabla e_c^{i,n},e_c^{i,n})\leq \frac{d_i\tau}{8} \|\nabla e_c^{i,n}\|_0^2 +C\tau \|e_c^{i,n}\|_0^2.
\end{align*}
Using Cauchy-Schwarz inequality and Young's inequality, we have
\begin{align*}
	2\tau \nu_i z_i |( e_c^{i,n}\nabla\phi^{n},\nabla e_c^{i,n})|\leq C\tau\|e_c^{i,n}\|_0^2+\frac{\tau d_i}{4}\|\nabla e_c^{i,n}\|_0^2.
\end{align*}
\begin{align*}
	2\tau |(c^{i}(t_n)\nabla e_\phi^n,\nabla e_c^{i,n})|\leq C\tau\|e_\phi^{n}\|_0^2+\frac{\tau d_i}{4}\|\nabla e_c^{i,n}\|_0^2.
\end{align*}
\begin{align*}
	2\tau |(Tr_c^i,e_c^{i,n})|\leq C\tau^3+ \frac{\tau d_i}{4}\|\nabla e_c^{i,n}\|_0^2.
\end{align*}
Then, we arrive at
\begin{align*}
	&\|e_c^{i,n}\|_0^2-\|e_c^{i,n-1}\|_0^2 +\tau d_i\|\nabla e_c^{i,n}\|_0^2 \leq c\tau\|e_n^n\|_0^2 + C\tau\|e_c^{i,n}\|_0^2+ C\tau\|e_\phi^{n}\|_0^2+C\tau^3.
\end{align*}
Taking sum over all $i$, we have
\begin{align}\label{28}
	&\sum_{i=1}^M\|e_c^{i,n}\|_0^2-\sum_{i=1}^M\|e_c^{i,n-1}\|_0^2 +\tau\sum_{i=1}^M d_i\|\nabla e_c^{i,n}\|_0^2 \leq c\tau\|e_n^n\|_0^2 + C\tau\sum_{i=1}^M\|e_c^{i,n}\|_0^2+ C\tau\|e_\phi^{n}\|_0^2+C\tau^3.
\end{align}

Taking sum of (\ref{26}), (\ref{27}) and (\ref{28}), it yields that
\begin{align}\label{29}
	&\|e_\rho^n\|_0^2-\|e_\rho^{n-1}\|_0^2+\|e_\rho^n-e_\rho^{n-1}\|_0^2+\tau D_0 \|\nabla e_\rho^n\|_0^2 +\sum_{i=1}^M\|e_c^{i,n}\|_0^2-\sum_{i=1}^M\|e_c^{i,n-1}\|_0^2 \nonumber\\
	&+\tau\sum_{i=1}^M d_i\|\nabla e_c^{i,n}\|_0^2+\|e_u^n\|_0^2-\|e_u^{n-1}\|_0^2+\|e_u^n-e_u^{n-1}\|_0^2
	+ \tau \mu \|\nabla e_u^n\|_0^2 \nonumber\\
\leq &  C\tau \sum_{i=1}^M \|e_c^{i,n-1}\|_0^2+C\tau \| e_\rho^{n-1}\|_0^2+C\tau \|e_u^n\|_0^2+ C\tau\|e_\phi^{n}\|_0^2+C\tau^3.
\end{align}
Taking sum of (\ref{29}) over all $n$ and using Gronwall's formula, when $\tau$ is sufficient small, we derive that
\begin{align*}
	\|e_\rho^n\|_0^2+\tau D_0 \sum_{n=1}^N \|\nabla e_\rho^n\|_0^2 + \|e_c^{i,n}\|_0^2 +\tau\sum_{n=1}^N \sum_{i=1}^M d_i\|\nabla e_c^{i,n}\|_0^2+\|e_u^n\|_0^2
	+ \tau \nu \sum_{n=1}^N \|\nabla e_u^n\|_0^2
\leq  C\tau^2.
\end{align*}

To obtain an $H^1$-estimate, we multiply \refe{eu1} by $2\tau D_{\tau}e_u^{n+1}$
and integrate it over $\Omega$
to get
\begin{align}
&\mu(\|e_u^{n+1}\|_1^2-\|e_u^n\|_1^2)+2\tau\|D_{\tau}e_u^{n+1}\|_0^2\nn\\
\leq& 2\tau |(u^{n-1}\cdot\nabla e_u^n,D_\tau e_u^n)| +2\tau |(e_u^{n-1}\cdot\nabla u(t_n),D_\tau e_u^n)| +2\tau |((u(t_{n-1})-u(t_n))\cdot\nabla u(t_n) , D_\tau e_u^n)|\nonumber\\
&+2\tau |(e_\rho^{n} \nabla \phi^{n},D_\tau e_u^n)| + 2\tau |(\rho^{e}(t_n)\nabla e_{\phi}^n,D_\tau e_u^n)|
+2\tau|(Tr_u,D_{\tau}e^{n+1})|.
\end{align}
Then, we have
\begin{align*}
	2\tau |(u^{n-1}\cdot\nabla e_u^n,D_\tau e_u^n)| \leq & C\tau \|u^{n-1}\|_{\infty}\| \nabla e_u^n \|_0 \| D_\tau e_u^n\|_0\\
	\leq & C\tau \| \nabla e_u^n \|_0^2 + \frac{\tau }{8}\|D_\tau e_u^n\|_0^2.
\end{align*}
Using Cauchy-Schwarz and Young's inequality, there holds that
\begin{align*}
	2\tau |(e_u^{n-1}\cdot\nabla u(t_n),D_\tau e_u^n)|\leq & C\tau \| e_u^{n-1}\|_0\| u(t_n)\|_{W^{1,\infty}}\| D_{\tau} e_u^n\|_0\\
	\leq & C\tau \|e_u^{n-1}\|_0^2 + \frac{\tau}{8} \| D_\tau e_u^n\|_0^2.
\end{align*}
Using Cauchy-Schwarz, Young's inequality and Taylor's formula, it follows that
\begin{align*}
	2\tau |((u(t_{n-1})-u(t_n))\cdot\nabla u(t_n) , D_\tau e_u^n)| \leq & C\|u(t_{n-1})-u(t_n)\|_0\| u(t_n)\|_{W^{1,\infty}}\| D_\tau e_u^n\|_0\\
	\leq & C\tau \| u(t_{n-1})-u(t_n)\|_0^2 +\frac{\tau}{8} \| D_\tau e_u^n\|_0^2\\
	\leq & C\tau^3+\frac{\tau}{8} \| D_\tau e_u^n\|_0^2.
\end{align*}
Using (\ref{10}), we can deduce that $\| \phi^n\|_\infty\leq \|\phi^n\|_2\leq \|\rho^{e,n}\|_0$. Then, there holds that
\begin{align*}
	2\tau |(e_\rho^{n} \nabla \phi^{n},D_\tau e_u^n)| \leq & C\tau \|\nabla e_\rho^{n}\|_0\| \phi^n\|_\infty \|D_\tau e_u^n\|_0\\
	\leq & C\tau\|\nabla e_\rho^{n}\|_0^2 +\frac{\tau}{8} \| D_\tau e_u^n\|_0^2.
\end{align*}
\begin{align*}
	2\tau |(\rho^{e}(t_n)\nabla e_{\phi}^n,D_\tau e_u^n)|\leq & C\tau \|\rho^{e}(t_n)\|_\infty \|\nabla e_{\phi	}^n\|_0\|D_\tau e_u^n\|_0\\
	\leq & C\tau \| \nabla e_\phi^n\|_0^2+ \frac{\tau}{8} \| D_\tau e_u^n\|_0^2.
\end{align*}
It follows that
\begin{align*}
\mu \|e_u^{n+1}\|_1^2-\mu \|e_u^n\|_1^2+\tau\|D_{\tau}e_u^{n+1}\|_0^2
\leq C\tau \| \nabla e_u^n \|_0^2  + C\tau \| \nabla e_u^{n-1} \|_0^2+  C\tau\|\nabla e_\rho^{n}\|_0^2+ C\tau \| \nabla e_\phi^n\|_0^2 +C\tau^3,
\end{align*}
which in turn produces
\begin{align}
\max_{1\leq m\leq N}\mu \|e_u^{m}\|_1^2+\tau\sum_{m=1}^N\|D_{\tau}e_u^{m}\|_0^2
\leq C \tau^2.
\label{e-h1-error}
\end{align}

Moreover, applying Lemma \ref{Lem2.1} to the equations
\refe{eu1}-\refe{eu2} with $p=2$, we arrive at
\begin{align*}
&\| e_u ^{n}\|_2+\|e_p^n\|_1\nn\\
\leq&
C\|D_{\tau}e_u^{n}\|_0+C\|u^{n-1}\cdot\nabla u^n-u(t_n)\cdot\nabla u(t_n)\|_0
+C\|\rho^{e,n} \nabla \phi^{n}-\rho^{e}(t_n) \nabla \phi(t_n)\|_0+C\|Tr_u^n\|_0\nn\\
\leq&
C\|D_{\tau}e^{n+1}\|_0+C\|e_u^n\|_1\|u^{n-1}\|_{L^{\infty}}
+C\|\nabla u(t_n)\|_{L^4}\|u(t_n)-u(t_{n-1})\|_{L^4}
\nonumber\\
&+C\|e_u^{n-1}\|_{L^4}\|\nabla u(t_n)\|_{L^4}+C\|e_\rho^{e,n}\|_1\|\phi(t_n)\|_{L^{\infty}}+C\|\nabla \rho^{e}(t_n) \|_{L^4}\|e_\phi^n\|_{L^4}
+C \|Tr_u^{n}\|_0
\nn\\
\leq& C\|D_{\tau}e^{n+1}\|_0+C\|e_u^n\|_1+C\|e_u^{n-1}\|_1 +C\|e_\rho^{e,n}\|_1 +C\|e_\phi^n\|_1+C\tau,
\nn
\end{align*}
which together with \refe{e-h1-error} implies
\begin{align}
\tau\sum_{m=1}^N\|e_u^{m}\|_2^2
+\tau\sum_{m=1}^N\| e_p^m \|_1^2\leq C\tau^2. \label{e-h2-error}
\end{align}
From \refe{e-h1-error} and \refe{e-h2-error},
we can see that
when $\tau \le \tau_3$ for some $\tau_3>0$,
\begin{align}
& \max_{0\leq m\leq n}\|u^{m}\|_2
\leq \max_{0\leq m\leq n}(\|u(t_m)\|_2+\|e_u^{m}\|_2)\leq C
\, , \label{Uh2}\\
& \max_{0\leq m\leq n}\|p^{m}\|_1
\leq \max_{0\leq m\leq n}(\|p(t_m)\|_1+\|e_p^m\|_1)\leq C
\, ,\\
& \tau\sum_{m=0}^{n} \| D_{\tau} u^{m} \|_2^2
\leq 2\tau\sum_{m=0}^{n} (\| D_{\tau} u(t_m) \|_2^2+\| D_{\tau} e_u^{m+1} \|_2^2)\leq C
\, , \label{dtU}\\
& \tau\sum_{m=0}^{n} \| D_{\tau} p^{m} \|_1^2
\le 2\tau\sum_{m=0}^{n} (\| D_{\tau} p(t_m) \|_1^2+\| D_{\tau}e_p^m \|_1^2)
\le C
\, .
\end{align}

Again, we apply Lemma \ref{Lem2.1} to the Stokes equation (\ref{step21}) and (\ref{step22}) with $p=d^*$,
and we get
\begin{align*}
\|u^{n}\|_{W^{2,d^*} }+ \|p^{n}\|_{W^{1,d^*}}
&
\leq C\|D_\tau u^{n+1}\|_{L^{d^*}}+\|u^{n-1}\cdot\nabla u^n\|_{L^{d^*}}
+\frac{C}{\tau}\left\| \rho^{e,n} \nabla \phi^{n} \right\|_{L^{d^*}}\nn\\
&
\leq C\|D_\tau u^{n}\|_{L^{d^*}} +  C \| \nabla u^n \|_{L^{d^*}} \| u^{n-1} \|_{L^\infty}+C\|\rho^{e,n}\|_{\infty}\|\nabla \phi^n\|_{L^{d^*}}\, .
\end{align*}
By \refe{Uh2} and \refe{dtU}, it yields that
\begin{eqnarray}
\tau\sum_{m=1}^N(\|u^{m}\|_{W^{2,d^*}}^2+\|p^{m}\|_{W^{1,d^*}}^2)\leq C\, .
\end{eqnarray}
By \refe{e-h2-error} and the above inequality, there exists $\tau_4>0$ such that when
$\tau\leq \tau_4$,
$$
\| e_u^{n} \|_2 + \tau^{3/4} \| u^{n} \|_{W^{2,d^*}} \le 1 \, .
$$

Similarly, we can prove
\begin{align*}
	& \max_{1\leq m\leq N}\|\rho^{e,m}\|_2
	\leq C
	\, , \\
	& \max_{1\leq m\leq N}\|c^{i,m}\|_1
	\leq C, i=1,\ldots, M\
	\, ,\\
	& \tau\sum_{m=1}^{N} \| D_{\tau} \rho^{e,m} \|_2^2
	\leq C
	\, , \\
	& \tau\sum_{m=1}^{N} \| D_{\tau} c^{i,m} \|_1^2
	\leq C, i=1,\ldots, M
	\, ,\\
	& \tau\sum_{m=1}^N\|\rho^{e,m}\|_{W^{2,d^*}}^2+\tau\sum_{m=1}^N\|c^{i,m}\|_{W^{2,d^*}}^2\leq C\, .
\end{align*}
Furthermore, there holds that
\begin{align*}
&\| e_\rho^{e,n} \|_2 + \tau^{3/4} \| \rho^{e,n} \|_{W^{2,d^*}} \le 1 \, ,\\
&\| e_c^{i,n} \|_2 + \tau^{3/4} \| c^{i,n} \|_{W^{2,d^*}} \le 1, i=1,\ldots, M.	
\end{align*}
Thus, the induction is closed.
\end{proof}

\subsection{Theoretical analysis for the finite element algorithm}

In order to give the error estimation, we define the error as follows
\begin{align*}
	e_{h,\rho}^n=\rho_h^{e,n}-R_{h,\rho}(\rho^{e,n}); &~~ e_{h,\phi}^n=\phi_h^n-R_{h,\phi}(\phi^n);& e_{h,c^i}^{n}= c_h^i-R_{h,c} (c^{i,n});\\
	e_n^n=u_h^n-R_h(u^n,p^n);&~~ e_{h,p}^n=p_h^n-Q_h(u^n, p^n).
\end{align*}
\begin{theorem}\label{Th3.1}
	Suppose $c_h^{i,n}\in W_h$, $\phi_h^n\in W_h$, $u_h^n\in V$ and $p_h^n\in Q_h$ be the solutions of (\ref{AL2.1}), when $\tau$ is sufficient small, we get the error estimates as follows
	\begin{align}
		\|e_{h,\rho}^{n}\|_0^2+D_0 \tau \sum_{n=1}^N\|\nabla e_{h,\rho}^{n}\|_0^2
		&+ \|e_{h,u}^N\|_0^2+2\mu\tau \sum_{n=1}^N \|\nabla e_{h,u}^n \|_0^2\nonumber\\
		&+\sum_{i=1}^M\|e_{h,c^i}^{n}\|_0^2+\tau d_i\sum_{n=1}^N\sum_{i=1}^M\|\nabla e_{h,c^i}^{n}\|_0^2
		\leq  C h^{4}.
	\end{align}
	Furthermore, we have
	\begin{align}\label{3.55}
		\max_{0\leq m\leq N}\| u_h^{m} \|_{L^\infty}+\tau \sum_{m=0}^N\| u_h^{m} \|_{W^{1,\infty}}^2 &\le C,\\
		\max_{0\leq m\leq N}\| \rho_h^{e,m} \|_{L^\infty}+\tau \sum_{m=0}^N\| \rho_h^{e,m} \|_{W^{1,\infty}}^2 &\le C,\\
		\max_{0\leq m\leq N}\| c_h^{i,m} \|_{L^\infty}+\tau \sum_{m=0}^N\| c_h^{i,m} \|_{W^{1,\infty}}^2 &\le C, i=1, \ldots, M.\label{3.57}
	\end{align}
\end{theorem}
\begin{proof}
	Now, we prove this theorem by mathematical induction. It is easy to see that
\refe{3.55} to \refe{3.57} hold at the initial time step.
We assume that there holds for $0\leq n\leq k$
for some integer $n\ge 0$
\begin{align*}
	\max_{0\leq m\leq n-1}\| u_h^{m} \|_{L^\infty}+\tau \sum_{m=0}^{n-1}\| u_h^{m} \|_{W^{1,\infty}}^2 &\le C,\\
	\max_{0\leq m\leq n-1}\| \rho_h^{e,m} \|_{L^\infty}+\tau \sum_{m=0}^{n-1}\| \rho_h^{e,m} \|_{W^{1,\infty}}^2 &\le C,\\
	\max_{0\leq m\leq n-1}\| c_h^{i,m} \|_{L^\infty}+\tau \sum_{m=0}^{n-1}\| c_h^{i,m} \|_{W^{1,\infty}}^2 &\le C, i=1, \ldots, M.
\end{align*}

	Using Green's formula, we can deduce the weak form of (\ref{9}) as follows
	\begin{align}\label{37}
		&(D_\tau \rho^{e,n},\psi)+D_0 (\nabla \rho^{e,n},\nabla \psi ) +\sum_{i=1}^M \nu_iz_i^2((c^{i,n-1}\nabla \phi^{n-1}), \nabla \psi )+ (u^{n-1}\cdot\nabla \rho^{e,n}, \psi )=0, \forall \psi \in W.
	\end{align}
	Subtracting (\ref{37}) with $\psi=\psi_h $ from (\ref{22}) and using the Ritz projection, it yields that
	\begin{align}
		&\left(D_\tau e_{h,\rho}^{n},\psi_h\right)+D_0 (\nabla (e_{h,\rho}^{n}),\nabla \psi_h) +\sum_{i=1}^M \nu_iz_i^2(c^{i,n-1}_h\nabla \phi_h^{n-1},\nabla \psi_h)\nn\\
		&-\sum_{i=1}^M \nu_iz_i^2((c^{i,n-1}\nabla \phi^{n-1}), \nabla \psi )+(u_h^{n-1}\cdot\nabla \rho^{e,n}_h, \psi_h)-(u^{n-1}\cdot\nabla \rho^{e,n}, \psi )\nonumber\\
		&=\left(\frac{(\rho^{e,n}-R_{h,\rho}(\rho^{e,n}))-(\rho^{e,n-1}-R_{h,\rho}(\rho^{e,n-1}))}{\tau},\psi_h\right), \forall \psi_h\in W_h.
	\end{align}
	Taking $\psi_h =2\tau e_{h,\rho}^{n}$ and using $2(a-b,b)=\|a\|_0^2-\|b\|_0^2+\|a-b\|_0^2$, there holds that
	\begin{align*}
		&\|e_{h,\rho}^{n}\|_0^2-\|e_{h,\rho}^{n-1}\|_0^2+\|e_{h,\rho}^{n}-e_{h,\rho}^{n-1}\|_0^2+2\tau D_0 \|\nabla e_{h,\rho}^{n}\|_0^2 \\
		&+2\tau\sum_{i=1}^M  \nu_iz_i^2(c^{i,n-1}_h\nabla \phi_h^{n-1},\nabla e_{h,\rho}^{n})-2\tau \sum_{i=1}^M \nu_iz_i^2(c^{i,n-1}\nabla \phi^{n-1}, \nabla e_{h,\rho}^{n} )\\
		&+2\tau (u_h^{n-1}\cdot\nabla \rho^{e,n}_h, e_{h,\rho}^{n})-2\tau (u^{n-1}\cdot\nabla \rho^{e,n}, e_{h,\rho}^{n} )\\
		=&2\left((\rho^{e,n}-R_{h,\rho}(\rho^{e,n}))-(\rho^{e,n-1}-R_{h,\rho}(\rho^{e,n-1})),e_{h,\rho}^{n}\right).
	\end{align*}
	Adding and subtracting some terms, we get
	\begin{align*}
		&\sum_{i=1}^M \nu_iz_i^2(c^{i,n-1}_h\nabla \phi_h^{n-1},\nabla e_{h,\rho}^{n})-\sum_{i=1}^M \nu_iz_i^2(c^{i,n-1}\nabla \phi^{n-1}, \nabla e_{h,\rho}^{n} )\\
		=&\sum_{i=1}^M \nu_iz_i^2(e_{h,c^i}^{n-1}\nabla \phi_h^{n-1},\nabla e_{h,\rho}^{n})+ \sum_{i=1}^M \nu_iz_i^2((R_{h,c}(c^{i,n-1})-c^{i,n-1})\nabla \phi_h^{n-1},\nabla e_{h,\rho}^{n})\\
		&+\sum_{i=1}^M \nu_iz_i^2(c_h^{i,n-1}\nabla e_{h,\phi}^{n-1},\nabla e_{h,\rho}^{n})+\sum_{i=1}^M \nu_iz_i^2(c_h^{i,n-1}\nabla (R_{h,\phi}(\phi^{e, n-1})-\phi^{e,n-1}),\nabla e_{h,\rho}^{n})
	\end{align*}
	Using Cauchy-Schwarz and Young's inequality, it follows by
	\begin{align*}
		|\sum_{i=1}^M \nu_iz_i^2(e_{h,c^i}^{n-1}\nabla \phi_h^{n-1},\nabla e_{h,\rho}^{n})|\leq C\sum_{i=1}^M\|e_{h,c^i}^{n-1}\|_0^2 +\frac{D_0}{8}\|\nabla e_{h,\rho}^{n}\|_0^2.
	\end{align*}
	Using Cauchy-Schwarz, Young's inequality, and the properties of the Ritz projection, it follows by
	\begin{align*}
		|\sum_{i=1}^M \nu_iz_i^2((R_{h,c}(c^{i,n-1})-c^{i,n-1})\nabla \phi_h^{n-1},\nabla e_{h,\rho}^{n})|\leq & \sum_{i=1}^M C\|R_{h,c}(c^{i,n-1})-c^{i,n-1}\|_0^2+\frac{D_0}{8}\|\nabla e_{h,\rho}^{n}\|_0^2\\
		\leq& Ch^{4}+\frac{D_0}{8}\|\nabla e_{h,\rho}^{n}\|_0^2.
	\end{align*}
	Using Cauchy-Schwarz and Young's inequality, we derive that
	\begin{align*}
		|\sum_{i=1}^M \nu_iz_i^2(c_h^{i,n-1}\nabla e_{h,\phi}^{n-1},\nabla e_{h,\rho}^{n})|\leq C\|e_{h,\phi}^{n-1}\|_0^2+\frac{D_0}{8}\|\nabla e_{h,\rho}^{n}\|_0^2.
	\end{align*}
	Using Cauchy-Schwarz, Young's inequality, and the properties of the Ritz projection, it yields that
	\begin{align*}
		| \sum_{i=1}^M \nu_iz_i^2(c_h^{i,n-1}\nabla (R_{h,\phi}(\phi^{e, n-1})-\phi^{e,n-1}),\nabla e_{h,\rho}^{n})| \leq&  C \| R_{h,\phi}(\phi^{e, n-1})-\phi^{e,n-1}\|_0^2 + \frac{D_0}{8}\|\nabla e_{h,\rho}^{n}\|_0^2\\
		\leq& Ch^{4}+\frac{D_0}{8}\|\nabla e_{h,\rho}^{n}\|_0^2 .
	\end{align*}
	Adding and Subtracting some terms, it holds that
	\begin{align*}
		&(u_h^{n-1}\cdot\nabla \rho^{e,n}_h, e_{h,\rho}^{n})-(u^{n-1}\cdot\nabla \rho^{e,n}, e_{h,\rho}^{n} )\\
		=&(u_h^{n-1}\cdot\nabla e_{h,\rho}^{n}, e_{h,\rho}^{n})+(u_h^{n-1}\cdot\nabla (R_{h,\rho }(\rho^{en})-\rho^{e,n}), e_{h,\rho}^{n})\\
		&+ (e_{u,h}^{n-1}\cdot\nabla \rho^{e,n}, e_{h,\rho}^{n} )+ ((R_h(u^{n-1},p^{n-1}-u^{n-1}))\cdot\nabla\rho^{e,n} , e_{h,\rho}^{n} )
	\end{align*}
	Using Cauchy-Schwarz and Young's inequality, we derive
	\begin{align*}
		|(e_{h,u}^{n-1}\cdot\nabla \rho^{e,n}, e_{h,\rho}^{n})|\leq C\|e_{h,u}^{n-1}\|_0^2+\frac{D_0}{8}\|\nabla e_{h,\rho}^{n}\|_0^2 .
	\end{align*}
	Using Cauchy-Schwarz, Young's inequality, and the properties of the Stokes projection, there holds that
	\begin{align*}
		|((R_h(u^{n-1},p^{n-1})-u^{n-1})\cdot\nabla \rho^{e,n}_h, e_{h,\rho}^{n})| \leq Ch^{4}+\frac{D_0}{8}\|\nabla e_{h,\rho}^{n}\|_0^2 .
	\end{align*}
	Noting $\nabla \cdot u^{n-1}=0$, there holds that
	\begin{align*}
		(u^{n-1}\cdot\nabla e_{h,\rho}^{n}, e_{h,\rho}^{n} )=0.
	\end{align*}
	Using Cauchy-Schwarz, Young's inequality, and the properties of the Stokes projection, it follows by
	\begin{align*}
		|(u^{n-1}\cdot\nabla (R_{h,\rho}(\rho^{e,n})-\rho^{e,n}), e_{h,\rho}^{n} )|\leq Ch^{4}+\frac{D_0}{8}\|\nabla e_{h,\rho}^{n}\|_0^2 .
	\end{align*}
	Using Cauchy-Schwarz, Young's inequality and the  properties of the Ritz projection, we derive that
	\begin{align*}
		&|2\left((\rho^{e,n}-R_{h,\rho}(\rho^{e,n}))-(\rho^{e,n-1}-R_{h,\rho}(\rho^{e,n-1})),e_{h,\rho}^{n}\right)|\\
		 =& |2\left((\rho^{e,n}-\rho^{e,n-1})-R_{h,\rho}(\rho^{e,n}-\rho^{e,n-1}),e_{h,\rho}^{n}\right)|\\
		\leq& C\|(\rho^{e,n}-\rho^{e,n-1})-R_{h,\rho}(\rho^{e,n}-\rho^{e,n-1})\|_0^2\|\nabla  e_{h,\rho}^{n}\|_0^2\\
		\leq& C\tau h^{2}\| D_\tau \rho^{e,n}\|_2\|\nabla e_{h,\rho}^{n}\|_0\\
		\leq& C\tau h^{4}+\frac{\tau D_0}{8}\|\nabla e_{h,\rho}^{n}\|_0^2 .
	\end{align*}
	Then, we arrive at
\begin{align}
	&\|e_{h,\rho}^{n}\|_0^2-\|e_{h,\rho}^{n-1}\|_0^2+\|e_{h,\rho}^{n}-e_{h,\rho}^{n-1}\|_0^2+D_0 \tau \|\nabla e_{h,\rho}^{n}\|_0^2\nonumber\\
	\leq&  C\tau \sum_{i=1}^M\|e_{h,c^{i}}^{n-1}\|_0^2+ C\tau \|e_{h,\phi}^{n-1}\|_0^2+ C\tau \|e_{h,u}^{n-1}\|_0^2+C\tau h^{4}.\label{er_rho}
\end{align}

	Using Green's formula, we deduce the weak form of (\ref{step21})
	\begin{align*}
		&\left(\frac{u^n-u^{n-1}}{\tau}, v\right)
	+\mu \left(\nabla u^n,\nabla v \right)+\left(u^{n-1}\cdot\nabla u^n,v \right)- \left(p^{n}, \nabla \cdot v \right)-\left(\rho^{e,n} \nabla \phi^{n},v\right)=0,\\
	&\left(\nabla \cdot u^n,q \right)=0,\ \forall v\in V, q\in M.
		\end{align*}
	Taking $v=v_h$, subtracting it from (\ref{24}) and using the definition of the Stokes projection, we get the error equation for $u$ as follows
	\begin{align*}
		&\left(\frac{e_{h,u}^n-e_{h,u}^{n-1}}{\tau
		},v_h\right)+\mu(\nabla e_{h,u}^n ,\nabla v_h)+(u_h^{n-1}\cdot\nabla u_h^n,v_h)-\left(u^{n-1}\cdot\nabla u^n,v_h \right)- (\nabla\cdot v_h,e_{h,p}^{n})\nonumber\\
		&- (\rho_h^{e,n-1}\nabla \phi_h^n,v_h)+\left(\rho^{e,n} \nabla \phi^{n},v_h\right) \\
		=&\left(\frac{(u^n-R_h(u^n,p^n))-(u^{n-1}-R_h(u^{n-1},p^{n-1}))}{\tau
		},v_h\right), \forall v_h\in X_h,\\
		&(\nabla \cdot u_h^n,q_h)=0,\forall q_h\in M_h.
	\end{align*}
Letting $v_h=2\tau e_{h,u}^n$ and using $2(a-b,a)=\|a\|_0^2-\|b\|_0^2+\|a-b\|_0^2$, we derive that
	\begin{align*}
		&\|e_{h,u}^n\|_0^2-\|e_{h,u}^{n-1}\|_0^2+\|e_{h,u}^n-e_{h,u}^{n-1}\|_0^2
		+2\tau \mu\|\nabla e_{h,u}^n \|_0^2
		+2\tau (u_h^{n-1}\cdot\nabla u_h^n,e_{h,u}^n)\nonumber\\
		&-2\tau \left(u^{n-1}\cdot\nabla u^n,e_{h,u}^n \right)-2\tau  (\rho_h^{e,n}\nabla \phi_h^n,e_{h,u}^n)+2\tau \left(\rho^{e,n} \nabla \phi^{n},e_{h,u}^n\right) \\
		=&2 \tau \left(\frac{(u^n-R_h(u^n,p^n))-(u^{n-1}-R_h(u^{n-1},p^{n-1}))}{\tau
		},e_{h,u}^n\right).
	\end{align*}
	Adding and subtracting some terms, there holds that
	\begin{align*}
		&2\tau (u_h^{n-1}\cdot\nabla u_h^n,e_{h,u}^n)-2\tau \left(u^{n-1}\cdot\nabla u^n,e_{h,u}^n \right)\\
		=& 2\tau (u_h^{n-1}\cdot\nabla (u_h^n-u^n),e_{h,u}^n)+2\tau \left((u_h^{n-1}-u^{n-1})\cdot\nabla u^n,e_{h,u}^n \right).
	\end{align*}
	Using $\|\nabla u_h^{n-1}\|_\infty <+\infty$, we have
	\begin{align*}
		|(u_h^{n-1}\cdot\nabla (u_h^n-u^n),e_{h,u}^n)|\leq & |(u_h^{n-1}\cdot\nabla e_{u,h}^n,e_{h,u}^n)|+|(u_h^{n-1}\cdot\nabla (u^n-R_h(u^n,p^n)),e_{h,u}^n)|\\
		\leq & C\tau h^{4}+\frac{\mu \tau }{8}\|\nabla e_{h,u}^n\|_0^2 .
	\end{align*}
	Using Cauchy-Schwarz, Young's inequality, and the properties of Stokes projection, we deduce that
	\begin{align*}
		| ((u_h^{n-1}-u^{n-1})\cdot\nabla u^n,e_{h,u}^n)| \leq & | (e_{h,u}^{n-1}\cdot\nabla u^n,e_{h,u}^n)|+| ((R_h(u^{n-1},p^{n-1})-u^{n-1})\cdot\nabla u^n,e_{h,u}^n)|\\
		\leq& C\tau h^{4}+C\tau\|e_{h,u}^{n-1}\|_0^2 +\frac{\mu \tau }{8}\|\nabla e_{h,u}^n\|_0^2 .
	\end{align*}
	Adding and subtracting some terms, it yields that
	\begin{align*}
	& 2\tau  (\rho_h^{e,n}\nabla \phi_h^n,e_{h,u}^n)-2\tau \left(\rho^{e,n} \nabla \phi^{n},e_{h,u}^n\right)\\
	=&	2\tau  (e_{h,\rho}^{n}\nabla \phi_h^n,e_{h,u}^n)+
	2\tau  ((R_{h,\rho}(\rho^{e,n})-\rho^{e,n})\nabla \phi_h^n,e_{h,u}^n)\\
	&+2\tau \left(\rho^{e,n} \nabla e_{h,\phi}^n,e_{h,u}^n\right)
	+2\tau \left(\rho^{e,n} \nabla (R_{h,\phi}(\phi ^n)-\phi^{n}),e_{h,u}^n\right).
	\end{align*}
	Using (\ref{23}), we can deduce that $\|\phi_h^n\|_2\leq \|\nabla\rho_{h}^{e,n}\|_0$. Then, we have
	\begin{align*}
		2\tau |(e_{h,\rho}^{n}\nabla \phi_h^n,e_{h,u}^n)|\leq & C\tau  \|e_{h,\rho}^{n}\|_0\|\phi_h^n\|_2\|\nabla e_{h,u}^n\|_0\\
		\leq & C\tau  \|e_{h,\rho}^{n}\|_0\|\nabla\rho_h^n\|_0\|\nabla e_{h,u}^n\|_0\\
		\leq &C\tau  \|e_{h,\rho}^{n}\|_0^2+\frac{\mu \tau }{8}\|\nabla e_{h,u}^n\|_0^2.
	\end{align*}
	Using Cauchy-Schwarz, Young's inequality, and the properties of Ritz projection, it holds that
	\begin{align*}
		2\tau | ((R_{h,\rho}(\rho^{e,n})-\rho^{e,n})\nabla \phi_h^n,e_{h,u}^n)|\leq& C\tau \| R_{h,\rho}(\rho^{e,n})-\rho^{e,n}\|_0 \| \phi_h^n\|_2\|\nabla e_{h,u}^n\|_0\\
		\leq & C\tau \| R_{h,\rho}(\rho^{e,n})-\rho^{e,n}\|_0 \|\nabla\rho_h^n\|_0 \|\nabla e_{h,u}^n\|_0\\
		\leq & C\tau h^4+\frac{\nu \tau}{8} \|\nabla e_{h,u}^n\|_0^2.
	\end{align*}
	Taking $\psi=\psi_h$ in (\ref{10}) and subtracting it  from (\ref{23}), we have
	\begin{align}\label{3.58}
		(\nabla (\phi_h^{n}-\phi^{n}),\nabla \psi_h)=(\rho_h^{e,n}-\rho^{e,n},\psi_h).
	\end{align}
	Taking $\psi_h=e_{h,\psi}^n$, it yields that
	\begin{align*}
		\|\nabla e_{h,\phi}^{e,n}\|_0^2 \leq Ch^4+C\|e_{h,\rho}^{n}\|_0^2.
	\end{align*}
	Then, there holds that
	\begin{align*}
		2\tau |\left(\rho^{e,n} \nabla e_{h,\phi}^n,e_{h,u}^n\right)|\leq& C\tau \|\nabla \rho^{e,n}\|_0\| \nabla e_{h,\phi}^n\|_0\|\nabla e_{h,u}^n\|_0\\
		\leq & C\tau h^4+C\tau \|e_{h,\rho}^{n}\|_0^2+\frac{\nu \tau}{8} \|\nabla e_{h,u}^n\|_0^2.
	\end{align*}
	Using Cauchy-Schwarz, Young's inequality, and the properties of Stokes projection, we have
	\begin{align*}
		2\tau |\left(\rho^{e,n} \nabla (R_{h,\phi}(\phi ^n)-\phi^{n}),e_{h,u}^n\right)|\leq & C\tau \|\nabla \rho^{e,n}\|_0\|R_{h,\phi}(\phi ^n)-\phi^{n}\|_0\|\nabla e_{h,u}^n\|_0\\
		\leq &C\tau h^4+\frac{\nu \tau}{8} \|\nabla e_{h,u}^n\|_0^2,
	\end{align*}
	and 	
	\begin{align*}
		2 \tau \left|\left(\frac{(u^n-R_h(u^n,p^n))-(u^{n-1}-R_h(u^{n-1},p^{n-1}))}{\tau
		},e_{h,u}^n\right)\right|&\leq C\tau h^{2}(\| D_\tau u^{n}\|_2+\|D_\tau p^n\|_0^2)\|\nabla e_{h,u}^{n}\|_0\\
		&\leq C\tau h^{4}+\frac{\tau D}{8}\|\nabla e_{h,u}^{n}\|_0^2 .
	\end{align*}
Then, we arrive at
\begin{align}
	&\|e_{h,u}^n\|_0^2-\|e_{h,u}^{n-1}\|_0^2+\|e_{h,u}^n-e_{h,u}^{n-1}\|_0^2
	+2\tau \nu\|\nabla e_{h,u}^n \|_0^2\nonumber\\
	\leq & C\tau h^4+C\tau \|e_{h,\rho}^{n}\|_0^2+C\tau\|e_{h,u}^{n-1}\|_0^2.\label{er_u}
\end{align}

We can take the weak form of (\ref{12}) as follows
\begin{align*}
	\left(D_\tau c^{i,n},\zeta \right)+d_i(\nabla c^{i,n},\nabla \zeta) +(u^{n}\cdot\nabla c^{i,n},\zeta)&+ \nu_i z_i ( c^{i,n}\nabla\phi^{n-1} ,\nabla \zeta)=0,\\ &i=1,\ldots,M, \forall \zeta\in W.
\end{align*}
Subtracting it with $\zeta=\zeta_h$ from (\ref{25}), it follows that
	\begin{align*}
		&\left(D_\tau e_{h,c^i}^{n},\zeta_h \right)+d_i(\nabla e_{h,c^i}^{n},\nabla \zeta_h) +(u_h^{n}\cdot\nabla c_h^{i,n},\zeta_h)-(u^{n}\cdot\nabla c^{i,n},\zeta_h)+ \nu_i z_i ( c_h^{i,n}\nabla\phi_h^{n-1} ,\nabla \zeta_h)\\
		&-\nu_i z_i ( c^{i,n}\nabla\phi^{n-1} ,\nabla \zeta_h)=\left(D_\tau (R_{h,c^i}(c^{i,n})-c^{i,n}),\zeta_h \right), i=1,\ldots,M, \forall \zeta_h\in W_h.
	\end{align*}
Taking $\zeta_h=2\tau e_{h,c^i}^{n}$ and using $2(a-b,b)=\|a\|_0^2-\|b\|_0^2+\|a-b\|_0^2$, it can be deduced that
\begin{align*}
	&\|e_{h,c^i}^{n}\|_0^2-\|e_{h,c^i}^{n-1}\|_0^2+ \|e_{h,c^i}^{n}-e_{h,c^i}^{n-1}\|_0^2+2\tau d_i\|\nabla e_{h,c^i}^{n}\|_0^2 +2\tau (u_h^{n}\cdot\nabla c_h^{i,n},e_{h,c^i}^{n})\\
	&-2\tau (u^{n}\cdot\nabla c^{i,n},e_{h,c^i}^{n})+ 2\tau \nu_i z_i ( c_h^{i,n}\nabla\phi_h^{n-1} ,\nabla e_{h,c^i}^{n})-2\tau \nu_i z_i ( c^{i,n}\nabla\phi^{n-1} ,\nabla e_{h,c^i}^{n}) \\
	&=2\tau \left(D_\tau (R_{h,c^i}(c^{i,n})-c^{i,n}),e_{h,c^i}^{n} \right),&i=1,\ldots,M.
\end{align*}
Adding and subtracting some terms, we have
\begin{align*}
	&2\tau (u_h^{n}\cdot\nabla c_h^{i,n},e_{h,c^i}^{n})-2\tau (u^{n}\cdot\nabla c^{i,n},e_{h,c^i}^{n}) \\
	=& 2\tau (u_h^{n}\cdot\nabla (c_h^{i,n}-c^{i,n}),e_{h,c^i}^{n})+2\tau ((u_h^{n}-u^n)\cdot\nabla c^{i,n},e_{h,c^i}^{n})\\
	=& 2\tau (u_h^{n}\cdot\nabla e_{h,c^i}^{n},e_{h,c^i}^{n})+
	2\tau (e_{h,u}^{n}\cdot\nabla (R_{h,c}(c^{i,n})-c^{i,n}),e_{h,c^i}^{n})+2\tau (e_{h,u}^n\cdot\nabla c^{i,n},e_{h,c^i}^{n})\\
	&+
	2\tau ((R_h(u^n,p^n)-u^n)\cdot\nabla c^{i,n},e_{h,c^i}^{n})+
	2\tau (R_h(u^n,p^n)\cdot\nabla (R_{h,c}(c^{i,n})-c^{i,n}),e_{h,c^i}^{n}).
\end{align*}
Noting $(\nabla\cdot u_h^n, q_h)=0, \forall q_h\in M_h$, there holds that
\begin{align*}
	2\tau (u_h^{n}\cdot\nabla e_{h,c^i}^{n},e_{h,c^i}^{n})=0.
\end{align*}
Using Cauchy-Schwarz and Young's inequality, we deduce that
\begin{align*}
	2\tau |(e_{h,u}^{n}\cdot\nabla (R_{h,c}(c^{i,n})-c^{i,n}),e_{h,c^i}^{n})| \leq& C\tau \| e_{h,u}^{n}\|_0 (\|R_{h,c}(c^{i,n})\|_2+\|c^{i,n}\|_2)\|\nabla
	e_{h,c^i}^{n}\|_0\\
\leq& C\tau \| e_{h,u}^{n}\|_0^2+ \frac{d_i\tau}{8}\|\nabla
e_{h,c^i}^{n}\|_0^2.
\end{align*}
Using Cauchy-Schwarz and Young's inequality, it follows that
\begin{align*}
	2\tau |(e_{h,u}^n\cdot\nabla c^{i,n},e_{h,c^i}^{n})|\leq & C\tau \|e_{h,u}^n\|_0\|\nabla c^{i,n}\|_\infty \|\nabla e_{h,c^i}^{n}\|_0\\
	\leq & C\tau \| e_{h,u}^{n}\|_0^2+ \frac{d_i\tau}{8}\|\nabla
	e_{h,c^i}^{n}\|_0^2.
\end{align*}
Using Cauchy-Schwarz, Young's inequality and the properties of Stokes projection, we derive
\begin{align*}
	2\tau |((R_h(u^n,p^n)-u^n)\cdot\nabla c^{i,n},e_{h,c^i}^{n})|\leq & C\tau \|R_h(u^n,p^n)-u^n\|_0\| c^{i,n}\|_2\| \nabla e_{h,c^i}^{n}\|_0\\
	\leq & C\tau h^4+ \frac{d_i\tau}{8}\|\nabla
	e_{h,c^i}^{n}\|_0^2.
\end{align*}
Adding and Subtracting some terms, it yields that
\begin{align*}
	&2\tau \nu_i z_i ( c_h^{i,n}\nabla\phi_h^{n-1} ,\nabla e_{h,c^i}^{n})-2\tau \nu_i z_i ( c^{i,n}\nabla\phi^{n-1} ,\nabla e_{h,c^i}^{n})\\
	=& 2\tau \nu_i z_i ( (c_h^{i,n}-c^{i,n})\nabla\phi_h^{n-1} ,\nabla e_{h,c^i}^{n})+2\tau \nu_i z_i ( c^{i,n}\nabla(\phi_h^{n-1}-\phi^{n-1} ),\nabla e_{h,c^i}^{n})\\
	=& 2\tau \nu_i z_i ( e_{h,c^i}^{n}\nabla\phi_h^{n-1} ,\nabla e_{h,c^i}^{n})+2\tau \nu_i z_i ( (R_{h,c}(c^{i,n})-c^{i,n})\nabla\phi_h^{n-1} ,\nabla e_{h,c^i}^{n})\\
	&+2\tau \nu_i z_i ( c^{i,n}\nabla e_{h,\phi}^{n-1},\nabla e_{h,c^i}^{n})
	+2\tau \nu_i z_i ( c^{i,n}\nabla(R_{h,\phi}(\phi^{n-1})-\phi^{n-1} ),\nabla e_{h,c^i}^{n}).
\end{align*}
Using Cauchy-Schwarz and Young's inequality, it follows by
\begin{align*}
	2\tau \nu_i z_i |( e_{h,c^i}^{n}\nabla\phi_h^{n-1} ,\nabla e_{h,c^i}^{n})|\leq& C\tau \| e_{h,c^i}^{n}\|_0\|\nabla\phi_h^{n-1}\|_\infty \|\nabla e_{h,c^i}^{n}\|\\
	\leq & C\tau \| e_{h,c^i}^{n}\|_0^2+\frac{d_i\tau}{8}\|\nabla
	e_{h,c^i}^{n}\|_0^2.
\end{align*}
Using Cauchy-Schwarz, Young's inequality and the Ritz projection, we deduce that
\begin{align*}
	2\tau \nu_i z_i |( (R_{h,c}(c^{i,n})-c^{i,n})\nabla\phi_h^{n-1} ,\nabla e_{h,c^i}^{n})|\leq & C\tau \|R_{h,c}(c^{i,n})-c^{i,n}\|_0\|\nabla\phi_h^{n-1}\|_{\infty}
	\|\nabla e_{h,c^i}^{n}\|_0\\
	\leq & C\tau h^4+\frac{d_i\tau}{8}\|\nabla
	e_{h,c^i}^{n}\|_0^2.
\end{align*}
Using (\ref{3.58}) and the properties of Ritz projection, we have
\begin{align*}
	(\nabla e_{h, \phi}^{n-1},\nabla \psi_h)=(\rho_h^{e,n-1}-\rho^{e,n-1},\psi_h).
\end{align*}
Taking $\psi_h= e_{h,\phi}^{n-1}$, we have
\begin{align*}
	\|\nabla e_{h,\phi}^{n-1}\|_0\leq Ch^2+\|e_{h,\rho}^{n-1}\|_0.
\end{align*}
Then, we can deduce
\begin{align*}
	2\tau \nu_i z_i |( c^{i,n}\nabla e_{h,\phi}^{n-1},\nabla e_{h,c^i}^{n})|\leq &C\tau \|c^{i,n}\|_\infty \|\nabla e_{h,\phi}^{n-1}\|_0\|\nabla e_{h,c^i}^{n}\|_0\\
	\leq  & C\tau h^4+C\tau \|e_{h,\rho}^{n-1}\|_0^2+\frac{d_i\tau}{8}\|\nabla
	e_{h,c^i}^{n}\|_0^2.
\end{align*}
Using Cauchy-Schwarz, Young's inequality and the Ritz projection, there holds that
\begin{align*}
	2\tau \nu_i z_i |( c^{i,n}\nabla(R_{h,\phi}(\phi^{n-1})-\phi^{n-1} ),\nabla e_{h,c^i}^{n})|\leq&  C\tau \|c^{i,n}\|_2\|R_{h,\phi}(\phi^{n-1})-\phi^{n-1}\|_0\| \nabla e_{h,c^i}^{n}\|_0\\
	\leq  & C\tau h^4+\frac{d_i\tau}{8}\|\nabla
	e_{h,c^i}^{n}\|_0^2,
\end{align*}
and
\begin{align*}
	2\tau |\left(D_\tau (R_{h,c^i}(c^{i,n})-c^{i,n}),e_{h,c^i}^{n} \right)| &\leq C\tau h^{2}\| D_\tau c^{i,n}\|_2\|\nabla e_{h,c}^{n}\|_0\\
	&\leq C\tau h^{4}+\frac{d_i \tau}{8}\|\nabla e_{h,c^i}^{n}\|_0^2 .
\end{align*}
Then, we arrive at
\begin{align*}
	&\|e_{h,c^i}^{n}\|_0^2-\|e_{h,c^i}^{n-1}\|_0^2+ \|e_{h,c^i}^{n}-e_{h,c^i}^{n-1}\|_0^2+\tau d_i\|\nabla e_{h,c^i}^{n}\|_0^2 \leq C\tau \| e_{h,c^i}^{n}\|_0^2+C\tau \| e_{h,u}^{n}\|_0^2+ C\tau h^{4}.
\end{align*}
Taking sum of it over all $i$, it yields that
\begin{align}\label{er_c}
	&\sum_{i=1}^M\|e_{h,c^i}^{n}\|_0^2-\sum_{i=1}^M\|e_{h,c^i}^{n-1}\|_0^2+\tau d_i\sum_{i=1}^M\|\nabla e_{h,c^i}^{n}\|_0^2 \leq C\tau \sum_{i=1}^M\| e_{h,c^i}^{n}\|_0^2+C\tau \| e_{h,u}^{n}\|_0^2+ C\tau h^{4}.
\end{align}

Combining (\ref{er_rho}), (\ref{er_u}) and (\ref{er_c}), we get
\begin{align*}
	&\|e_{h,\rho}^{n}\|_0^2-\|e_{h,\rho}^{n-1}\|_0^2+D_0\tau \|\nabla e_{h,\rho}^{n}\|_0^2  + \|e_{h,u}^n\|_0^2-\|e_{h,u}^{n-1}\|_0^2
	+2\tau \mu\|\nabla e_{h,u}^n \|_0^2\\
	&+\sum_{i=1}^M\|e_{h,c^i}^{n}\|_0^2-\sum_{i=1}^M\|e_{h,c^i}^{n-1}\|_0^2+\tau d_i\sum_{i=1}^M\|\nabla e_{h,c^i}^{n}\|_0^2\\
	\leq&  C\tau \sum_{i=1}^M\|e_{h,c^{i}}^{n-1}\|_0^2+ C\tau \|e_{h,\phi}^{n-1}\|_0^2+ C\tau \|e_{h,u}^{n}\|_0^2+C\tau \|e_{h,\rho}^{n}\|_0^2+C\tau\|e_{h,u}^{n-1}\|_0^2+C\tau h^{4}.
\end{align*}
Taking it over all $n$ and using Gronwall's lemma when $\tau$ is sufficient small, we get
\begin{align}\label{error}
	&\|e_{h,\rho}^{n}\|_0^2+D_0 \tau \sum_{n=1}^N\|\nabla e_{h,\rho}^{n}\|_0^2
	+ \|e_{h,u}^N\|_0^2+2\tau \sum_{n=1}^N \mu\|\nabla e_{h,u}^n \|_0^2\nonumber\\
	&+\sum_{i=1}^M\|e_{h,c^i}^{n}\|_0^2+\tau d_i\sum_{n=1}^N\sum_{i=1}^M\|\nabla e_{h,c^i}^{n}\|_0^2
	\leq  C h^{4}.
\end{align}

Secondly from \refe{error}, we see that
\begin{align*}
\max_{0\leq m\leq N}\| u_h^{m} \|_{L^\infty}
&\le
\max_{0\leq m\leq N}(\| R_h^{m} \|_{L^\infty} + \| e_h^{m} \|_{L^\infty})\\
&\le
C \max_{0\leq m\leq N}\|U^m\|_2+C\max_{1\leq m\leq N}\|P^m\|_1
+ C h^{-d/2} \max_{0\leq m\leq N}\| e_h^{m} \|_0\\
&\leq C\, ,
\end{align*}and
\begin{align*}
\tau \sum_{m=0}^N\| u_h^{m} \|_{W^{1,\infty}}^2
&\le 2\tau \sum_{m=0}^N(\| R_h^{m} \|_{W^{1,\infty}}^2
+ \| e_h^{m} \|_{W^{1,\infty}}^2)\\
&\le C\tau\left(\sum_{m=0}^N\|U^m\|_{W^{1,\infty}}^2
+\sum_{m=1}^N\|P^m\|_{L^{\infty}}^2\right)
+ C \tau h^{-d} \sum_{m=0}^N\| e_h^{m} \|_1 ^2\\
&\leq C\, .
\end{align*}
Similarly, we have
\begin{align*}
	\max_{0\leq m\leq N}\| \rho_h^{e,m} \|_{L^\infty}+\tau \sum_{m=0}^N\| \rho_h^{e,m} \|_{W^{1,\infty}}^2 &\le C,\\
	\max_{0\leq m\leq N}\| c_h^{i,m} \|_{L^\infty}+\tau \sum_{m=0}^N\| c_h^{i,m} \|_{W^{1,\infty}}^2 &\le C, i=1, \ldots, M.
\end{align*}

\end{proof}

\subsection{The optimal error estimate for the finite element method}
To give the optimal error estimation, we define the errors as follows
\begin{align*}
	\e_{h,\rho}^n=\rho_h^{e,n}-R_{h,\rho}(\rho^{e}(t_n)); &~~ \e_{h,\phi}^n=\phi_h^n-R_{h,\phi}(\phi(t_n));& \e_{h,c^i}^{n}= c_h^{i,n}-R_{h,c} (c^{i}(t_n));\\
	\e_{h,u}^n=u_h^n-R_h(u(t_n),p(t_n));&~~\e_{h,p}^n=p_h^n-Q_h(u(t_n), p(t_n)).
\end{align*}
\begin{theorem}
	Suppose $c_h^{i,n}\in W_h$, $\phi_h^n\in W_h$, $u_h^n\in V$ and $p_h^n\in Q_h$ be the solutions of (\ref{AL2.1}), when $\tau$ is sufficient small, we get the error estimates as follows
	\begin{align*}
		\|\rho^e(T)-\rho_h^{e,N}\|_0^2&+D_0 \tau h\sum_{n=1}^N\|\nabla (\rho^e(t_n)-\rho_h^{e,n})\|_0^2 +\|u(T)-u_h^N\|_0^2 + \tau h\mu\sum_{n=1}^N\|\nabla (u(t_n)-u_h^n)\|_0^2\nonumber\\
		&+\sum_{i=1}^M\|c^i(T)-c_h^{i,N}\|_0^2+ \tau h\sum_{n=1}^N \sum_{i=1}^M d_i\|\nabla c^i(t_n)-c_h^{i,n}\|_0^2
		\leq C(\tau^2+h^{2r+2}).
	\end{align*}
	Furthermore, we have
	\begin{eqnarray*}
		\sum_{n=1}^{N} \tau \| p(t_n)-p_h^n\|_0^2\leq C(\tau^2+h^{2r}).
	\end{eqnarray*}	
\end{theorem}
\begin{proof}
We deduce the weak form of (\refe{requ}) at $t=t_n$, as follows
\begin{align}
	&(D_\tau\rho^e(t_n),\psi)+D_0 (\nabla \rho^e(t_n),\nabla \psi) +\sum_{i=1}^M \nu_iz_i((c^i(t_n)\nabla \phi(t_n)),\nabla \psi)\nonumber\\
	&+(u(t_n)\cdot\nabla \rho^e(t_n),\psi)=(Tr_\rho^n,\psi),&\forall \psi\in W, \label{3.66a}\\
	&(D_\tau u(t_n),v)+\mu (\nabla u(t_n), \nabla v)+((u(t_n)\cdot\nabla)u(t_n),v) \nonumber\\
	&- (p(t_n), \nabla \cdot v) -(\rho^e(t_n)\nabla \phi(t_n),v)=(Tr_u^n,v),&\forall v\in X, \label{3.66b}\\
	&(\nabla\cdot u(t_n),q)=0,&\forall q\in M, \label{3.66c}\\
	&(D_\tau c^i_t(t_n),\zeta)+ d_i(\nabla c^i , \nabla \zeta)+\nu_i z_i (c^i(t_n)\nabla\phi(t_h), \nabla \zeta )\nonumber\\
	&+(u\cdot\nabla c^i,\zeta)=(Tr_{c^i}^n, \zeta ),i=1,\ldots,M,&\forall \zeta \in W, \label{3.66d}\\
	&\varepsilon (\nabla \phi(t_n),\nabla \theta )-(\rho^e(t_n),\theta)=0,& \forall \theta \in W. \label{3.66e}
	\end{align}

Taking $\psi=\psi_h$ in (\ref{3.66a}) and subtracting it form (\ref{22}), we get error equation of $\rho_{h}^{e,n}$ between $\rho^e(t_n)$ as follows
\begin{align}
	&\left(D_\tau\e_{h,\rho}^n,\psi_h\right)+D_0 (\nabla \e_{h,\rho}^n,\nabla \psi_h)+\sum_{i=1}^M \nu_iz_i^2 (c^{i,n-1}_h\nabla \phi_h^{n-1},\nabla \psi_h)-\sum_{i=1}^M \nu_iz_i((c^i(t_n)\nabla \phi(t_n)),\nabla \psi_h)\nn\\
	&+(u_h^{n-1}\cdot\nabla \rho^{e,n}_h, \psi_h)-(u(t_n)\cdot\nabla \rho^e(t_n),\psi_h) =(Tr_\rho^n,\psi_h), \forall \psi_h\in W_h.
\end{align}

Taking $\psi_h =2\tau \e_{h,\rho}^{n}$ and using $2(a-b,b)=\|a\|_0^2-\|b\|_0^2+\|a-b\|_0^2$, there holds that
	\begin{align*}
		&\|\e_{h,\rho}^{n}\|_0^2-\|\e_{h,\rho}^{n-1}\|_0^2+\|\e_{h,\rho}^{n}-\e_{h,\rho}^{n-1}\|_0^2+2\tau D_0 \|\nabla \e_{h,\rho}^{n}\|_0^2 \\
		&+2\tau\sum_{i=1}^M  \nu_iz_i^2(c^{i,n-1}_h\nabla \phi_h^{n-1},\nabla \e_{h,\rho}^{n})-2\tau \sum_{i=1}^M \nu_iz_i^2(c^{i}(t_n)\nabla \phi(t_n), \nabla \e_{h,\rho}^{n} )\nn\\
		&+2\tau (u_h^{n-1}\cdot\nabla \rho^{e,n}_h, \e_{h,\rho}^{n})-2\tau (u(t_n)\cdot\nabla \rho^{e}(t_n), \e_{h,\rho}^{n} )\\
		=&2\tau \left(D_\tau (\rho^{e,n}-R_{h,\rho}(\rho^{e}(t_n))),\e_{h,\rho}^{n}\right)+2\tau (Tr_\rho^n,\e_{h,\rho}^{n}).
	\end{align*}
	Adding and subtracting some terms and using Cauchy-Schwarz, Young's inequality, and Taylor's formula, it yields that
	\begin{align*}
		&2\tau |\sum_{i=1}^M \nu_iz_i^2(c^{i,n-1}_h\nabla \phi_h^{n-1},\nabla \e_{h,\rho}^{n})-\sum_{i=1}^M \nu_i z_i^2(c^{i}(t_n)\nabla \phi(t_n), \nabla \e_{h,\rho}^{n} )|\\
		=&2\tau |\sum_{i=1}^M \nu_iz_i^2(\e_{h,c^i}^{n-1}\nabla \phi_h^{n-1},\nabla \e_{h,\rho}^{n})+ \sum_{i=1}^M \nu_iz_i^2((R_{h,c}(c^{i}(t_{n-1}))-c^{i}(t_{n-1}))\nabla \phi_h^{n-1},\nabla \e_{h,\rho}^{n})|\\
		&+2\tau |\sum_{i=1}^M \nu_iz_i^2((c^{i}(t_n)-c^{i}(t_{n-1}))\nabla \phi_h^{n-1},\nabla \e_{h,\rho}^{n})+\sum_{i=1}^M \nu_iz_i^2(c_h^{i}(t_n)\nabla \e_{h,\phi}^{n-1},\nabla \e_{h,\rho}^{n})|\\
		&+2\tau |\sum_{i=1}^M \nu_iz_i^2(c^{i}(t_n)\nabla (R_{h,\phi}(\phi^{e}(t_{n-1}))-\phi(t_{n-1})),\nabla \e_{h,\rho}^{n})+\sum_{i=1}^M \nu_iz_i^2(c^{i}(t_n)\nabla (\phi(t_{n})-\phi(t_{n-1})),\nabla \e_{h,\rho}^{n})|\\
		\leq & C\tau  \sum_{i=1}^M \|\e_{h,c^i}^{n-1}\|_0\|\nabla \phi_h^{n-1}\|_{\infty}\|\nabla \e_{h,\rho}^{n}\|_0+ C\tau\sum_{i=1}^M  \|R_{h,c}(c^{i}(t_{n-1}))-c^{i}(t_{n-1})\|_0\|\nabla \phi_h^{n-1}\|_{\infty}\|\nabla \e_{h,\rho}^{n}\|_0 \\
		& +  C\tau \sum_{i=1}^M \|c^{i}(t_n)-c^{i}(t_{n-1})\|_0\|\nabla \phi_h^{n-1}\|_{\infty}\|\nabla \e_{h,\rho}^{n}\|_0 +C \tau \sum_{i=1}^M \|c_h^{i}(t_n)\|_{W^{1,\infty}}\|  \e_{h,\phi}^{n-1}\|_0\|\nabla \e_{h,\rho}^{n}\|_0\\
		&+ C\tau\sum_{i=1}^M  \|c^{i}(t_n)\|_{W^{1,\infty}} \|R_{h,\phi}(\phi^{e}(t_{n-1}))-\phi(t_{n-1})\|_0\|\nabla \e_{h,\rho}^{n}\|_0 \\
		&+C\tau \sum_{i=1}^M  \|c^{i}(t_n)\|_{W^{1,\infty}} \|\phi(t_{n})-\phi(t_{n-1})\|_0\|\nabla \e_{h,\rho}^{n}\|_0\\
		\leq & C\tau \sum_{i=1}^M  \|\e_{h,c^i}^{n-1}\|_0^2+C \tau \| \e_{h,\rho}^{n-1}\|_0^2 + C\tau h^{r+1}+C\tau^3+\frac{D_0 \tau}{4}\|\nabla \e_{h,\rho}^{n}\|_0^2.
	\end{align*}
	Adding and subtracting some terms and using Cauchy-Schwarz, Young's inequality, and Taylor's formula, there holds that
	\begin{align*}
		&|2\tau (u_h^{n-1}\cdot\nabla \rho^{e,n}_h, \e_{h,\rho}^{n})-2\tau (u(t_n)\cdot\nabla \rho^{e}(t_n), \e_{h,\rho}^{n} )|\\
		\leq & 2\tau |(\e_{h,u}^{n-1}\cdot \nabla \rho^{e,n}_h, \e_{h,\rho}^{n}) |+2\tau |((R_h(u(t_{n-1}),p(t_{n-1}))-u(t_{n-1}))\cdot \nabla \rho^{e,n}_h, \e_{h,\rho}^{n}) |\\
		&+ 2\tau |(u(t_{n-1})-u(t_{n})\cdot \nabla \rho^{e,n}_h, \e_{h,\rho}^{n})|
		+2\tau |(u(t_n)\cdot\nabla \e_{h,\rho}^{n -1}, \e_{h,\rho}^{n} ) |\\
		&+2\tau |(u(t_n)\cdot\nabla  (R_{h,\rho}(\rho(t_{n-1}))-\rho(t_{n-1})) , \e_{h,\rho}^{n} ) |+2\tau |(u(t_n)\cdot\nabla  (\rho(t_{n-1})-\rho(t_{n})) , \e_{h,\rho}^{n} ) |\\
		\leq & C\tau \|\e_{h,u}^{n-1}\|_0\|\rho^{e,n}_h\|_{W^{1,\infty}} \|\nabla \e_{h,\rho}^{n}\|_0 + C\tau \|R_h(u(t_{n-1}),p(t_{n-1}))-u(t_{n-1})\|_0\|\rho^{e,n}_h\|_{W^{1,\infty}} \|\nabla \e_{h,\rho}^{n}\|_0\\
		&+C\tau \|u(t_{n-1})-u(t_{n})\|_0\|\rho^{e,n}_h\|_{W^{1,\infty}} \|\nabla \e_{h,\rho}^{n}\|_0+C\tau \|u(t_n)\|_{W^{1,\infty}}\|\e_{h,\rho}^{n -1}\|_0\|\nabla  \e_{h,\rho}^{n}\|_0\\
		&+C\tau \|u(t_n)\|_{W^{1,\infty}}  \|(R_{h,\rho}(\rho(t_{n-1}))-\rho(t_{n-1})\|_0\|\nabla  \e_{h,\rho}^{n}\|_0+C\tau \|u(t_n)\|_{W^{1,\infty}}  \|\rho(t_{n-1})-\rho(t_{n})\|_0\|\nabla \e_{h,\rho}^{n}\|_0\\
		\leq & C\tau \|\e_{h,u}^{n-1}\|_0^2+C\tau \|\e_{h,\rho}^{n -1}\|_0^2 + C\tau h^{r+1}+C\tau^3+\frac{D_0 \tau}{4}\|\nabla \e_{h,\rho}^{n}\|_0^2.
	\end{align*}
	Using the properties of  Ritz projection, we derive that
	\begin{align*}
		2\tau\left| \left(D_\tau (\rho^{e,n}-R_{h,\rho}(\rho^{e}(t_n))),e_{h,\rho}^{n}\right)\right| &\leq C\tau h^{r+1}\|\rho_t^{e,n}\|_{r+1}\|\nabla e_{h,\rho}^{n}\|_0\\
		&\leq C\tau h^{2r+2}+\frac{D_0 \tau }{4}\|\nabla \e_{h,\rho}^{n}\|_0^2 .
	\end{align*}
	Using Cauchy-Schwarz and Young's inequality, we have
	\begin{align*}
		2\tau |(Tr_\rho^n,\e_{h,\rho}^{n})|\leq&  C\tau \|Tr_\rho^n\|_0\|\nabla \e_{h,\rho}^{n}\|_0\\
		\leq & C\tau^3+\frac{D_0 \tau}{4}\|\nabla \e_{h,\rho}^{n}\|_0^2.
	\end{align*}
	Then, we arrive at
\begin{align}
	&\|\e_{h,\rho}^{n}\|_0^2-\|\e_{h,\rho}^{n-1}\|_0^2+\|\e_{h,\rho}^{n}-\e_{h,\rho}^{n-1}\|_0^2+D_0 \tau  \|\nabla \e_{h,\rho}^{n}\|_0^2\nonumber\\
	\leq&  C\tau \sum_{i=1}^M\|\e_{h,c^{i}}^{n-1}\|_0^2+ C\tau \|\e_{h,\phi}^{n-1}\|_0^2+ C\tau \|\e_{h,u}^{n-1}\|_0^2+C\tau h^{2r+2}+C\tau^3.\label{Er_r}
\end{align}

Taking $v=v_h$ in (\ref{3.66b}) and $q=q_h$ in (\ref{3.66c}), subtracting them from (\ref{24}) and (\ref{24b}) respectively and using the Stokes projection, we get the error equations for $u$ and $p$ as follows
\begin{align}\label{3.68}
	&(D_\tau \e_{h,u}^n,v_h)+\mu (\nabla \e_{h,u}^n, \nabla v_h)+(u_h^{n-1}\cdot\nabla u_h^n,v_h)-(u(t_n)\cdot\nabla u(t_n),v_h) \nonumber\\
	&- (\e_{h,p}^n, \nabla\cdot v_h)+ (\rho_h^{e,n}\nabla \phi_h^{n-1},v_h)-(\rho^e(t_n)\nabla \phi(t_n),v_h)\nonumber\\
	&=(D_\tau (R_h(u(t_n),p(t_n))-u(t_n)),v_h)+(Tr_u^n,v_h),\forall v_h\in X,\\
	&(\nabla\cdot \e_{h,u}^n,q_h)=0,\forall q_h\in M.\label{3.69}
\end{align}
Taking $v_h=2\tau \e_{h,u}^n$ and $q_h=2\tau \e_{h,p}^n$, we can deduce that
\begin{align*}
	&\|\e_{h,u}^n\|_0^2 - \|\e_{h,u}^{n-1}\|_0^2+\|\e_{h,u}^n - \e_{h,u}^{n-1}\|_0^2 +2\tau \mu\|\nabla \e_{h,u}^n\|_0^2\\
	&+2\tau ( u_h^{n-1}\cdot\nabla u_h^n,\e_{h,u}^n)-2\tau ( u(t_n)\cdot\nabla u(t_n),\e_{h,u}^n) \nonumber\\
	&+ 2\tau (\rho_h^{e,n}\nabla \phi_h^{n-1},\e_{h,u}^n)-2\tau (\rho^e(t_n)\nabla \phi(t_n),\e_{h,u}^n)\\
	&=2\tau (D_\tau (R_h(u(t_n),p(t_n))-u(t_n)),\e_{h,u}^n)+2\tau (Tr_u^n,\e_{h,u}^n).
\end{align*}
Adding and subtracting some terms and using Cauchy-Schwarz, Young's inequality, Taylor's formula, it yields that
\begin{align*}
	&|2\tau ((u_h^{n-1}\cdot\nabla)u_h^n,\e_{h,u}^n)-2\tau ((u(t_n)\cdot\nabla)u(t_n),\e_{h,u}^n)|\\
	\leq& 2\tau |((\e_{h,u}^{n-1}\cdot\nabla)u_h^n,\e_{h,u}^n)|+2\tau |((R_h(u(t_{n-1}), p(t_{n-1}))-u(t_{n-1}))\cdot\nabla u_h^n,\e_{h,u}^n)|\\
	&+2\tau |((u(t_{n-1})-u(t_{n}))\cdot\nabla)u_h^n,\e_{h,u}^n)|+2\tau |((u(t_n)\cdot\nabla)\e_{h,u}^n,\e_{h,u}^n)|\\
	&+2\tau |((u(t_n)\cdot\nabla)(R_h(u(t_{n}), p(t_{n}))-u(t_{n})) ,\e_{h,u}^n)|\\
	\leq & C\tau \|\e_{h,u}^{n-1}\|_0\|u_h^n\|_{W^{1,\infty}}\|\nabla \e_{h,u}^n\|_0
	+C\tau  \|R_h(u(t_{n-1}), p(t_{n-1}))-u(t_{n-1})\|_0\|u_h^n\|_{W^{1,\infty}}\|\nabla \e_{h,u}^n\|_0\\
	&+C\tau \|u(t_{n-1})-u(t_{n})\|_0\|u_h^n\|_{W^{1,\infty}}\|\nabla \e_{h,u}^n\|_0
	+C\tau \|u(t_n)\|_{W^{1,\infty}} \|R_h(u(t_{n}), p(t_{n}))-u(t_{n})\|_0 \|\nabla \e_{h,u}^n\|_0\\
	\leq & C\tau \|\e_{h,u}^{n-1}\|_0^2 +C\tau h^{2r+2}+C\tau^3+\frac{\nu\tau}{4}\|\nabla \e_{h,u}^n\|_0^2.
\end{align*}
Adding and subtracting some terms and using Cauchy-Schwarz, Young's inequality, Taylor's formula, there holds that
\begin{align*}
	&|2\tau (\rho_h^{e,n}\nabla \phi_h^{n-1},\e_{h,u}^n)-2\tau (\rho^e(t_n)\nabla \phi(t_n),\e_{h,u}^n)| \\
	\leq & 2\tau |(\e_{h,\rho}^n \nabla \phi_h^{n-1},\e_{h,u}^n) |
	+2\tau | ((R_{h,\rho}(\rho^e(t_n))-\rho^e(t_n)) \nabla \phi_h^{n-1},\e_{h,u}^n) |\\
	&+2\tau |( \rho^e(t_n)\nabla \e_{h,\phi}^{n-1} ,\e_{h,u}^n ) |
	+2\tau |( \rho^e(t_n)\nabla (R_{h,\rho}(\phi(t_{n-1}))-\phi(t_{n-1})) ,\e_{h,u}^n ) |\\
	&+2\tau |( \rho^e(t_n)\nabla (\phi(t_{n-1})-\rho^e(t_{n})) ,\e_{h,u}^n ) |\\
	\leq & C\tau \| \e_{h,\rho}^n\|_0  \| \phi_h^{n-1}\|_{W^{1,\infty}}\|\nabla \e_{h,u}^n\|_0 + C\tau \|R_{h,\rho}(\rho^e(t_n))-\rho^e(t_n)\|_0  \| \phi_h^{n-1}\|_{W^{1,\infty}}\|\nabla \e_{h,u}^n\|_0  \\
	&+C\tau \| \rho^e(t_n)\|_{\
	\infty}\|\nabla \e_{h,\phi}^{n-1}\|_0 \|\nabla \e_{h,u}^n\|_0+C\tau \| \rho^e(t_n)\|_{W^{1,\infty }}\|R_{h,\phi}(\phi(t_{n-1}))-\phi^e(t_{n-1})\|_0 \|\nabla \e_{h,u}^n\|_0\\
	& + C\tau \|\rho^e(t_n)\|_{W^{1,\infty }} \|\phi(t_{n-1})-\phi(t_{n})\|_0\|\nabla \e_{h,u}^n\|_0\\
	\leq & C\tau \| \e_{h,\rho}^n\|_0^2 + C\tau \|\e_{h,\rho}^{n-1}\|_0^2+C\tau h^{2r+2}+C\tau^3+\frac{\nu\tau}{4}\|\nabla \e_{h,u}^n\|_0^2.
\end{align*}
Using the properties of Stokes projection, Cauchy-Schwarz and Young's inequality, we deduce that
\begin{align*}
	2\tau |(D_\tau (R_h(u(t_n),p(t_n))-u(t_n)),\e_{h,u}^n)|\leq & C\tau h^{r+1}\|u_t^{n}\|_{r+1}\|\nabla \e_{h,u}^{n}\|_0\\
	\leq & C\tau h^{2r+2}+\frac{\nu\tau}{4}\|\nabla \e_{h,u}^n\|_0^2.	
\end{align*}
Using Cauchy-Schwarz and Young's inequality, there holds that
\begin{align*}
	2\tau |(Tr_u^n,\e_{h,u}^{n})|\leq&  C\tau \|Tr_u^n\|_0\|\nabla \e_{h, u}^{n}\|_0\\
	\leq & C\tau^3+\frac{\nu \tau}{4}\|\nabla \e_{h,u}^{n}\|_0^2.
\end{align*}
Then, we arrive at
\begin{align}
	&\|\e_{h,u}^n\|_0^2 - \|\e_{h,u}^{n-1}\|_0^2+\|\e_{h,u}^n - \e_{h,u}^{n-1}\|_0^2 + \tau \mu\|\nabla \e_{h,u}^n\|_0^2\nonumber\\
	\leq & C\tau \| \e_{h,\rho}^n\|_0^2 + C\tau \|\e_{h,\phi}^{n-1}\|_0^2+C\tau\| \e_{h,u}^{n-1}\|_0^2 +C\tau h^{2r+2}+C\tau^3.\label{Er_u}
\end{align}

Taking $\zeta=\zeta_h$ in (\ref{3.66d}) and subtracting it for (\ref{25}), it yields the error equation for $c^i$ as follows
\begin{align*}
	& (D_\tau \e_{h,c^i}^n,\zeta_h)+ d_i(\nabla \e_{h,c^i}^n , \nabla \zeta_h)+(u_h^{n}\cdot\nabla c_h^{i,n},\zeta_h)- (u(t_n)\cdot\nabla c^i(t_n),\zeta_h)\nonumber\\
	&+ \nu_i z_i ( c_h^{i,n}\nabla\phi_h^{n-1} ,\nabla \zeta_h)-\nu_i z_i (c^i(t_n)\nabla\phi(t_h), \nabla \zeta_h )\\
	&=2\tau (D_\tau (R_{h,c^i}(c^i(t_n))-c^i(t_n)),\zeta_h)+(Tr_{c^i}^n, \zeta_h ),i=1,\ldots,M,&\forall \zeta_h \in W_h,
\end{align*}
Taking $\zeta_h=2\tau \e_{h,c^i}^n$, we have
\begin{align*}
	&\|\e_{h,c^i}^n\|_0^2-\|\e_{h,c^i}^{n-1}\|_0^2+\|\e_{h,c^i}^n-\e_{h,c^i}^{n-1}\|_0^2+ 2\tau d_i\|\nabla \e_{h,c^i}^n\|_0^2+2\tau (u_h^{n}\cdot\nabla c_h^{i,n},\e_{h,c^i}^n)\nonumber\\
	&- 2\tau (u(t_n)\cdot\nabla c^i(t_n),\e_{h,c^i}^n)+ 2\tau \nu_i z_i ( c_h^{i,n}\nabla\phi_h^{n-1} ,\nabla \e_{h,c^i}^n)-2\tau \nu_i z_i (c^i(t_n)\nabla\phi(t_h), \nabla \e_{h,c^i}^n )\\
	&=2\tau (D_\tau (R_{h,c^i}(c^i(t_n))-c^i(t_n)),\e_{h,c^i}^n)+(Tr_{c^i}^n, \e_{h,c^i}^n ),i=1,\ldots,M.
\end{align*}
Adding and subtracting some terms, using Cauchy-Schwarz and Young's inequality, there holds that
\begin{align*}
	&|2\tau (u_h^{n}\cdot\nabla c_h^{i,n},\e_{h,c^i}^n)- 2\tau (u(t_n)\cdot\nabla c^i(t_n),\e_{h,c^i}^n)|\\
	\leq & C\tau |(\e_{h,u}^n\cdot\nabla c_h^{i,n},\e_{h,c^i}^n) |+ C\tau |((R_h(u(t_n),p(t_n))-u(t_n))\cdot\nabla c_h^{i,n},\e_{h,c^i}^n) |\\
	&+C\tau |(u(t_n)\cdot\nabla \e_{h,c^i}^n,\e_{h,c^i}^n) |+C\tau |(u(t_n)\cdot\nabla (R_{h,c}(c^i(t_n))-c^i(t_n)),\e_{h,c^i}^n) |\\
	\leq & C\tau \|\e_{h,u}^n\|_0\|c_h^{i,n}\|_{W^{1,\infty}}\|\nabla \e_{h,c^i}^n\|_0
	+C\tau \|R_h(u(t_n),p(t_n))-u(t_n)\|_0\| c_h^{i,n}\|_{W^{1,\infty}}\|\nabla \e_{h,c^i}^n\|_0\\
	&+ C\tau \|u(t_n)\|_{W^{1,\infty}}\|(R_{h,c}(c^i(t_n))-c^i(t_n))\|_0\|\nabla \e_{h,c^i}^n\|_0\\
	\leq& C\tau \|\e_{h,u}^n\|_0^2+ C\tau h^{2r+2}+\frac{d_i\tau }{4}\|\nabla \e_{h,c^i}^n\|_0.
\end{align*}
Adding and subtracting some terms, using Cauchy-Schwarz and Young's inequality, we derive that
\begin{align*}
	&|2\tau \nu_i z_i ( c_h^{i,n}\nabla\phi_h^{n-1} ,\nabla \e_{h,c^i}^n)-2\tau \nu_i z_i (c^i(t_n)\nabla\phi(t_h), \nabla \zeta_h )|\\
	\leq & 2\tau \nu_i z_i |(e_{h,c^i}^{n} \nabla\phi_h^{n-1} ,\nabla \e_{h,c^i}^n)|
	+2\tau \nu_i z_i |((R_{h,c}(c(t_n))-c(t_n) )\nabla\phi_h^{n-1} ,\nabla \e_{h,c^i}^n)|\\
	&+ 2\tau \nu_i z_i |(c^i(t_n) \nabla\e_{h,\phi}^{n-1} ,\nabla \e_{h,c^i}^n)|+2\tau \nu_i z_i |(c^i(t_n) \nabla(R_{h,\phi}(\phi(t_{n-1}))-\phi(t_{n-1})),\nabla \e_{h,c^i}^n)|\\
	&+2\tau \nu_i z_i |(c^i(t_n) \nabla(\phi(t_{n-1})-\phi(t_{n})),\nabla \e_{h,c^i}^n)|\\
	\leq & C\tau \| \e_{h,c^i}^{n}\|_0\|\phi_h^{n-1}\|_{W^{1,\infty}}\|\nabla \e_{h,c^i}^n\|_0+C\tau \|R_{h,c}(c(t_n))-c(t_n)\|_0\|\phi_h^{n-1}\|_{W^{1,\infty}}\|\nabla \e_{h,c^i}^n\|_0\\
	&	+C\tau \|c^i(t_n)\|_{W^{1,\infty}}\|\nabla \e_{h,\phi}^{n-1}\|_0\|\nabla \e_{h,c^i}^n\|_0
	+C\tau \|c^i(t_n)\|_{W^{1,\infty}} \|R_{h,\phi}(\phi(t_{n-1}))-\phi(t_{n-1})\|_0 \|\nabla \e_{h,c^i}^n\|_0\\
	\leq & C\tau \| \e_{h,c^i}^{n}\|_0^2+C\tau \|\e_{h,\rho}^{n-1}\|_0^2+C\tau h^{2r+2} +\frac{d_i\tau }{4}\|\nabla \e_{h,c^i}^n\|_0.
\end{align*}
Using Cauchy-Schwarz, Young's inequality, and the properties of Ritz projection, there holds that
\begin{align*}
	2\tau |(D_\tau (R_{h,c^i}(c^i(t_n))-c^i(t_n)),\e_{h,c^i}^n)|\leq&  C\tau h^{r+1}\|c^i_t(t_{n})\|_{r+1}\|\nabla  \e_{h,c^i}^n\|_0\\
	\leq & C\tau h^{2r+2}+\frac{d_i\tau}{4}\|\nabla  \e_{h,c^i}^n\|_0^2.
\end{align*}
Using Cauchy-Schwarz and Young's inequality, we derive that
\begin{align*}
	2\tau |(Tr_{c^i}^n, \e_{h,c^i}^n )|\leq C\tau^3+\frac{d_i\tau }{4}\|\nabla \e_{h,c^i}^n\|_0.
\end{align*}
Then, we arrive at
\begin{align*}
	&\|\e_{h,c^i}^n\|_0^2-\|\e_{h,c^i}^{n-1}\|_0^2+\|\e_{h,c^i}^n-\e_{h,c^i}^{n-1}\|_0^2+ \tau d_i\|\nabla \e_{h,c^i}^n\|_0^2\\
	\leq & C\tau \|\e_{h,u}^n\|_0^2+C\tau \| \e_{h,c^i}^{n}\|_0^2+C\tau \|\e_{h,\rho}^{n-1}\|_0^2+ C\tau^3+ C\tau h^{2r+2} ,i=1,\ldots,M.
\end{align*}
Taking sum of it over all $i$, we arrive at
\begin{align}
	&\sum_{i=1}^M\|\e_{h,c^i}^n\|_0^2-\sum_{i=1}^M\|\e_{h,c^i}^{n-1}\|_0^2+ \tau \sum_{i=1}^M d_i\|\nabla \e_{h,c^i}^n\|_0^2\nonumber\\
	\leq & C\tau \|\e_{h,u}^n\|_0^2+C\tau \sum_{i=1}^M\| \e_{h,c^i}^{n}\|_0^2+C\tau \|\e_{h,\rho}^{n-1}\|_0^2+ C\tau^3+ C\tau h^{2r+2}.\label{Er_c}
\end{align}
Taking sum of (\ref{Er_r}), (\ref{Er_u}) and (\ref{Er_c}), it yields that
\begin{align*}
	&\|\e_{h,\rho}^{n}\|_0^2-\|\e_{h,\rho}^{n-1}\|_0^2+D_0 \tau \|\nabla \e_{h,\rho}^{n}\|_0^2 +\|\e_{h,u}^n\|_0^2 - \|\e_{h,u}^{n-1}\|_0^2+ \tau \mu\|\nabla \e_{h,u}^n\|_0^2\nonumber\\
	&+\sum_{i=1}^M\|\e_{h,c^i}^n\|_0^2-\sum_{i=1}^M\|\e_{h,c^i}^{n-1}\|_0^2+ \tau \sum_{i=1}^M d_i\|\nabla \e_{h,c^i}^n\|_0^2\nonumber\\
	\leq&  C\tau \sum_{i=1}^M\|\e_{h,c^{i}}^{n-1}\|_0^2+C\tau \| \e_{h,\rho}^n\|_0^2+ C\tau \|\e_{h,\rho}^{n-1}\|_0^2+ C\tau \|\e_{h,u}^{n-1}\|_0^2 + C\tau \|\e_{h,u}^n\|_0^2\nonumber\\
	&+C\tau \sum_{i=1}^M\| \e_{h,c^i}^{n}\|_0^2+C\tau h^{2r+2}+C\tau^3.
\end{align*}
Taking sum of it over all $i$, and using Gronwall's lemma, we get
\begin{align*}
	\|\e_{h,\rho}^{n}\|_0^2+D_0 \tau \sum_{n=1}^N\|\nabla \e_{h,\rho}^{n}\|_0^2& +\|\e_{h,u}^N\|_0^2 + \tau \mu\sum_{n=1}^N\|\nabla \e_{h,u}^n\|_0^2\nonumber\\
	&+\sum_{i=1}^M\|\e_{h,c^i}^N\|_0^2+ \tau\sum_{n=1}^N \sum_{i=1}^M d_i\|\nabla \e_{h,c^i}^n\|_0^2
	\leq  Ch^{2r+2}+C\tau^2.
\end{align*}
Using triangle inequality, we have
\begin{align*}
	\|\rho^e(T)-\rho_h^{e,N}\|_0^2&+D_0 \tau h\sum_{n=1}^N\|\nabla (\rho^e(t_n)-\rho_h^{e,n})\|_0^2 +\|u(T)-u_h^N\|_0^2 + \tau h\mu\sum_{n=1}^N\|\nabla (u(t_n)-u_h^n)\|_0^2\nonumber\\
	&+\sum_{i=1}^M\|c^i(T)-c_h^{i,N}\|_0^2+ \tau h\sum_{n=1}^N \sum_{i=1}^M d_i\|\nabla c^i(t_n)-c_h^{i,n}\|_0^2
	\leq  Ch^{2r+2}+C\tau^2.
\end{align*}

By (\ref{3.69}), there holds that
\begin{align*}
	(\nabla D_\tau \e_{h,u}^n,q_h)=0,\forall q_h\in M.
\end{align*}
Taking $q_h= \e_{h,p}^n$, we have
\begin{align*}
	(\nabla \cdot D_\tau \e_{h,u}^n,\e_{h,p}^n)=0.
\end{align*}
Taking $v_h=2\tau D_\tau \e_{h,u}^n$ in (\ref{3.68}), it follows that
\begin{align*}
	&2\tau \|D_\tau \e_{h,u}^n\|_0^2+\mu \|\nabla \e_{h,u}^n\|_0-\mu \|\nabla \e_{h,u}^{n-1}\|_0^2+ \mu \|\nabla \e_{h,u}^n-\nabla \e_{h,u}^{n-1}\|_0^2 \\
	&+2 \tau ((u_h^{n-1}\cdot\nabla)u_h^n,D_\tau \e_{h,u}^n)-2\tau ((u(t_n)\cdot\nabla)u(t_n),D_\tau \e_{h,u}^n) \nonumber\\
	&- 2\tau (\e_{h,p}^n, \nabla\cdot D_\tau \e_{h,u}^n)+ 2\tau (\rho_h^{e,n}\nabla \phi_h^{n-1},D_\tau \e_{h,u}^n)-2\tau (\rho^e(t_n)\nabla \phi(t_n),D_\tau \e_{h,u}^n)\nonumber\\
	=&2\tau (D_\tau (R_h(u(t_n),p(t_n))-u(t_n)),D_\tau \e_{h,u}^n)+2\tau (Tr_u^n,D_\tau \e_{h,u}^n)
\end{align*}
Then, it yields that
\begin{align*}
	&2\tau \|D_\tau \e_{h,u}^n\|_0^2+\mu \|\nabla \e_{h,u}^n\|_0-\mu \|\nabla \e_{h,u}^{n-1}\|_0^2+ \mu \|\nabla \e_{h,u}^n-\nabla \e_{h,u}^{n-1}\|_0^2 \\
	&+2 \tau ((u_h^{n-1}\cdot\nabla)u_h^n,D_\tau \e_{h,u}^n)-2\tau ((u(t_n)\cdot\nabla)u(t_n),D_\tau \e_{h,u}^n) \nonumber\\
	&+ 2\tau (\rho_h^{e,n}\nabla \phi_h^{n-1},D_\tau \e_{h,u}^n)-2\tau (\rho^e(t_n)\nabla \phi(t_n),D_\tau \e_{h,u}^n)\nonumber\\
	=&2\tau (D_\tau (R_h(u(t_n),p(t_n))-u(t_n)),D_\tau \e_{h,u}^n)+2\tau (Tr_u^n,D_\tau \e_{h,u}^n)
\end{align*}
Adding and subtracting some terms, using Cauchy-Schwarz and Young's inequality, we derive that
\begin{align*}
	&|2\tau ((u_h^{n-1}\cdot\nabla)u_h^n,D_\tau \e_{h,u}^n)-2\tau ((u(t_n)\cdot\nabla)u(t_n),D_\tau \e_{h,u}^n)|\\
	\leq& 2\tau |((\e_{h,u}^{n-1}\cdot\nabla)u_h^n,D_\tau\e_{h,u}^n)|+2\tau |((R_h(u(t_{n-1}), p(t_{n-1}))-u(t_{n-1}))\cdot\nabla u_h^n,D_\tau\e_{h,u}^n)|\\
	&+2\tau |((u(t_{n-1})-u(t_{n}))\cdot\nabla u_h^n,D_\tau\e_{h,u}^n)|+2\tau |((u(t_n)\cdot\nabla)\e_{h,u}^n,D_\tau\e_{h,u}^n)|\\
	&+2\tau |((u(t_n)\cdot\nabla)(R_h(u(t_{n}), p(t_{n}))-u(t_{n})) ,D_\tau\e_{h,u}^n)|\\
	\leq & C\tau \|\e_{h,u}^{n-1}\|_0\|u_h^n\|_{w^{1,\infty}}\|D_\tau \e_{h,u}^n\|_0
	+C\tau  \|R_h(u(t_{n-1}), p(t_{n-1}))-u(t_{n-1})\|_0\|u_h^n\|_{W^{1,\infty}}\|D_\tau \e_{h,u}^n\|_0\\
	&+C\tau \|u(t_{n-1})-u(t_{n})\|_0\|u_h^n\|_{W^{1,\infty}}\|D_\tau \e_{h,u}^n\|_0
	+C\tau \|u(t_n)\|_{W^{1,\infty}} \|R_h(u(t_{n}), p(t_{n}))-u(t_{n})\|_0 \|D_\tau \e_{h,u}^n\|_0\\
	\leq & C\tau h^{2r+2}+C\tau^3+\frac{\tau}{4}\|D_\tau \e_{h,u}^n\|_0^2.
\end{align*}
Adding and subtracting some terms, using Cauchy-Schwarz and Young's inequality, it follows by
\begin{align*}
	&|2\tau (\rho_h^{e,n}\nabla \phi_h^{n-1},D_\tau\e_{h,u}^n)-2\tau (\rho^e(t_n)\nabla \phi(t_n),D_\tau\e_{h,u}^n)| \\
	\leq & 2\tau |(e_{h,\rho}^n \nabla \phi_h^{n-1},D_\tau\e_{h,u}^n) |
	+2\tau | ((R_{h,\rho}(\rho^e(t_n))-\rho^e(t_n)) \nabla \phi_h^{n-1},D_\tau\e_{h,u}^n) |\\
	&+2\tau |( \rho^e(t_n)\nabla \e_{h,\phi}^{n-1} ,D_\tau\e_{h,u}^n ) |
	+2\tau |( \rho^e(t_n)\nabla (R_{h,\rho}(\phi(t_{n-1}))-\rho^e(t_{n-1})) ,D_\tau\e_{h,u}^n ) |\\
	&+2\tau |( \rho^e(t_n)\nabla (\phi(t_{n-1})-\rho^e(t_{n})) ,D_\tau\e_{h,u}^n ) |\\
	\leq & C\tau \| \e_{h,\rho}^n\|_0  \| \phi_h^{n-1}\|_{W^{1,\infty}}\|D_\tau\e_{h,u}^n\|_0 + C\tau \|R_{h,\rho}(\rho^e(t_n))-\rho^e(t_n)\|_0  \| \phi_h^{n-1}\|_{W^{1,\infty}}\|D_\tau\e_{h,u}^n\|_0  \\
	&+C\tau \| \rho^e(t_n)\|_{\
	\infty}\|\nabla \e_{h,\phi}^{n-1}\|_0 \|D_\tau\e_{h,u}^n\|_0+C\tau \| \rho^e(t_n)\|_{W^{1,\infty }}\|R_{h,\phi}(\phi(t_{n-1}))-\phi^e(t_{n-1})\|_0 \|D_\tau\e_{h,u}^n\|_0\\
	& + C\tau \|\rho^e(t_n)\|_{W^{1,\infty }} \|\phi(t_{n-1})-\phi(t_{n})\|_0\|D_\tau \e_{h,u}^n\|_0\\
	\leq & C\tau^3+\frac{\tau}{4}\|D_\tau \e_{h,u}^n\|_0^2.
\end{align*}
Using Cauchy-Schwarz, Young's inequality and the properties of Ritz projection, we derive that
\begin{align*}
	2\tau |(D_\tau (R_h(u(t_n),p(t_n))-u(t_n)),D_\tau \e_{h,u}^n)|\leq & C\tau h^{r+1}\|u_t^{n}\|_{r+1}\|D_\tau  \e_{h,u}^{n}\|_0\\
	\leq & C\tau h^{2r+2}+\frac{\tau}{4}\|D_\tau \e_{h,u}^n\|_0^2.	
\end{align*}
Using Cauchy-Schwarz and Young's inequality, we have
\begin{align*}
	2\tau |(Tr_u^n,D_\tau \e_{h,u}^{n})|\leq&  C\tau \|Tr_u^n\|_0\|D_\tau  \e_{h, u}^{n}\|_0\\
	\leq & C\tau^3+\frac{\tau}{4}\|\nabla \e_{h,u}^{n}\|_0^2.
\end{align*}
Then, we have
\begin{align*}
	2\tau \|D_\tau \e_{h,u}^n\|_0^2+\mu \|\nabla \e_{h,u}^n\|_0-\mu \|\nabla \e_{h,u}^{n-1}\|_0^2 \leq C\tau h^{2r+2}+C\tau^3.
\end{align*}
Taking sum of it over all $n$, we deduce
\begin{align*}
	2\tau \sum_{n=1}^N\|D_\tau \e_{h,u}^n\|_0^2+\mu \|\nabla \e_{h,u}^N\|_0 \leq C h^{2r+2}+C\tau^2.
\end{align*}

In order to get the error estimation, (\ref{3.68}) can be rewritten as
\begin{align*}
	(\e_{h,p}^n,\nabla \cdot v_h)=&(D_\tau \e_{h,u}^n,v_h)+\mu (\nabla \e_{h,u}^n, \nabla v_h)+((u_h^{n-1}\cdot\nabla)u_h^n,v_h)\nonumber\\
	&-((u(t_n)\cdot\nabla)u(t_n),v_h)+ (\rho_h^{e,n}\nabla \phi_h^{n-1},v_h)-(\rho^e(t_n)\nabla \phi(t_n),v_h)\nonumber\\
	&-(D_\tau (R_h(u(t_n),p(t_n))-u(t_n)),v_h)+(Tr_u^n,v_h),\forall v_h\in X.
\end{align*}
By the discrete LBB condition \refe{LBB} and the estimates obtained above, we get
\begin{align*}
	\|\e_{h,p}^n\|_0 \leq &\|D_\tau \e_{h,u}^n\|_{-1}+\mu \|\nabla \e_{h,u}^n\|_0 + C\tau \| \e_{h,\rho}^n\|_0 + C\tau \|\e_{h,\phi}^{n-1}\|_0+C\tau\| \e_{h,u}^{n-1}\|_0\\
	\leq & Ch^r+C\tau.
\end{align*}
which in turn produces
\begin{eqnarray*}
\sum_{n=1}^{N} \tau \| \e_{h,p}^n \|_0^2\leq C(\tau^2+h^{2r}).
\end{eqnarray*}
By triangle inequality, we have
\begin{eqnarray*}
	\sum_{n=1}^{N} \tau \| p(t_n)-p_h^n\|_0^2\leq C(\tau^2+h^{2r}).
\end{eqnarray*}	
\end{proof}
 \section{Numerical results}	
 In this section, numerical computations are used to show the effectiveness of the proposed method. Firstly, the electro-osmotic flow in T-junction microchannels with different viscosity is numerically investigated. Secondly, the effect of the roughness in microchannels for the electro-osmotic flow is studied.

\subsection{Analytical problems}
In this subsection, we give some numerical results for the electroneutral micro-fluids  equation in the domain $\Omega = [0 , 1e-3]\times [0 , 1e-3] $.

\begin{align}\label{requ2}
	\left\{\begin {array}{rll}
	\rho^e_t-D_0\triangle \rho^e -\sum_{i=1}^M \nu_iz_i^2\cdot\nabla(c^i\nabla \phi)+u\cdot\nabla \rho^e=f_\rho,&\\
	u_t-\mu \triangle u+u\cdot\nabla u+\nabla p-\rho^e\nabla \phi=f_u,&\\
	\nabla\cdot u=0,&\\
	c^i_t-d_i\triangle c^i -\nu_i z_i \cdot\nabla(c^i\nabla\phi )+u\cdot\nabla c^i=f_{c^i},& i=1,2,3,\\
	-\varepsilon \triangle \phi-\rho^e=f_\phi.
	\end{array}\right.
	\end{align}
Here, we consider the analytical solution as follows
\begin{align*}
	\rho^e& =-\frac{2}{\pi} \cos(t)(\cos(\pi x)+\cos(\pi y)),\\
u_1&=(x^2 (y-1)^2+y) \cos(t),\\
u_2&=(-2.0 x(y-1)^3) \cos(t)/3.0+(2-\pi\sin(\pi x))\cos(t),\\
p&=(2-\pi\sin(\pi x)) \sin(0.5\pi y)\cos(t),\\
c_1&=\frac{t}{\pi }(\cos(\pi x)+\cos(\pi y)),\\
c_2&=\frac{10t}{\pi }(\cos(\pi x)+\cos(\pi y)),\\
c_3&=\frac{5t}{\pi }(\cos(\pi x)+\cos(\pi y)),\\
\phi&=(2-\pi\sin(\pi x)) (1-y-\cos(\pi y))\cos(t).
\end{align*}
where the velocity field is $u=(u_1,u_2)$,  the boundary and initial conditions in (\ref{requ2}) are set equal to the analytical solution,  the force terms $f_\rho$, $f_u$, $f_\phi$ and $f_{c^i}, i=1,2, 3$ are given by evaluating the momentum equation of problem (\ref{requ2}) for the analytical solution, the mobilities are $\nu_1=5\times 10^{-2}$,  $\nu_2=3\times 10^{-2}$ and $\nu_3=3\times 10^{-2}$, and the valence numbers are $z_1=1.0$, $z_2=-1.0 $ and $z_3=-2.0$. Firstly, we choose  $D_0= d_1=d_2=d_3=1e-6$. As the flow field will get into the steady state vary fast, we choose $T=1e-8$. The finite element spaces are Mini finite element spaces ($P1b-P1$) for fluids and $P1b$ for charge density of species. The finite element steps are chosen as $\tau=h^2$, $h=1e-3/G$, and $G=8,16,24,32,40,48$. The numerical results with $\mu =1$ are given in Table \ref{Ta1}, \ref{Ta2} and \ref{Ta3}.  And The numerical results with $\mu =1e-3$ are given in Table  \ref{Ta4}, \ref{Ta5} and \ref{Ta6}. We can see that the errors go smaller as the spacial step goes smaller, it is second-order convergence with respect to the spatial mesh refinement. To show the robustness of the model, we choose $D_0=d_1=d_2=d_3=1e-8$, $\mu = 1 $ and $1e-3$, respectively. The numerical results are given in Table \ref{Ta7} to Table \ref{Ta12}. The convergence rates are optimal for fluids and charge density of species in $H^1$-norm. From the numerical results, we can see that our numerical method is robust and has an optimal convergence order.

 \begin{table}[h!]
	\tabcolsep 0pt \caption{The numerical results for $d_1=d_2=d_3=1e-6$ and $\nu_1=\nu_2=\nu_3=1$ for different $h$ at $T=1e-8$.}\label{Ta1} \vspace*{-18pt}
	\begin{center}
	\def\temptablewidth{1.1\textwidth}
	{\rule{\temptablewidth}{1pt}}
	\begin{tabular*}{\temptablewidth}{@{\extracolsep{\fill}}llllllllllllllllllllllll}
	$\footnotesize G $ & $\footnotesize\|u_h^N-u(T)\|_0$ &  $\footnotesize\|\nabla(u_{h}^N-u(T))\|_0$  &$\footnotesize\|c_h^{1,N}-c^1(T)\|_0$ & $\footnotesize\|\nabla(c_h^{1,N}-c^1(T))\|_0$& $\footnotesize\|c_h^{2,N}-c^2(T)\|_0$\\
	\hline
	8  & 2.68703e-12  &  1.6669e-07 & 6.84005e-11  &   7.07065e-06 & 3.94352e-10  \\
16 & 6.4933e-13  &  8.07642e-08 & 1.78455e-11  &   3.76939e-06 & 1.02007e-10  \\
24  & 2.84653e-13  &  5.33769e-08 & 8.04204e-12  &   2.56408e-06 & 4.58042e-11  \\
32  & 1.59358e-13  &  3.98753e-08 & 4.55495e-12  &   1.9409e-06 & 2.58628e-11  \\
40  & 1.01433e-13  &  3.18295e-08 & 2.92811e-12  &   1.56073e-06 & 1.65828e-11  \\
48   & 7.0348e-14  &  2.6487e-08 & 2.04088e-12  &   1.30387e-06 & 1.15188e-11 
	\end{tabular*}
	{\rule{\temptablewidth}{1pt}}
	\end{center}
	\end{table}
	
	\begin{table}[h!]
	\tabcolsep 0pt \caption{The numerical results for $d_1=d_2=d_3=1e-6$ and $\nu_1=\nu_2=\nu_3=1$ for different $h$ at $T=1e-8$.}\label{Ta2} \vspace*{-18pt}
	\begin{center}
	\def\temptablewidth{1.1\textwidth}
	{\rule{\temptablewidth}{1pt}}
	\begin{tabular*}{\temptablewidth}{@{\extracolsep{\fill}}llllllllllllllllllllllll}
	$\footnotesize G $ &   
	 $\footnotesize\|\nabla(c_h^{2,N}-c^2(T))\|_0$ & $\footnotesize\|c_h^{3,N}-c^3(T)\|_0$& $\footnotesize\|\nabla(c_h^{3,N}-c^3(T))\|_0$ & $\footnotesize\|\rho^{e,N}-\rho^e (T)\|_0$& $\footnotesize\|\nabla(\rho^{e,N}-\rho^e (T))\|_0$\\
	\hline
	8  &  3.83571e-05 & 3.98279e-10  &   3.96121e-05  &   1.73877e-10  &   3.80776e-06 \\
16   &  2.05526e-05 & 1.03338e-10  &   2.11908e-05  &   4.63217e-11  &   1.42144e-06 \\
24  &   1.40019e-05 & 4.64605e-11  &   1.44297e-05  &   2.21401e-11  &   8.28182e-07 \\
32  &   1.06066e-05 & 2.62614e-11  &   1.09281e-05  &   1.3842e-11  &   5.92291e-07 \\
40  &   8.53273e-06 & 1.68533e-11  &   8.79019e-06  &   9.74172e-12  &   4.62255e-07\\
48   &   7.13043e-06 & 1.17203e-11  &   7.34491e-06  &   7.73827e-12  &   3.96696e-07
	\end{tabular*}
	{\rule{\temptablewidth}{1pt}}
	\end{center}
	\end{table}

\begin{table}[h!]
	\tabcolsep 0pt \caption{The convergence rates for for $d_1=d_2=d_3=1e-6$ and $\nu_1=\nu_2=\nu_3=1$ for different $h$ at $T=1e-8$.}\label{Ta3} \vspace*{-18pt}
	\begin{center}
	\def\temptablewidth{1.1\textwidth}
	{\rule{\temptablewidth}{1pt}}
	\begin{tabular*}{\temptablewidth}{@{\extracolsep{\fill}}llllllllllllllllllllllll}
	$\footnotesize G $ & \footnotesize $u-{L^2}$ &  \footnotesize  $u-{H^1}$  & \footnotesize  $c_1-{L^2}$ & \footnotesize $c_1-{H^1}$& \footnotesize $c_2-{L^2}$ & \footnotesize $c_2-{H^1} $ & \footnotesize  $c_3-{L^2}$ & \footnotesize  $c_3-{H^1}$ & \footnotesize  $\rho^e-{L^2}$ & \footnotesize  $\rho^e-{H^1}$\\
	\hline
	16 & 2.0490  &    1.0454  &    1.9384  &    0.9075  &    1.9508  &    0.9002  &    1.9464  &    0.9025  &    1.9083  &    1.4216  \\
    24 &    2.0339  &    1.0214  &    1.9658  &    0.9503  &    1.9747  &    0.9466  &    1.9716  &    0.9477  &    1.8207  &    1.3323 \\
    32 &    2.0165  &    1.0137  &    1.9760  &    0.9679  &    1.9868  &    0.9654  &    1.9831  &    0.9662  &    1.6326  &    1.1653 \\
    40 &    2.0245  &    1.0099  &    1.9801  &    0.9769  &    1.9917  &    0.9750  &    1.9878  &    0.9756  &    1.5743  &    1.1109 \\
    48 &    2.0071  &    1.0078  &    1.9799  &    0.9863  &    1.9986  &    0.9847  &    1.9922  &    0.9852  &    1.2628  &  0.8389
	\end{tabular*}
	{\rule{\temptablewidth}{1pt}}
	\end{center}
	\end{table}

	\begin{table}[h!]
		\tabcolsep 0pt \caption{The numerical results for $d_1=d_2=d_3=1e-6$ and $\nu_1=\nu_2=\nu_3=1e-3$ for different $h$ at $T=1e-8$.}\label{Ta4} \vspace*{-18pt}
		\begin{center}
		\def\temptablewidth{1.1\textwidth}
		{\rule{\temptablewidth}{1pt}}
		\begin{tabular*}{\temptablewidth}{@{\extracolsep{\fill}}llllllllllllllllllllllll}
		$\footnotesize G $ & $\footnotesize\|u_h^N-u(T)\|_0$ &  $\footnotesize\|\nabla(u_{h}^N-u(T))\|_0$  &$\footnotesize\|c_h^{1,N}-c^1(T)\|_0$ & $\footnotesize\|\nabla(c_h^{1,N}-c^1(T))\|_0$& $\footnotesize\|c_h^{2,N}-c^2(T)\|_0$\\
		\hline
		8   & 2.33348e-12  &  2.12195e-07 & 6.84005e-11  &   7.07065e-06 & 3.94352e-10  \\
		16   & 5.61143e-13  &  9.36172e-08 & 1.78455e-11  &   3.76939e-06 & 1.02007e-10  \\
		24  & 2.47508e-13  &  5.70626e-08 & 8.04204e-12  &   2.56408e-06 & 4.58042e-11  \\
		32   & 1.41076e-13  &  4.06433e-08 & 4.55495e-12  &   1.9409e-06 & 2.58628e-11  \\
		40  & 9.1234e-14  &  3.20474e-08 & 2.92811e-12  &   1.56073e-06 & 1.65828e-11  \\
		48  & 6.38807e-14  &  2.65549e-08 & 2.04088e-12  &   1.30387e-06 & 1.15188e-11 		 
		\end{tabular*}
		{\rule{\temptablewidth}{1pt}}
		\end{center}
		\end{table}
		
		\begin{table}[h!]
		\tabcolsep 0pt \caption{The numerical results for $d_1=d_2=d_3=1e-6$ and $\nu_1=\nu_2=\nu_3=1e-3$ for different $h$ at $T=1e-8$.}\label{Ta5} \vspace*{-18pt}
		\begin{center}
		\def\temptablewidth{1.1\textwidth}
		{\rule{\temptablewidth}{1pt}}
		\begin{tabular*}{\temptablewidth}{@{\extracolsep{\fill}}llllllllllllllllllllllll}
		$\footnotesize G $ &   
		 $\footnotesize\|\nabla(c_h^{2,N}-c^2(T))\|_0$ & $\footnotesize\|c_h^{3,N}-c^3(T)\|_0$& $\footnotesize\|\nabla(c_h^{3,N}-c^3(T))\|_0$ & $\footnotesize\|\rho^{e,N}-\rho^e (T)\|_0$& $\footnotesize\|\nabla(\rho^{e,N}-\rho^e (T))\|_0$\\
		\hline
8  &   3.83571e-05 & 3.98279e-10  &   3.96121e-05  &   1.73877e-10  &   3.80776e-06\\
16 &  2.05526e-05 & 1.03338e-10  &   2.11908e-05  &   4.63217e-11  &   1.42144e-06 \\
24  &  1.40019e-05 & 4.64605e-11  &   1.44297e-05  &   2.21401e-11  &   8.28182e-07 \\
32  &   1.06066e-05 & 2.62614e-11  &   1.09281e-05  &   1.3842e-11  &   5.92291e-07 \\
40   & 8.53273e-06 & 1.68533e-11  &   8.79019e-06  &   9.74172e-12  &   4.62255e-07 \\
48   &   7.13043e-06 & 1.17203e-11  &   7.34491e-06  &   7.73827e-12  &   3.96696e-07
		\end{tabular*}
		{\rule{\temptablewidth}{1pt}}
		\end{center}
		\end{table}

	\begin{table}[h!]
		\tabcolsep 0pt \caption{The convergence rates for for $d_1=d_2=d_3=1e-6$ and $\nu_1=\nu_2=\nu_3=1e-3$ for different $h$ at $T=1e-8$.}\label{Ta6} \vspace*{-18pt}
		\begin{center}
		\def\temptablewidth{1.1\textwidth}
		{\rule{\temptablewidth}{1pt}}
		\begin{tabular*}{\temptablewidth}{@{\extracolsep{\fill}}llllllllllllllllllllllll}
		$\footnotesize G $ & \footnotesize $u-{L^2}$ &  \footnotesize  $u-{H^1}$  & \footnotesize  $c_1-{L^2}$ & \footnotesize $c_1-{H^1}$& \footnotesize $c_2-{L^2}$ & \footnotesize $c_2-{H^1} $ & \footnotesize  $c_3-{L^2}$ & \footnotesize  $c_3-{H^1}$ & \footnotesize  $\rho^e-{L^2}$ & \footnotesize  $\rho^e-{H^1}$\\
		\hline
		16 &     2.0560  &    1.1805  &    1.9384  &    0.9075  &    1.9508  &    0.9002  &    1.9464  &    0.9025  &    1.9083  &    1.4216 \\
		24 & 2.0188  &    1.2210  &    1.9658  &    0.9503  &    1.9747  &    0.9466  &    1.9716  &    0.9477  &    1.8207  &    1.3323 \\
		32 & 1.9540  &    1.1795  &    1.9760  &    0.9679  &    1.9868  &    0.9654  &    1.9831  &    0.9662  &    1.6326  &    1.1653 \\
		40 &  1.9533  &    1.0649  &    1.9801  &    0.9769  &    1.9917  &    0.9750  &    1.9878  &    0.9756  &    1.5743  &    1.1109 \\
		48 & 1.9548  &    1.0312  &    1.9799  &    0.9863  &    1.9986  &    0.9847  &    1.9922  &    0.9852  &    1.2628  &  0.8389
		\end{tabular*}
		{\rule{\temptablewidth}{1pt}}
		\end{center}
		\end{table}	
		\begin{table}[h!]
			\tabcolsep 0pt \caption{The numerical results for $d_1=d_2=d_3=1e-8$ and $\nu_1=\nu_2=\nu_3=1$ for different $h$ at $T=1e-8$.}\label{Ta7} \vspace*{-18pt}
			\begin{center}
			\def\temptablewidth{1.1\textwidth}
			{\rule{\temptablewidth}{1pt}}
			\begin{tabular*}{\temptablewidth}{@{\extracolsep{\fill}}llllllllllllllllllllllll}
			$\footnotesize G $ & $\footnotesize\|u_h^N-u(T)\|_0$ &  $\footnotesize\|\nabla(u_{h}^N-u(T))\|_0$  &$\footnotesize\|c_h^{1,N}-c^1(T)\|_0$ & $\footnotesize\|\nabla(c_h^{1,N}-c^1(T))\|_0$& $\footnotesize\|c_h^{2,N}-c^2(T)\|_0$\\
			\hline
			8  & 2.68703e-12  &  1.6669e-07 & 6.84116e-11  &   7.0718e-06 & 3.94412e-10 \\
			16 & 6.4933e-13  &  8.07642e-08 & 1.78546e-11  &   3.77131e-06 & 1.02056e-10  \\
			24 & 2.84653e-13  &  5.33769e-08 & 8.05036e-12  &   2.56672e-06 & 4.58494e-11  \\
			32 & 1.59358e-13  &  3.98753e-08 & 4.56364e-12  &   1.94459e-06 & 2.59102e-11  \\
			40 & 1.01433e-13  &  3.18295e-08 & 2.93624e-12  &   1.56506e-06 & 1.66274e-11  \\
			48 & 7.03479e-14  &  2.6487e-08 & 2.04938e-12  &   1.30931e-06 & 1.15654e-11  
			\end{tabular*}
			{\rule{\temptablewidth}{1pt}}
			\end{center}
			\end{table}
			
			\begin{table}[h!]
			\tabcolsep 0pt \caption{The numerical results for $d_1=d_2=d_3=1e-8$ and $\nu_1=\nu_2=\nu_3=1$ for different $h$ at $T=1e-8$.}\label{Ta8} \vspace*{-18pt}
			\begin{center}
			\def\temptablewidth{1.1\textwidth}
			{\rule{\temptablewidth}{1pt}}
			\begin{tabular*}{\temptablewidth}{@{\extracolsep{\fill}}llllllllllllllllllllllll}
			$\footnotesize G $ &   
			 $\footnotesize\|\nabla(c_h^{2,N}-c^2(T))\|_0$ & $\footnotesize\|c_h^{3,N}-c^3(T)\|_0$& $\footnotesize\|\nabla(c_h^{3,N}-c^3(T))\|_0$ & $\footnotesize\|\rho^{e,N}-\rho^e (T)\|_0$& $\footnotesize\|\nabla(\rho^{e,N}-\rho^e (T))\|_0$\\
			\hline
			8   &   3.83633e-05 & 3.98341e-10  &   3.96186e-05  &   1.73878e-10  &   3.80802e-06 \\
			16 &   2.0563e-05 & 1.03389e-10  &   2.12015e-05  &   4.63224e-11  &   1.42174e-06 \\
			24  &   1.40163e-05 & 4.65073e-11  &   1.44446e-05  &   2.21407e-11  &   8.28527e-07 \\
			32 &   1.06268e-05 & 2.63104e-11  &   1.09489e-05  &   1.38425e-11  &   5.92728e-07 \\
			40 &   8.5564e-06 & 1.68993e-11  &   8.81459e-06  &   9.74214e-12  &   4.62735e-07 \\
			48 &   7.16019e-06 & 1.17684e-11  &   7.37559e-06  &   7.73869e-12  &   3.97282e-07		
			\end{tabular*}
			{\rule{\temptablewidth}{1pt}}
			\end{center}
			\end{table}

		\begin{table}[h!]
			\tabcolsep 0pt \caption{The convergence rates for for $d_1=d_2=d_3=1e-8$ and $\nu_1=\nu_2=\nu_3=1$ for different $h$ at $T=1e-8$.}\label{Ta9} \vspace*{-18pt}
			\begin{center}
			\def\temptablewidth{1.1\textwidth}
			{\rule{\temptablewidth}{1pt}}
			\begin{tabular*}{\temptablewidth}{@{\extracolsep{\fill}}llllllllllllllllllllllll}
			$\footnotesize G $ & \footnotesize $u-{L^2}$ &  \footnotesize  $u-{H^1}$  & \footnotesize  $c_1-{L^2}$ & \footnotesize $c_1-{H^1}$& \footnotesize $c_2-{L^2}$ & \footnotesize $c_2-{H^1} $ & \footnotesize  $c_3-{L^2}$ & \footnotesize  $c_3-{H^1}$ & \footnotesize  $\rho^e-{L^2}$ & \footnotesize  $\rho^e-{H^1}$\\
			\hline
			16 & 2.0490  &    1.0454  &    1.9379  &    0.9070  &    1.9503  &    0.8997  &    1.9459  &    0.9020  &    1.9083  &    1.4214 \\
			24 & 2.0339  &    1.0214  &    1.9645  &    0.9490  &    1.9734  &    0.9453  &    1.9703  &    0.9464  &    1.8206  &    1.3318 \\
			32 & 2.0165  &    1.0137  &    1.9730  &    0.9649  &    1.9839  &    0.9623  &    1.9801  &    0.9632  &    1.6326  &    1.1642 \\
			40 & 2.0245  &    1.0099  &    1.9763  &    0.9730  &    1.9879  &    0.9711  &    1.9839  &    0.9717  &    1.5742  &    1.1095 \\
			48 & 2.0071  &    1.0078  &    1.9723  &    0.9786  &    1.9912  &    0.9771  &    1.9847  &    0.9776  &    1.2628  &  0.8365
			\end{tabular*}
			{\rule{\temptablewidth}{1pt}}
			\end{center}
			\end{table}	
		\begin{table}[h!]
			\tabcolsep 0pt \caption{The numerical results for $d_1=d_2=d_3=1e-8$ and $\nu_1=\nu_2=\nu_3=1e-3$ for different $h$ at $T=1e-8$.}\label{Ta10} \vspace*{-18pt}
			\begin{center}
			\def\temptablewidth{1.1\textwidth}
			{\rule{\temptablewidth}{1pt}}
			\begin{tabular*}{\temptablewidth}{@{\extracolsep{\fill}}llllllllllllllllllllllll}
			$\footnotesize G $ & $\footnotesize\|u_h^N-u(T)\|_0$ &  $\footnotesize\|\nabla(u_{h}^N-u(T))\|_0$  &$\footnotesize\|c_h^{1,N}-c^1(T)\|_0$ & $\footnotesize\|\nabla(c_h^{1,N}-c^1(T))\|_0$& $\footnotesize\|c_h^{2,N}-c^2(T)\|_0$\\
			\hline
			8  & 2.33348e-12  &  2.12195e-07 & 6.84116e-11  &   7.0718e-06 & 3.94412e-10  \\
			16 & 5.61143e-13  &  9.36172e-08 & 1.78546e-11  &   3.77131e-06 & 1.02056e-10  \\
			24 & 2.47508e-13  &  5.70626e-08 & 8.05036e-12  &   2.56672e-06 & 4.58494e-11  \\
			32  & 1.41076e-13  &  4.06433e-08 & 4.56364e-12  &   1.94459e-06 & 2.59102e-11  \\
			40  & 9.1234e-14  &  3.20474e-08 & 2.93624e-12  &   1.56506e-06 & 1.66274e-11  \\
			48  & 6.38807e-14  &  2.65549e-08 & 2.04938e-12  &   1.30931e-06 & 1.15654e-11  
			\end{tabular*}
			{\rule{\temptablewidth}{1pt}}
			\end{center}
			\end{table}
			
			\begin{table}[h!]
			\tabcolsep 0pt \caption{The numerical results for $d_1=d_2=d_3=1e-8$ and $\nu_1=\nu_2=\nu_3=1$ for different $h$ at $T=1e-8$.}\label{Ta11} \vspace*{-18pt}
			\begin{center}
			\def\temptablewidth{1.1\textwidth}
			{\rule{\temptablewidth}{1pt}}
			\begin{tabular*}{\temptablewidth}{@{\extracolsep{\fill}}llllllllllllllllllllllll}
			$\footnotesize G $ &   
			 $\footnotesize\|\nabla(c_h^{2,N}-c^2(T))\|_0$ & $\footnotesize\|c_h^{3,N}-c^3(T)\|_0$& $\footnotesize\|\nabla(c_h^{3,N}-c^3(T))\|_0$ & $\footnotesize\|\rho^{e,N}-\rho^e (T)\|_0$& $\footnotesize\|\nabla(\rho^{e,N}-\rho^e (T))\|_0$\\
			\hline
			8 &   3.83633e-05 & 3.98341e-10  &   3.96186e-05  &   1.73878e-10  &   3.80802e-06\\
			16  &  2.0563e-05 & 1.03389e-10  &   2.12015e-05  &   4.63224e-11  &   1.42174e-06 \\
			24    &   1.40163e-05 & 4.65073e-11  &   1.44446e-05  &   2.21407e-11  &   8.28527e-07 \\
			32  &   1.06268e-05 & 2.63104e-11  &   1.09489e-05  &   1.38425e-11  &   5.92728e-07 \\
			40   &   8.5564e-06 & 1.68993e-11  &   8.81459e-06  &   9.74214e-12  &   4.62735e-07 \\
			48  &   7.16019e-06 & 1.17684e-11  &   7.37559e-06  &   7.73869e-12  &   3.97282e-07	
			\end{tabular*}
			{\rule{\temptablewidth}{1pt}}
			\end{center}
			\end{table}

		\begin{table}[h!]
			\tabcolsep 0pt \caption{The convergence rates for for $d_1=d_2=d_3=1e-8$ and $\nu_1=\nu_2=\nu_3=1$ for different $h$ at $T=1e-8$.}\label{Ta12} \vspace*{-18pt}
			\begin{center}
			\def\temptablewidth{1.1\textwidth}
			{\rule{\temptablewidth}{1pt}}
			\begin{tabular*}{\temptablewidth}{@{\extracolsep{\fill}}llllllllllllllllllllllll}
			$\footnotesize G $ & \footnotesize $u-{L^2}$ &  \footnotesize  $u-{H^1}$  & \footnotesize  $c_1-{L^2}$ & \footnotesize $c_1-{H^1}$& \footnotesize $c_2-{L^2}$ & \footnotesize $c_2-{H^1} $ & \footnotesize  $c_3-{L^2}$ & \footnotesize  $c_3-{H^1}$ & \footnotesize  $\rho^e-{L^2}$ & \footnotesize  $\rho^e-{H^1}$\\
			\hline
			16 & 2.0560  &    1.1805  &    1.9379  &    0.9070  &    1.9503  &    0.8997  &    1.9459  &    0.9020  &    1.9083  &    1.4214 \\
			24 & 2.0188  &    1.2210  &    1.9645  &    0.9490  &    1.9734  &    0.9453  &    1.9703  &    0.9464  &    1.8206  &    1.3318 \\
			32 & 1.9540  &    1.1795  &    1.9730  &    0.9649  &    1.9839  &    0.9623  &    1.9801  &    0.9632  &    1.6326  &    1.1642 \\
			40 & 1.9533  &    1.0649  &    1.9763  &    0.9730  &    1.9879  &    0.9711  &    1.9839  &    0.9717  &    1.5742  &    1.1095 \\
			48 & 1.9548  &    1.0312  &    1.9723  &    0.9786  &    1.9912  &    0.9771  &    1.9847  &    0.9776  &    1.2628  &  0.8365 
			\end{tabular*}
			{\rule{\temptablewidth}{1pt}}
			\end{center}
			\end{table}	
\subsection{Numerical results for  electro-osmotic flow in T-junction microchannels}

\begin{figure}
	\centering
	\subfigure[Geometry model]{\includegraphics[width=0.4\linewidth]{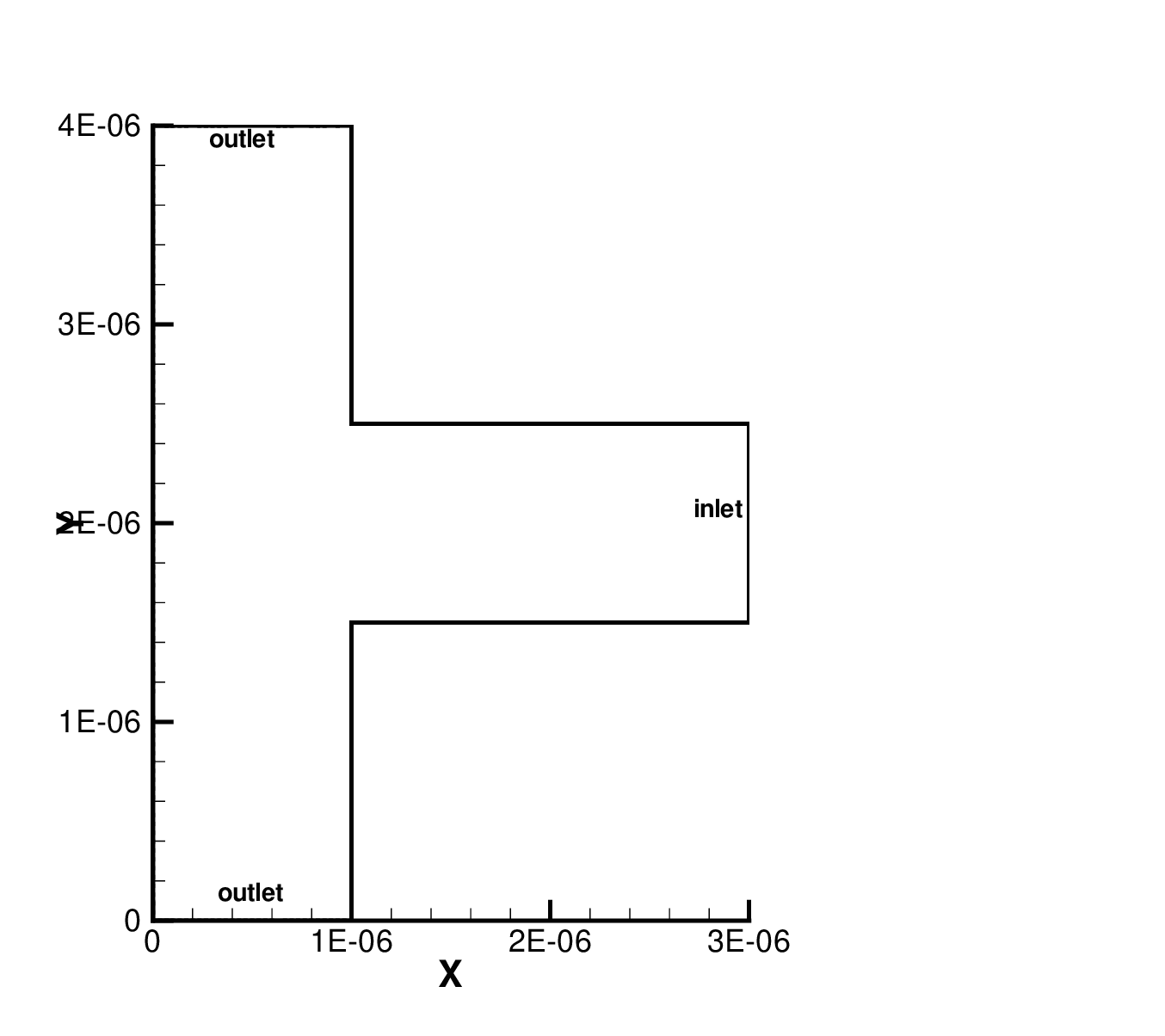}}
	\subfigure[Grid]{\includegraphics[width=0.4\linewidth]{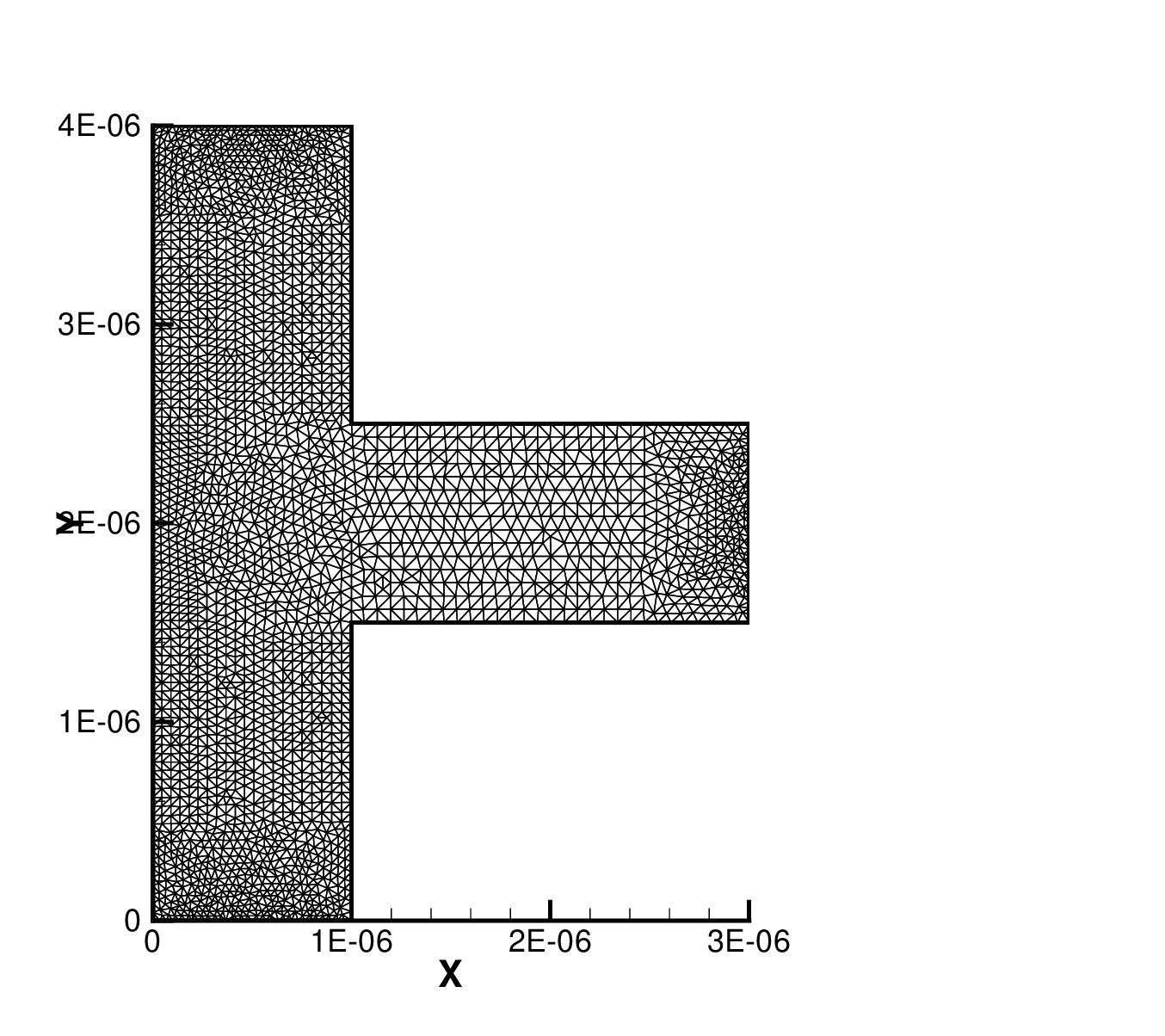}}
	\caption{Geometry model for the T-junction microchannel.}
	\label{fig:E1}
\end{figure}

In this subsection, we provide the numerical results for the electro-osmotic flow in a T-junction microchannel, the geometry of the model is given in Figure \ref{fig:E1}. The boundary condition on the inlet boundary is $u_1=-1,u_2=0$ and $c_i=1, i= 1, 2, 3$. We choose $\alpha=4\times 10^4$ in initial condition given as follows
\begin{align*}
c_1(x,0)&=\frac{1}{2} b_1(\gamma +1 -(\gamma-1) erf(\alpha x)),\\
c_3(x,0)&= b_3(1+erf(\alpha x)),\\
c_2(x,0)&=-\frac{z_1}{z_2} c_1(x,0)-\frac{z_3}{z_2} c_3(x,0),
\end{align*}
where $b_1=100, b_3=0.1$, $\gamma=50$. The slip boundary condition for the fluid is given by
\begin{align*}
u=-\xi \nabla \phi,
\end{align*}
where $\xi =1\times 10^{-6}$. The times step is $\tau = 1e-7$ and $T=6e-6$.
The finite element spaces are Mini finite element spaces for the fluid, $P1b$ (piecewise $1$ order polynomial with bubble function) element space for the molar concentrations and the electric potential. The parameters are given as follows.

\begin{center}
	\begin{table}[!ht]
		\caption{Parameters of the microfluidic model}
		\begin{center}
			\begin{tabular}{|l|l|l|l|}
				\hline
				mobility & $\nu_1=5\times 10^{-8}$ & $\nu_2=3\times 10^{-7}$ &  $\nu_3=3\times 10^{-8}$  \\
				\hline
				diffusivity & $d_1=2\times 10^{-10}$ & $d_2=3\times 10^{-10}$ & $d_3=2\times 10^{-10}$ \\
				\hline
				valence number & $z_1=1.0$ & $z_2=-1.0 $ & $z_3=-2.0$ \\
				\hline
			\end{tabular}
		\end{center}
	\end{table}
\end{center}

\begin{figure}
	\centering
	\subfigure[Contours of $c_1$]{	\includegraphics[width=0.44\linewidth]{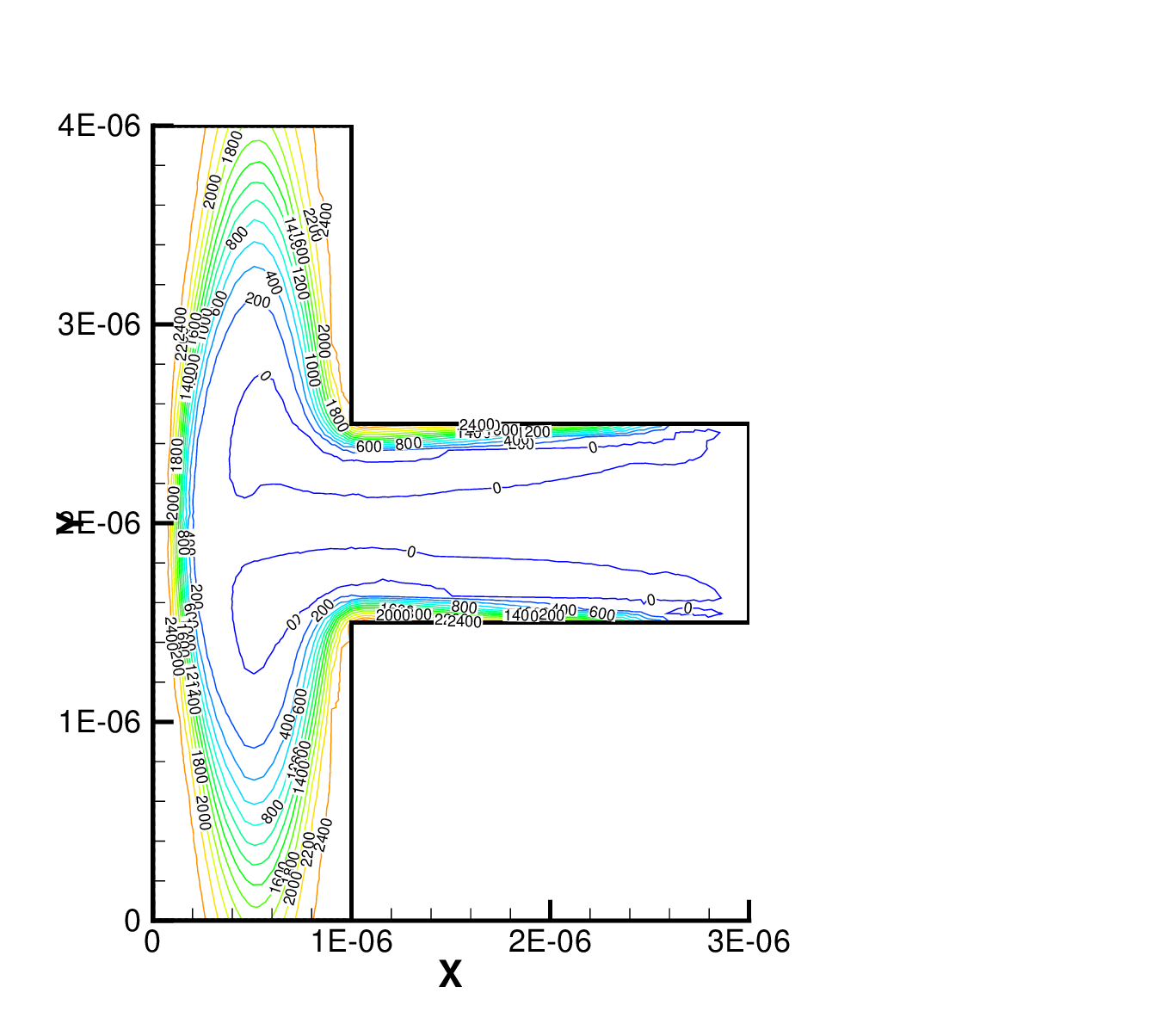}}
	\subfigure[Contours of $c_2$]{	\includegraphics[width=0.44\linewidth]{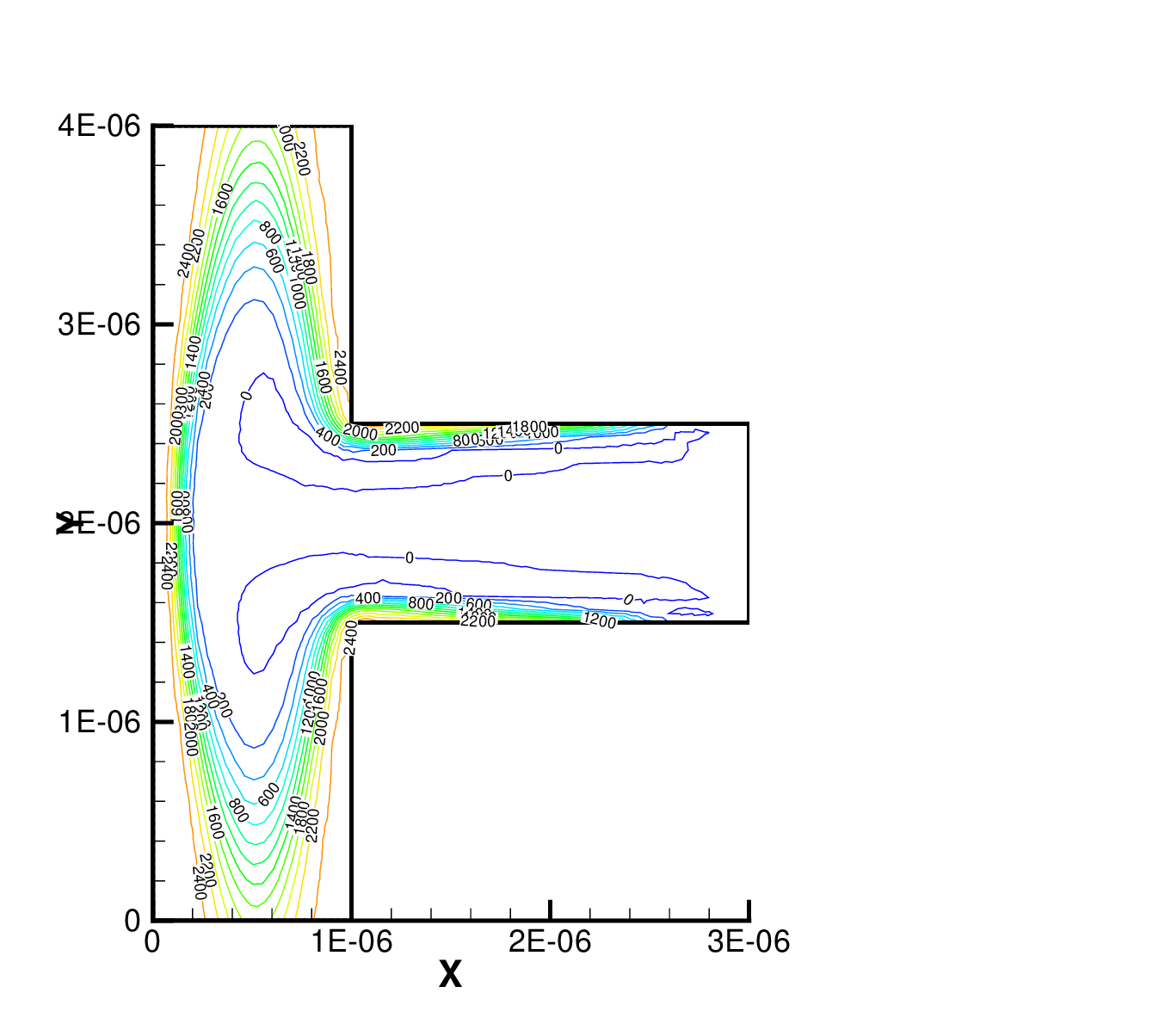}}\\
	\subfigure[Contours of $c_3$]{\includegraphics[width=0.44\linewidth]{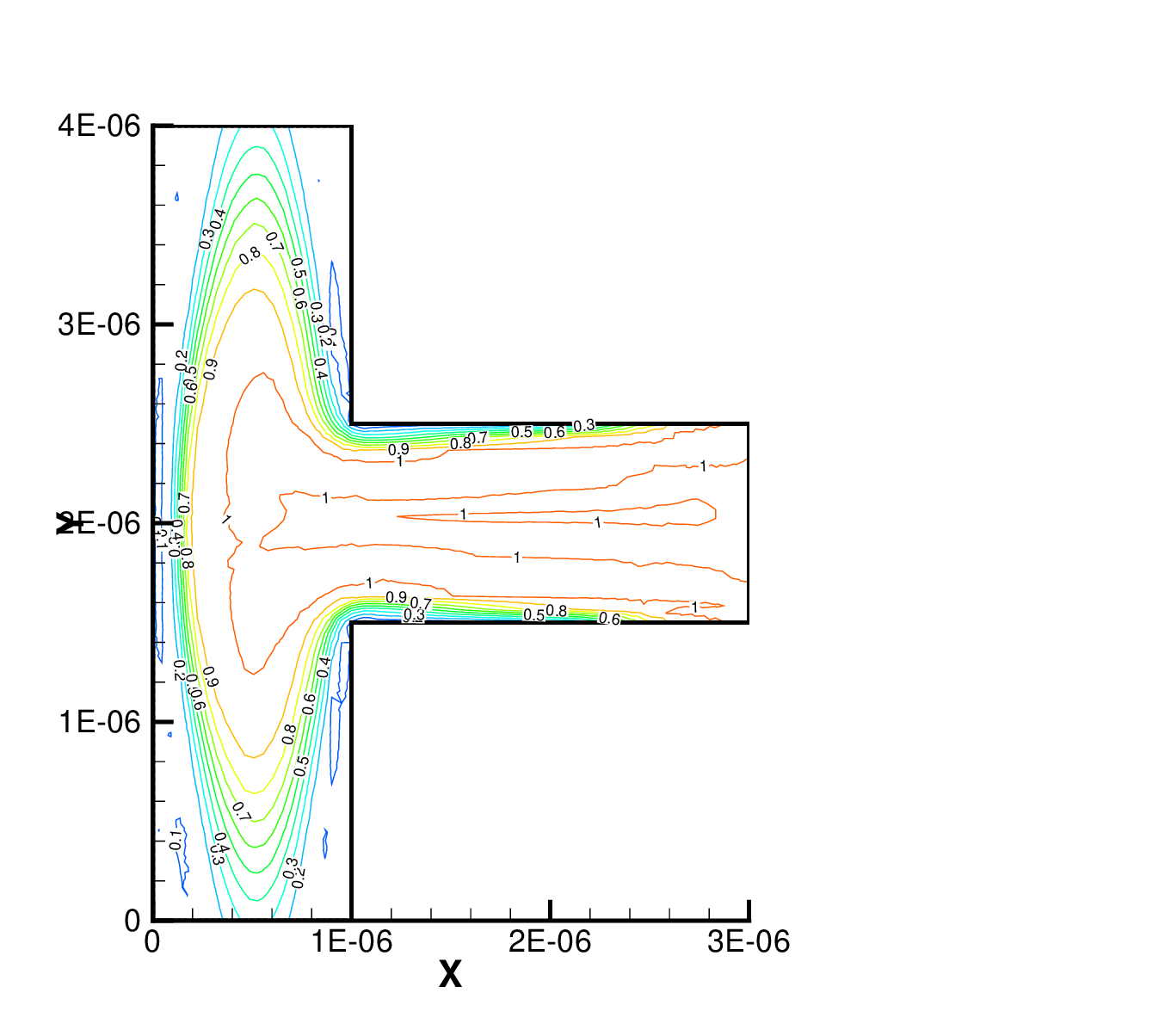}}\subfigure[Streamline]{\includegraphics[width=0.44\linewidth]{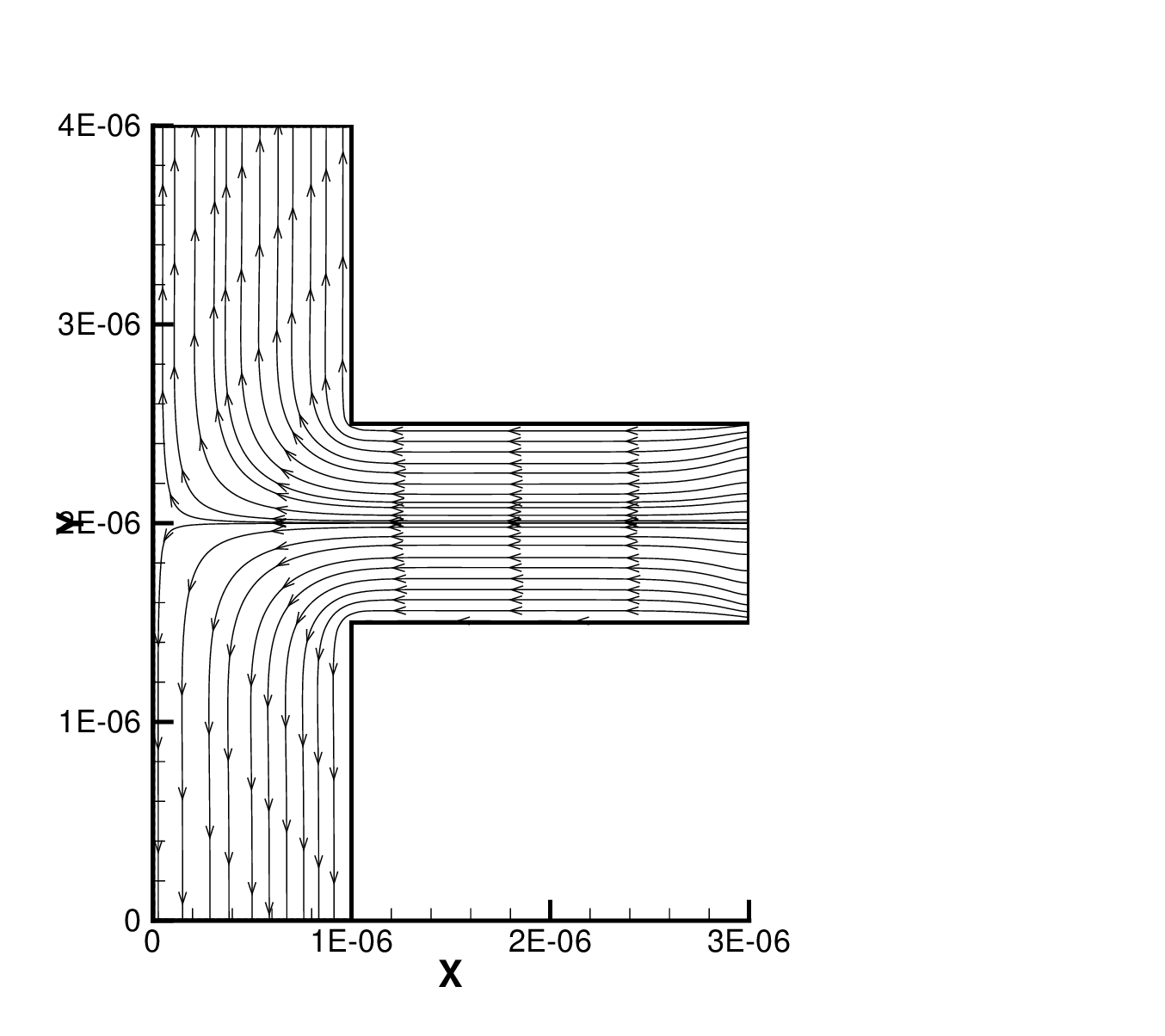}}
	\caption{Numerical results of the T-junction model with $\nu_1=\nu_2=\nu_3=1.0$ and $T=4e-6$.}
	\label{fig:e1-1}
\end{figure}

\begin{figure}
	\centering
	\subfigure[Contours of $c_1$]{	\includegraphics[width=0.44\linewidth]{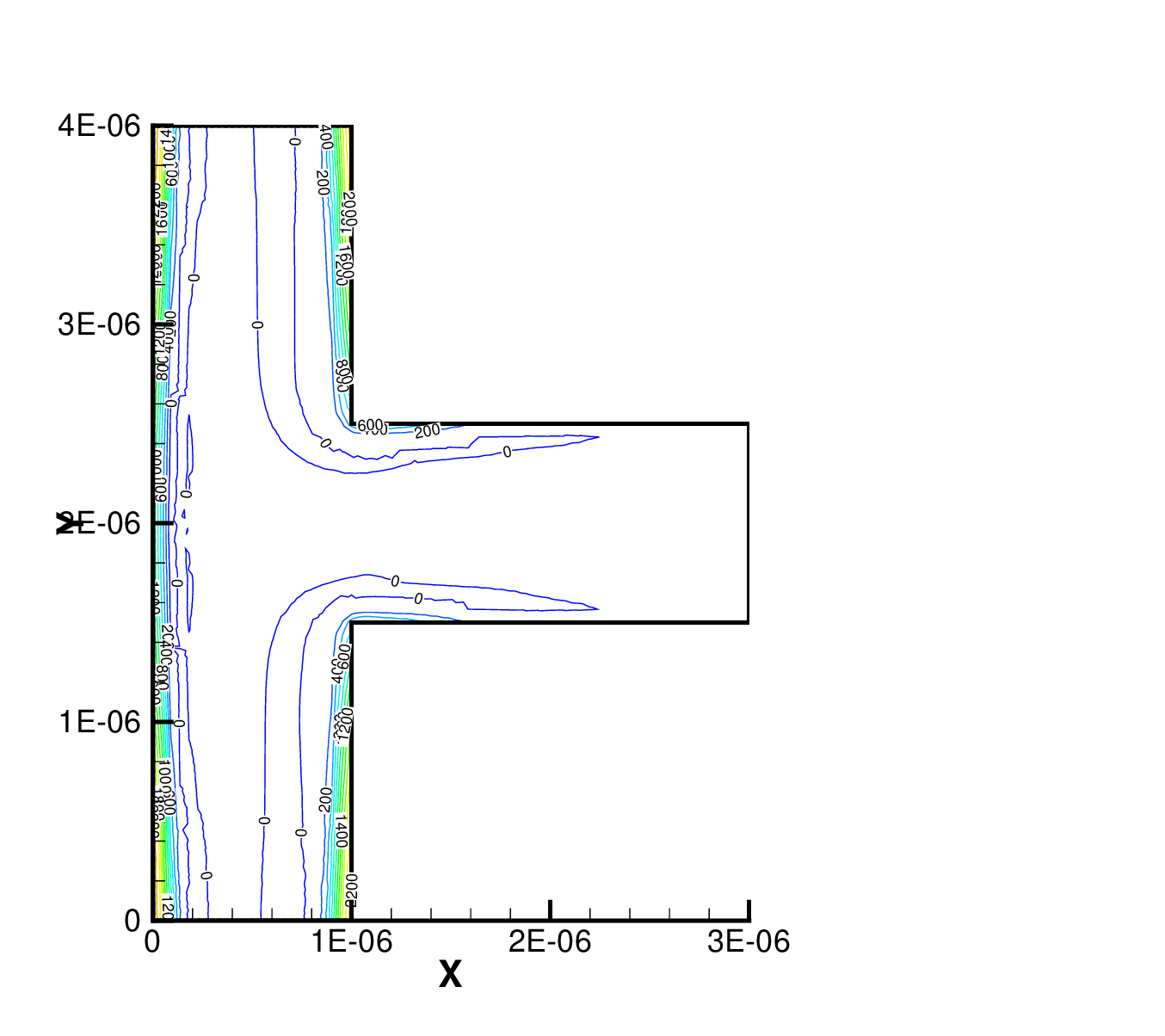}}
	\subfigure[Contours of $c_2$]{	\includegraphics[width=0.44\linewidth]{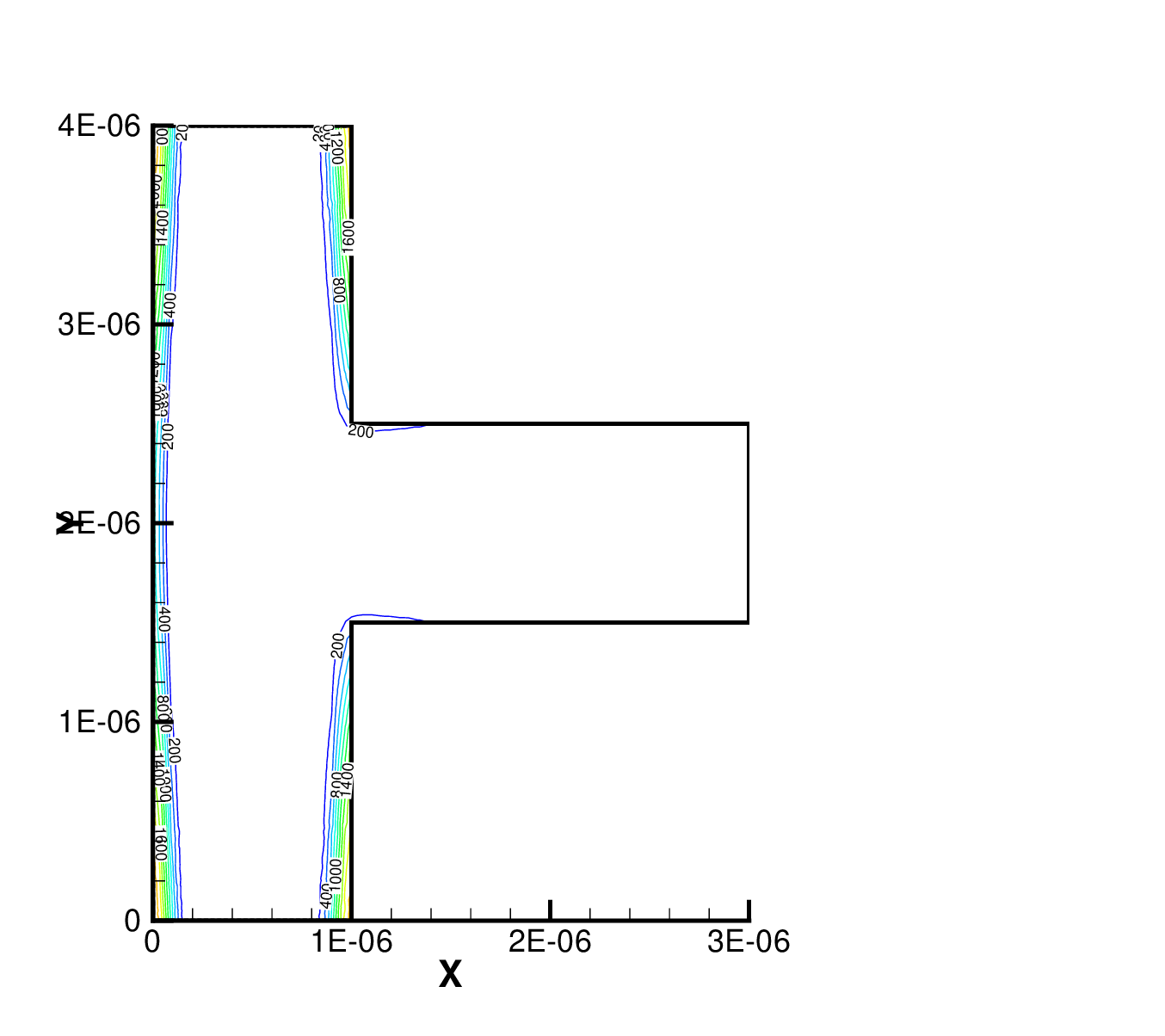}}\\
	\subfigure[Contours of $c_3$]{\includegraphics[width=0.44\linewidth]{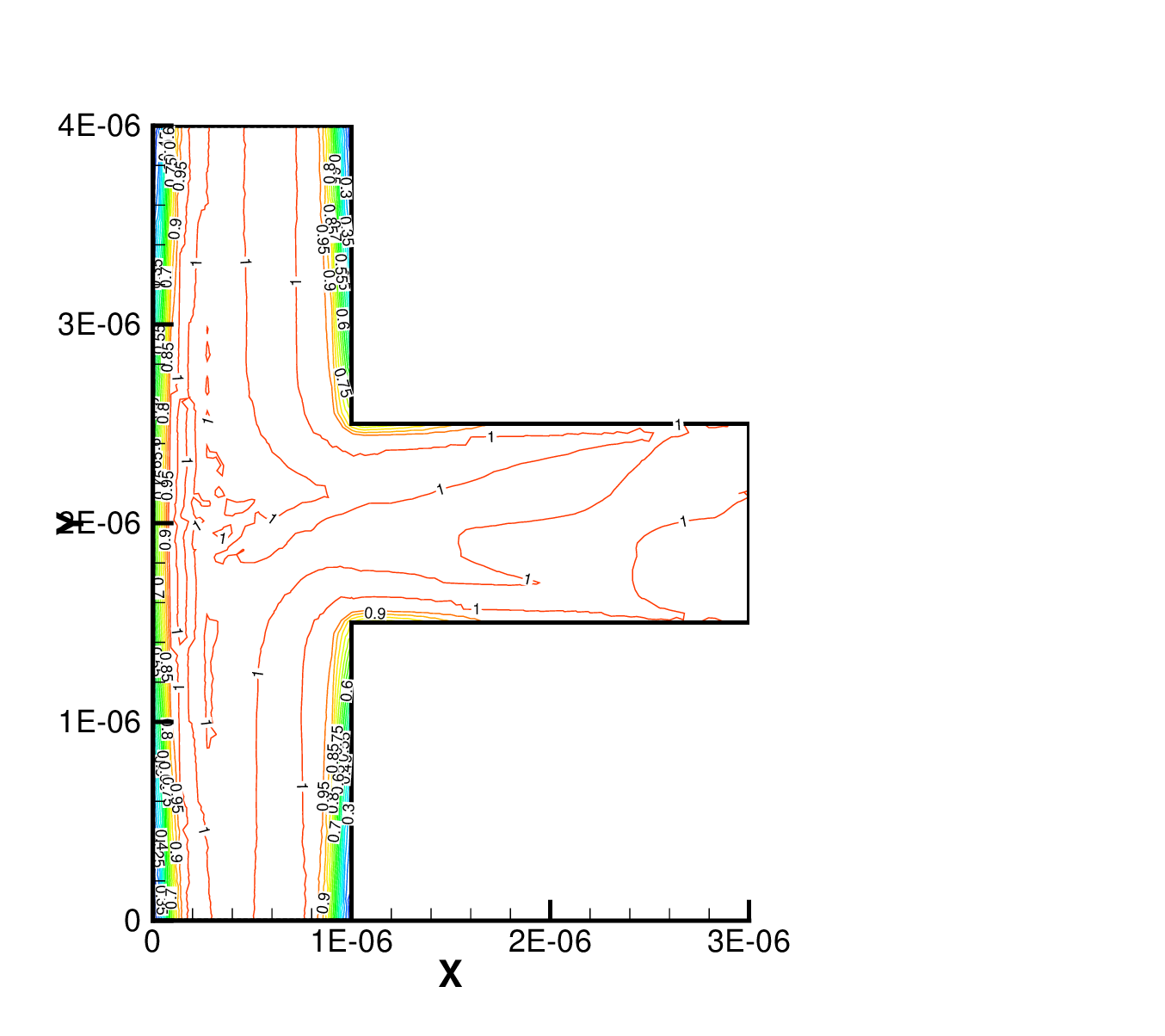}}\subfigure[Streamline]{\includegraphics[width=0.44\linewidth]{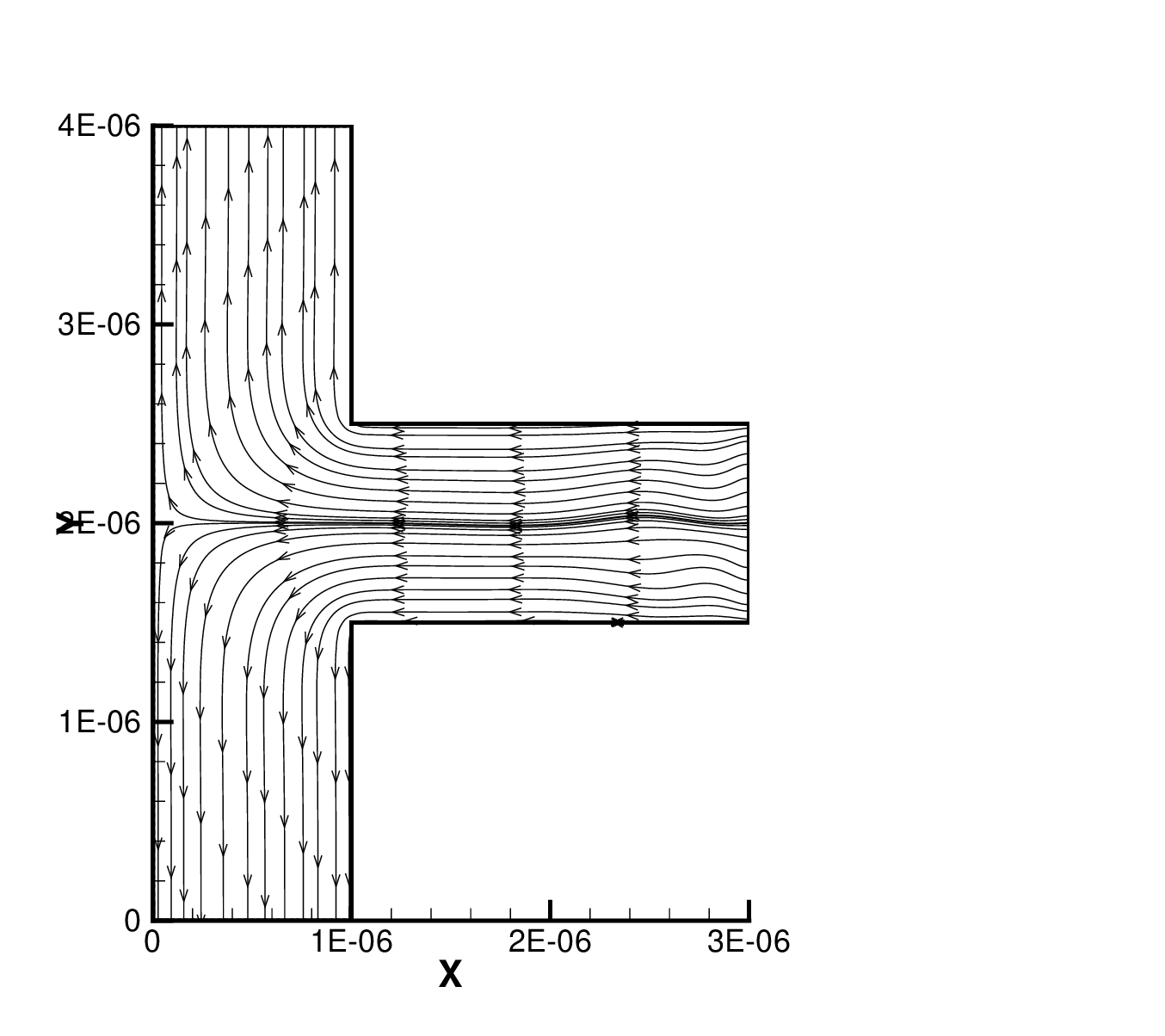}}
	\caption{Numerical results of the T-junction model with $\nu_1=\nu_2=\nu_3=1.0$ and $T=1.36e-5$.}
	\label{fig:e1-2}
\end{figure}

\begin{figure}
	\centering
	\subfigure[Contours of $c_1$]{	\includegraphics[width=0.44\linewidth]{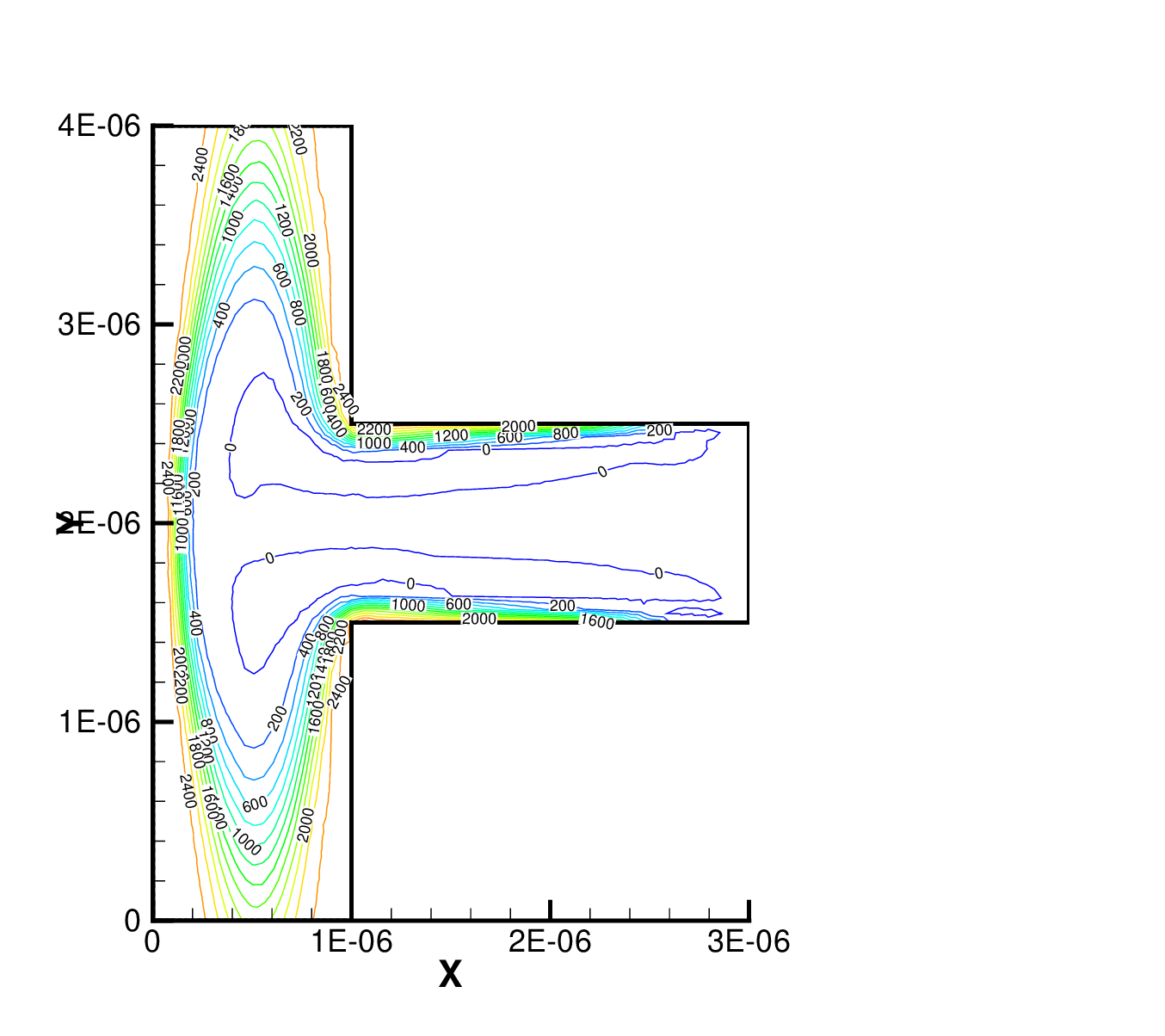}}
	\subfigure[Contours of $c_2$]{	\includegraphics[width=0.44\linewidth]{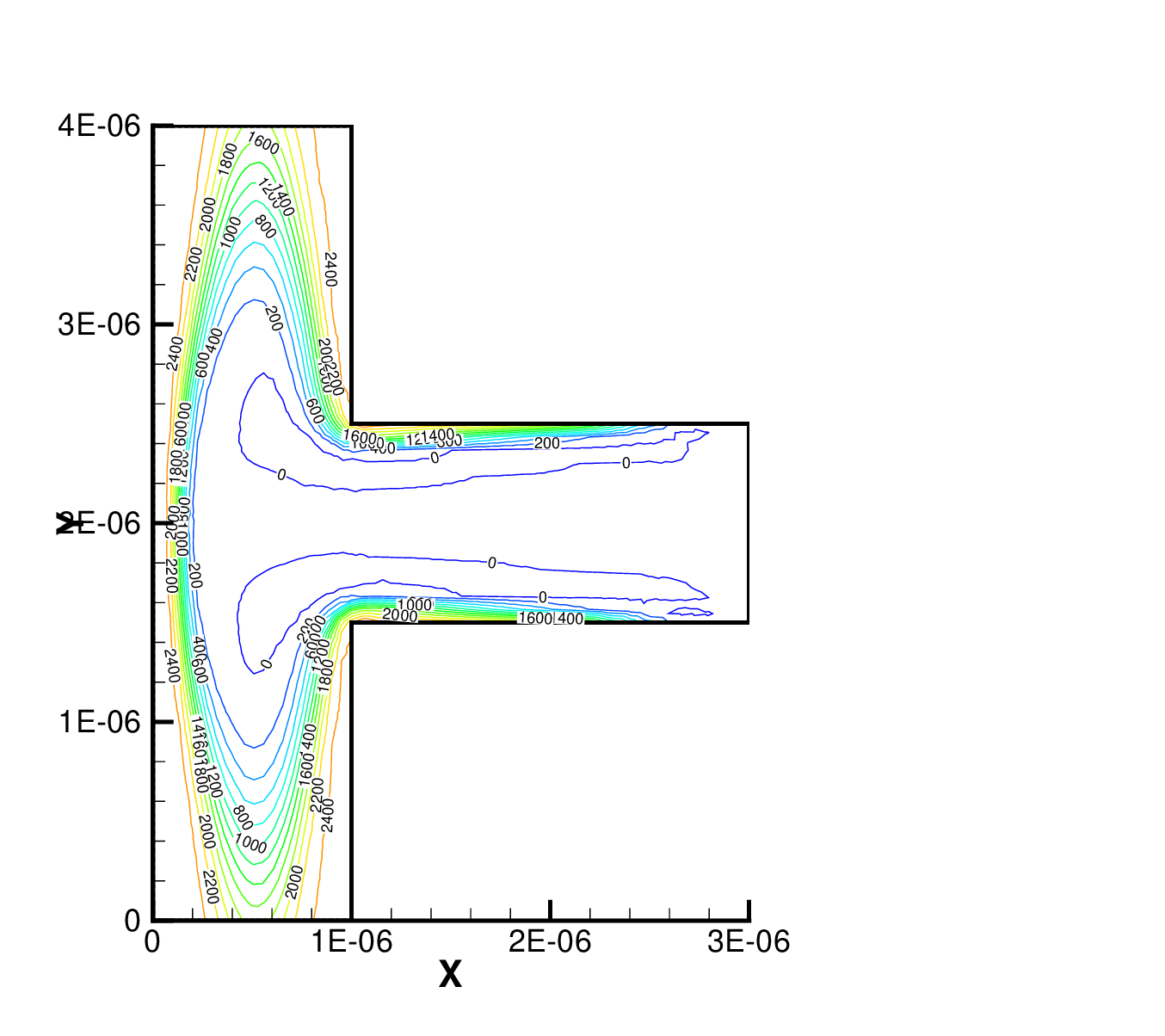}}\\
	\subfigure[Contours of $c_3$]{\includegraphics[width=0.44\linewidth]{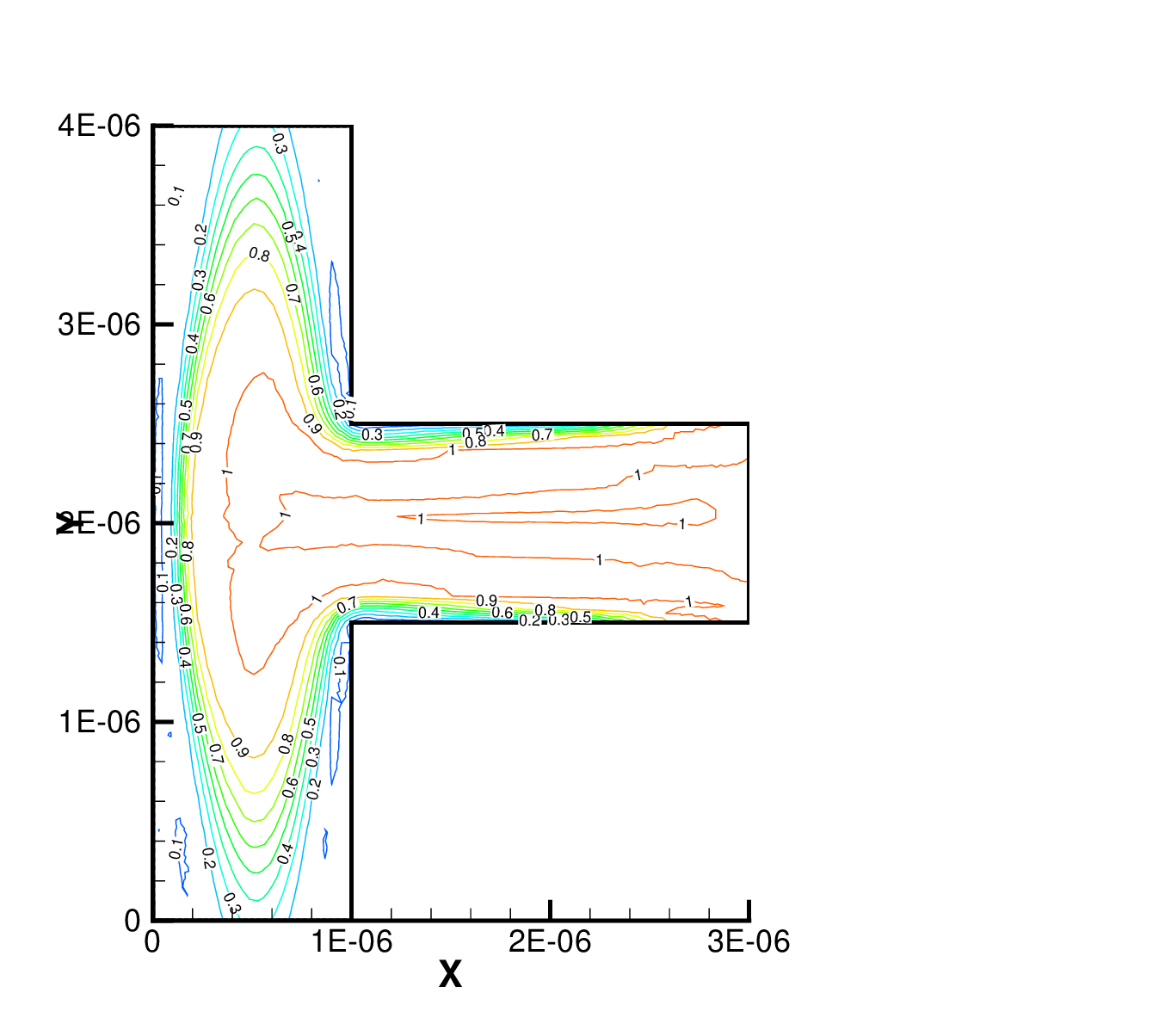}}\subfigure[Streamline]{\includegraphics[width=0.44\linewidth]{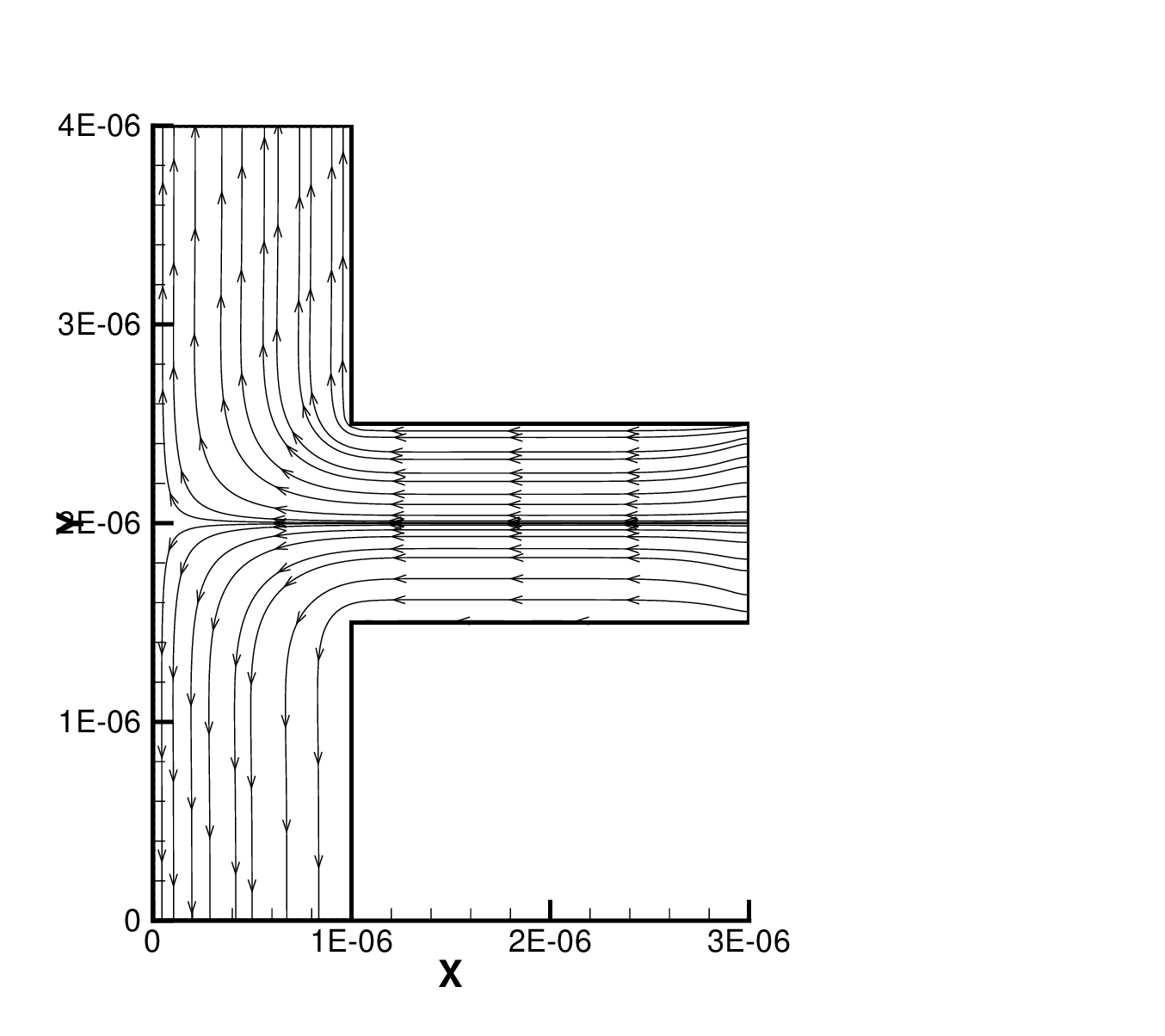}}
	\caption{Numerical results of the T-junction model with $\nu=0.1$ and $T=4e-6$.}
	\label{fig:e1-3}
\end{figure}

\begin{figure}
	\centering
	\subfigure[Contours of $c_1$]{	\includegraphics[width=0.44\linewidth]{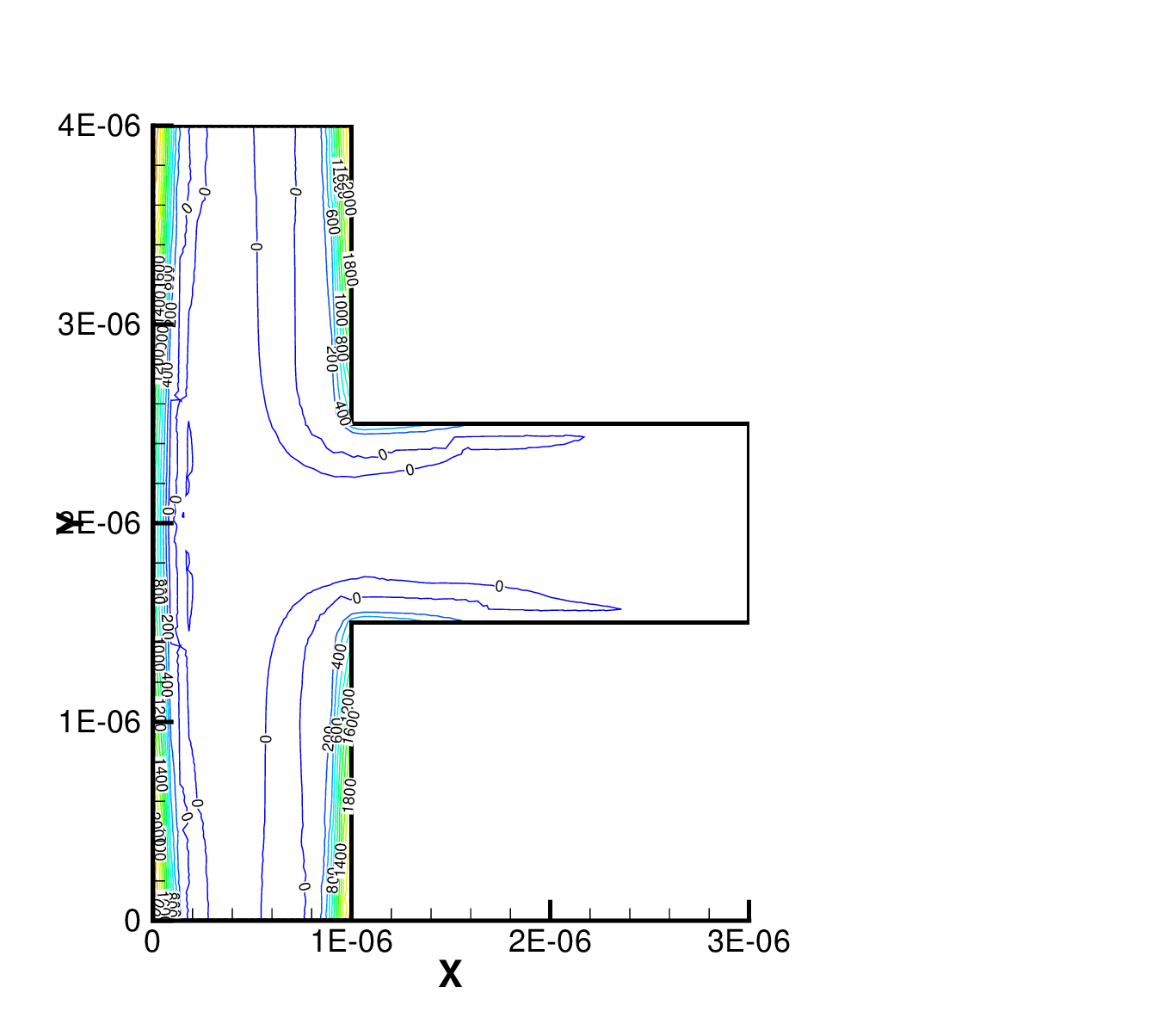}}
	\subfigure[Contours of $c_2$]{	\includegraphics[width=0.44\linewidth]{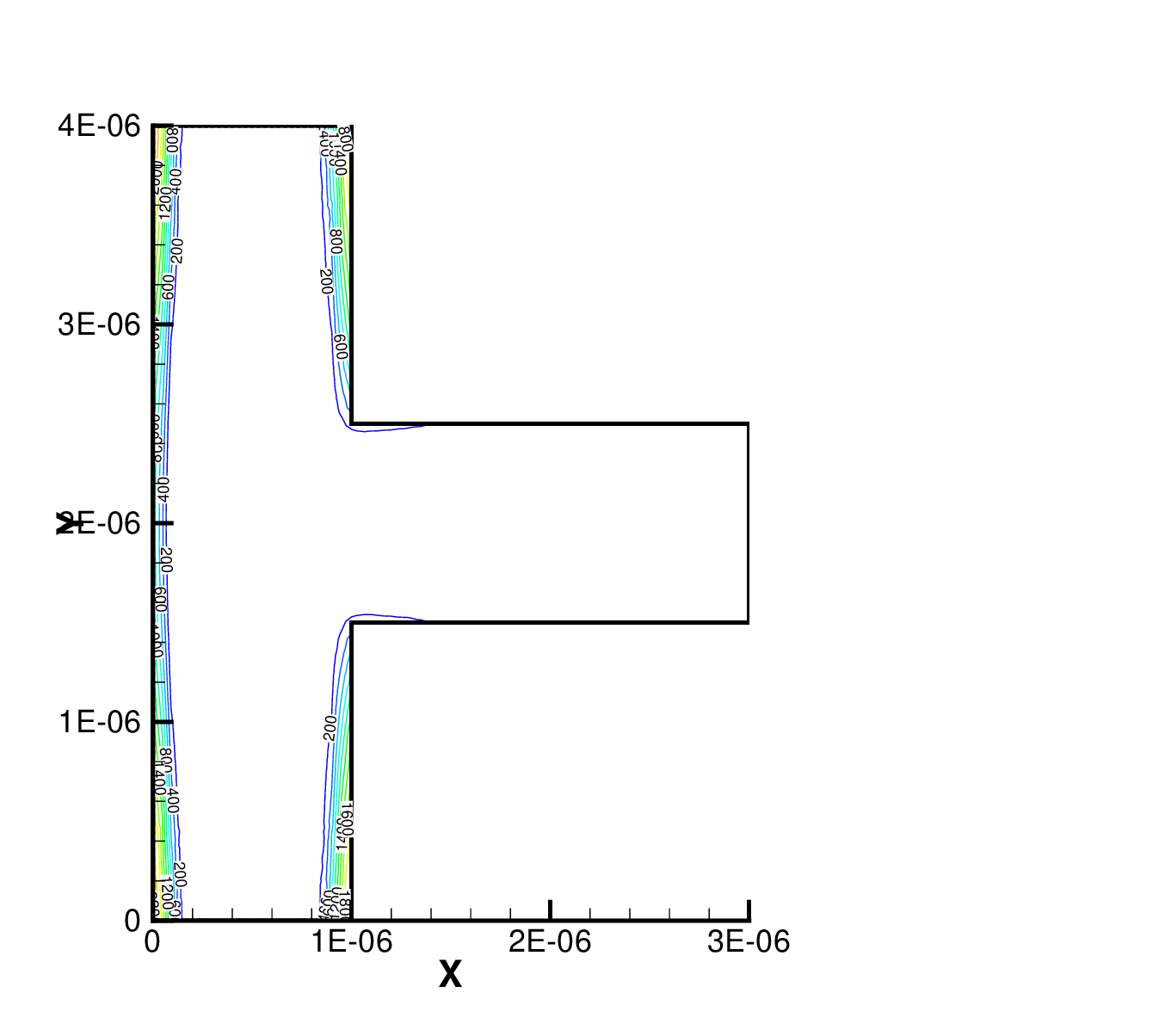}}\\
	\subfigure[Contours of $c_3$]{\includegraphics[width=0.44\linewidth]{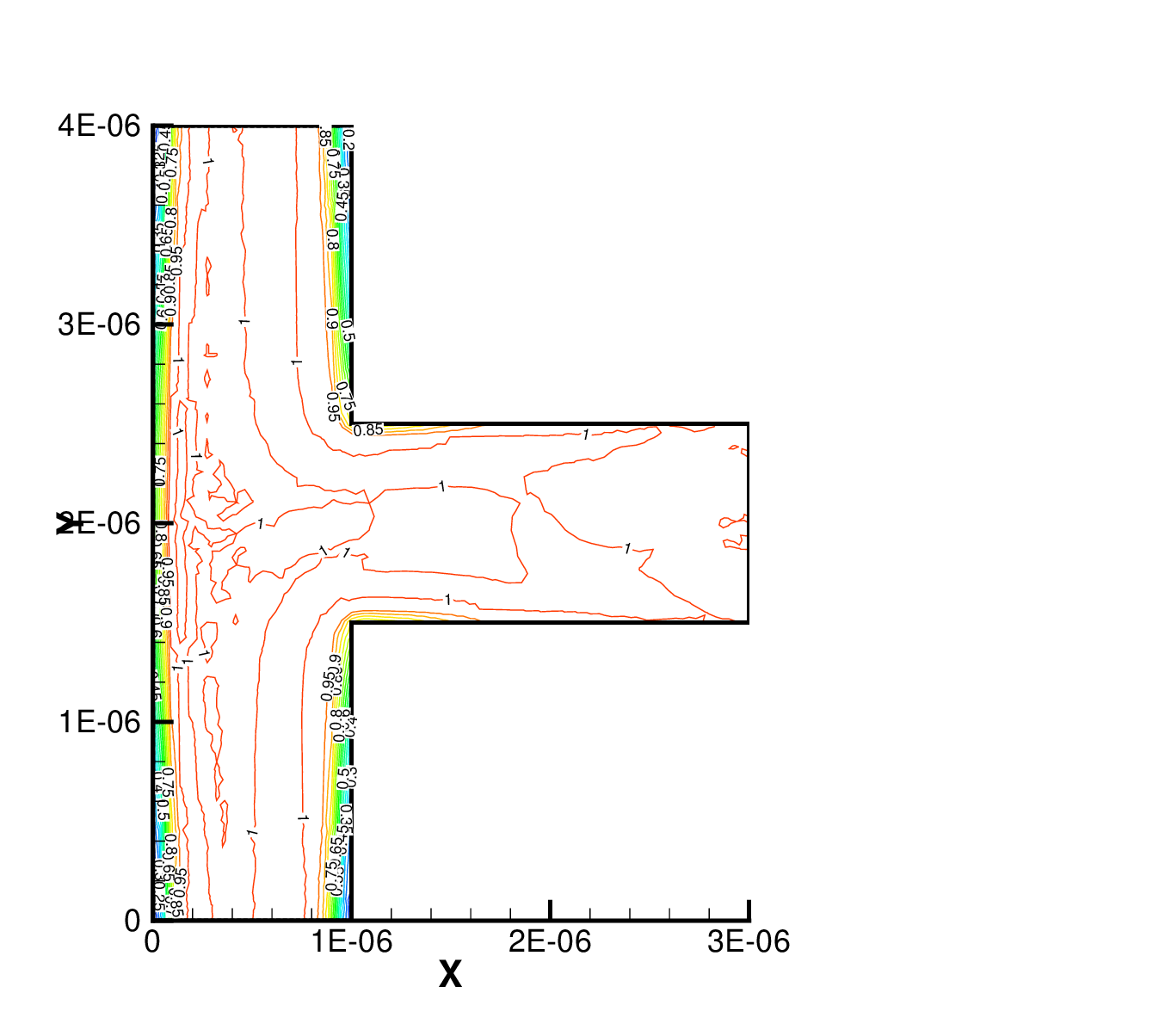}}\subfigure[Streamline]{\includegraphics[width=0.44\linewidth]{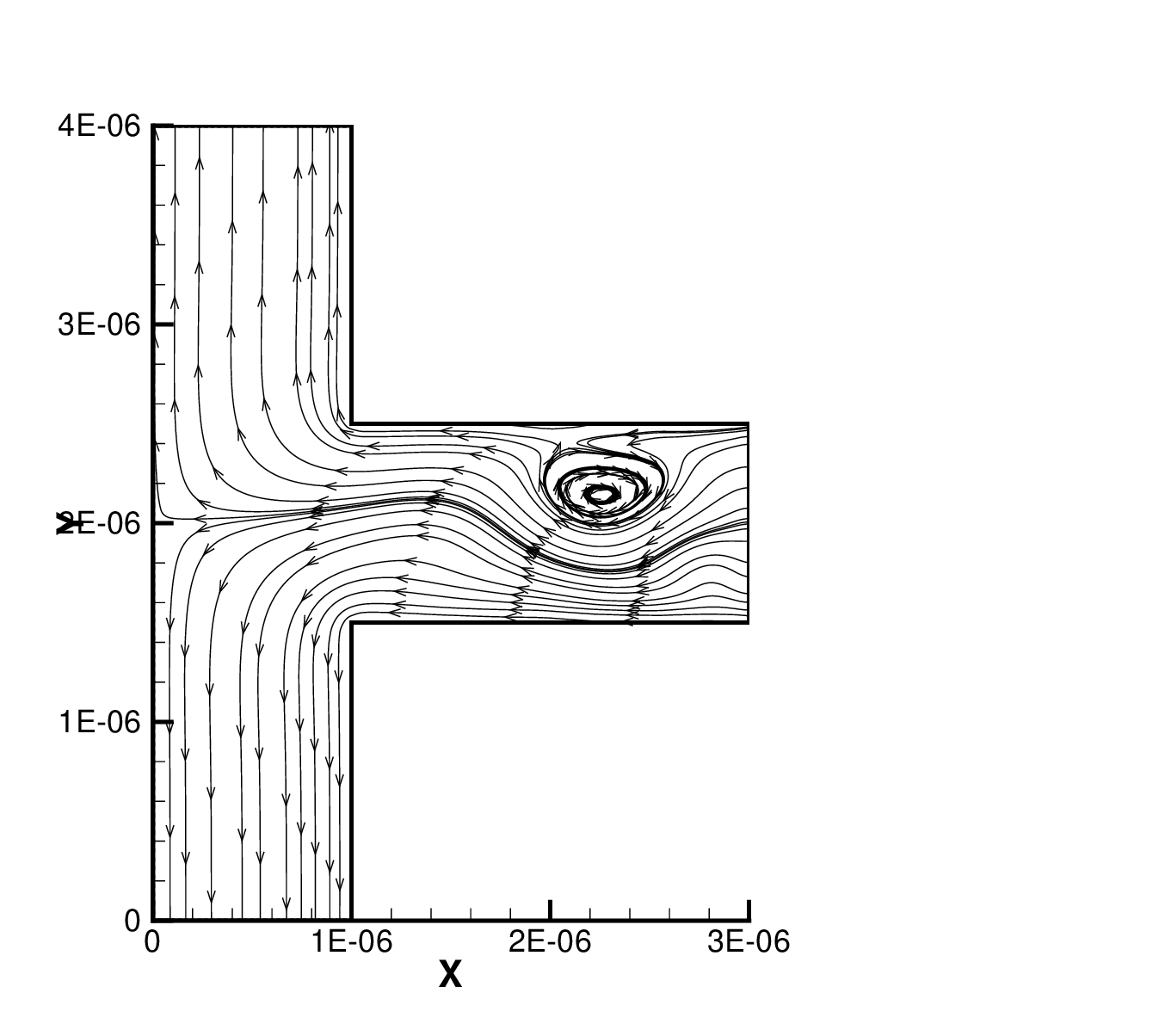}}
	\caption{Numerical results of the T-junction model with $\nu=0.1$ and $T=1.36e-5$.}
	\label{fig:e1-4}
\end{figure}

In Figure \ref{fig:e1-1}, we give the numerical results with $\nu_1=\nu_2=\nu_3=1.0$ at $T=4e-6$. In Figure \ref{fig:e1-2}, the numerical results with $\nu_1=\nu_2=\nu_3=1.0$ at $T=1.36e-5$ are presented. The flow is stable, and the molar concentrations change as time evolves. Then, numerical results with $\nu=0.1$ is shown in Figure  \ref{fig:e1-3} and \ref{fig:e1-4}. We can see that there are two vortexes near the inlet at $T=1.36e-5$.


\subsection{The effect of the roughness in microchannels}

\begin{figure}
	\centering
	\subfigure[Geometry model]{\includegraphics[width=0.44\linewidth]{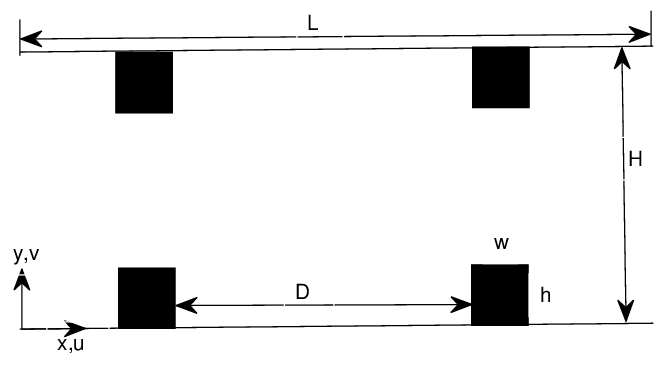}}
	\subfigure[Grid]{\includegraphics[width=0.44\linewidth]{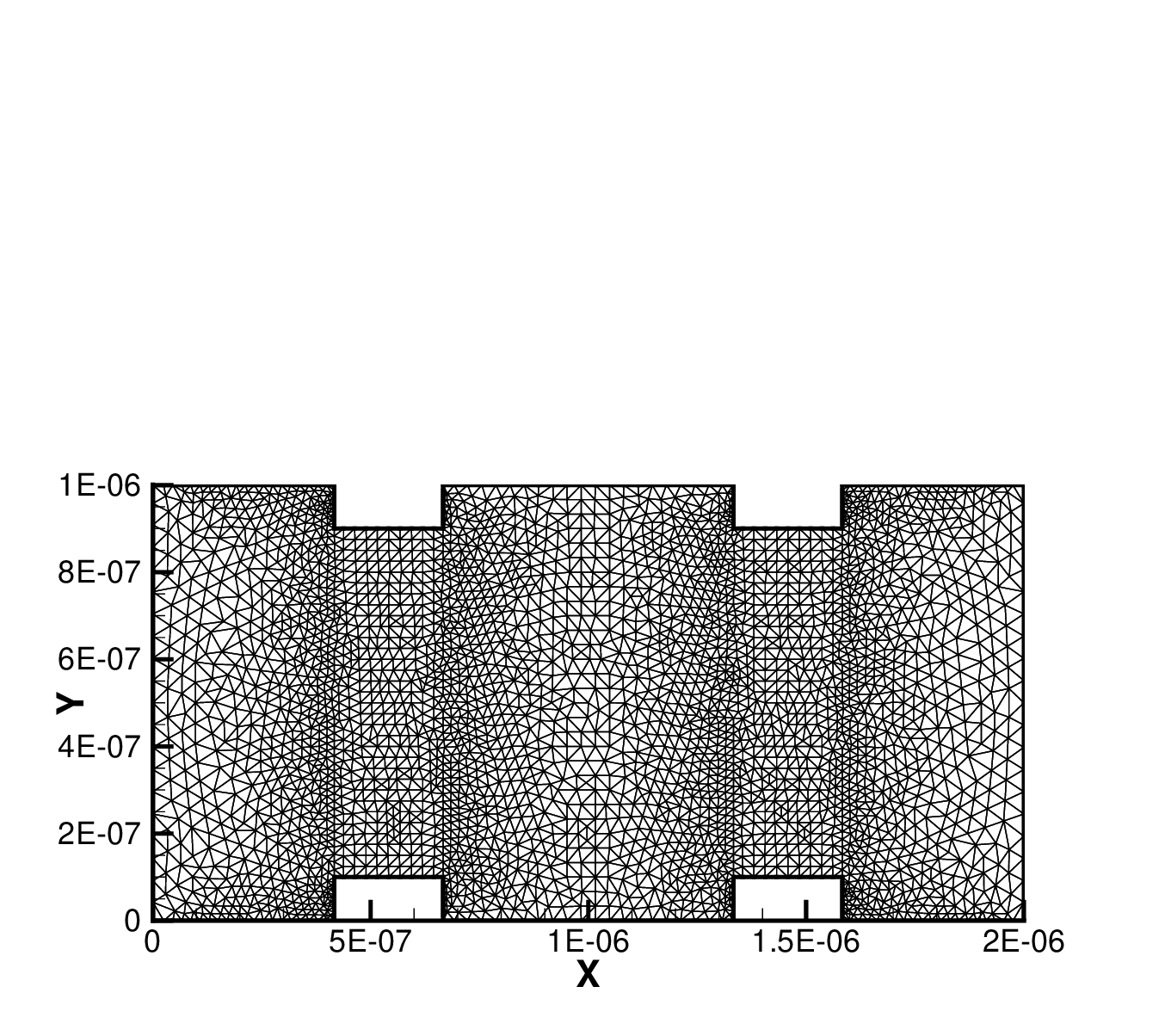}}		
	\caption{Geometry model for the microchannel with roughness.}
	\label{fig:model}
\end{figure}

An electro-osmotic flow in two-dimensional microchannel is numerically studied. The geometry of the model is given in Figure \ref{fig:model}.
The channel width is $H=1e-6$, the channel length is $L=2H$. The width of each roughness is $w=H/4$ and the space interval is $D=L/3$.
The roughness dimensions are $w$ and $h$ in $x$ and $h$ directions, respectively. The roughness is uniformly positioned in the channels with an interval space $D$.
The boundary condition on the top and bottom boundaries is the solid condition.
The boundary condition on the right boundary is the outlet condition. 
As same as the first example, we choose the same parameters and initial condition. The boundary condition for the fluid is slip boundary condition is chosen as $\xi =1\times 10^{-6}$. The times step is $\tau = 1e-7$ and $T=2e-5$.
The boundary conditions on the inlet boundary are given as
\begin{align*}
u_1=1e-2,u_2=0, c_i=1.0, i=1,2,3.
\end{align*}
The finite element spaces are Mini finite element spaces for the fluid, $P1b$ element space for the molar concentrations and the electric potential.

\begin{figure}
	\centering
	\subfigure[Contours of $c_1$]{	\includegraphics[width=0.44\linewidth]{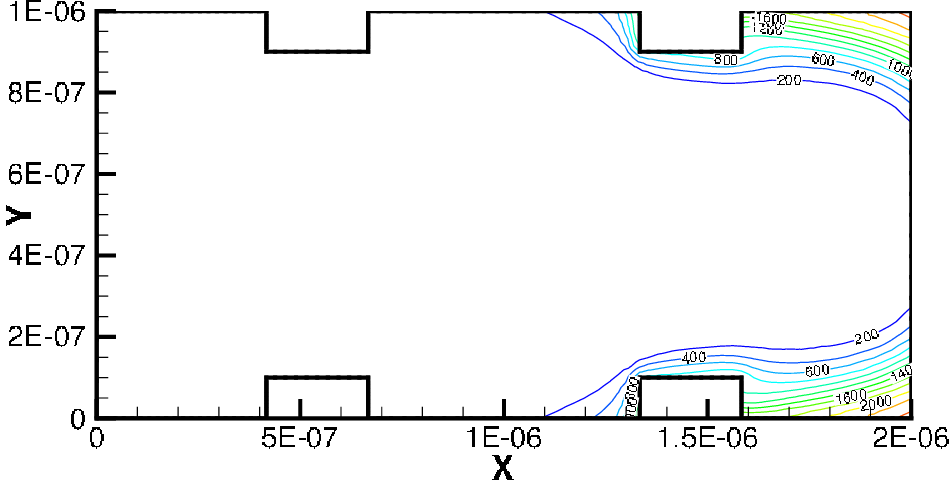}}
	\subfigure[Contours of $c_2$]{	\includegraphics[width=0.44\linewidth]{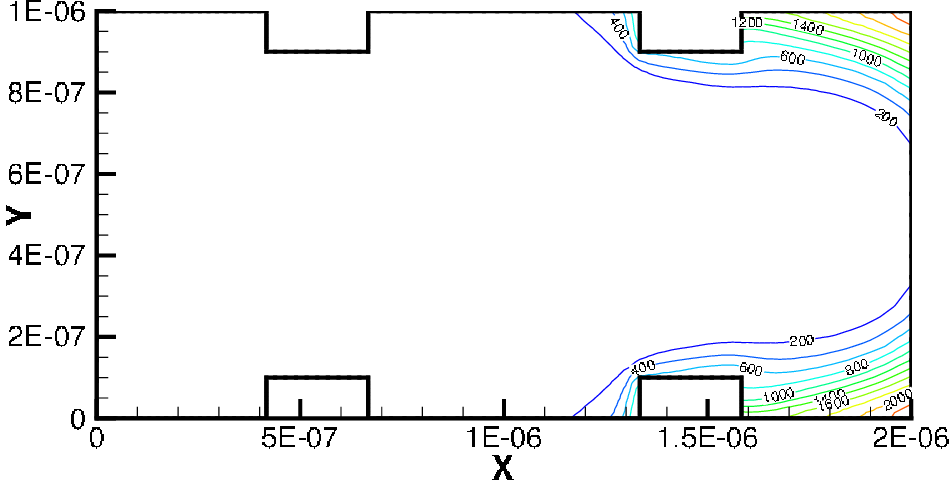}}\\
	\subfigure[Contours of $c_3$]{\includegraphics[width=0.44\linewidth]{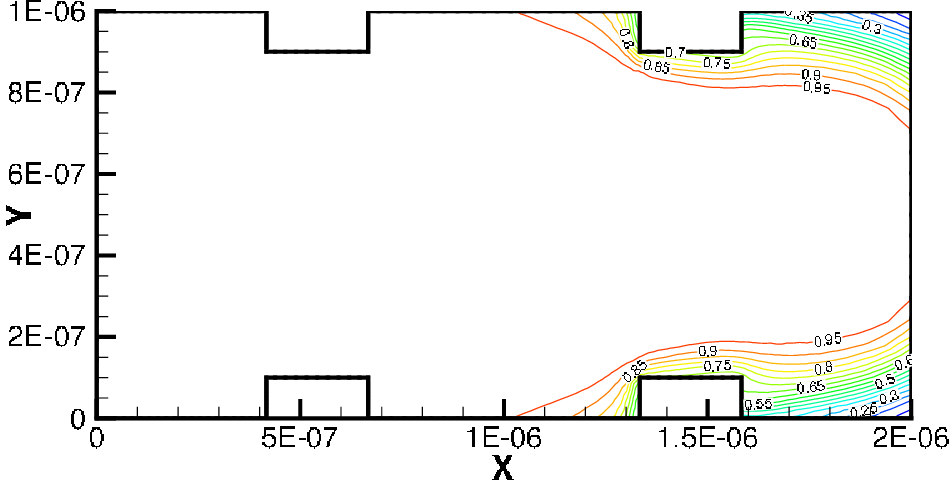}}
	\subfigure[Streamline]{\includegraphics[width=0.44\linewidth]{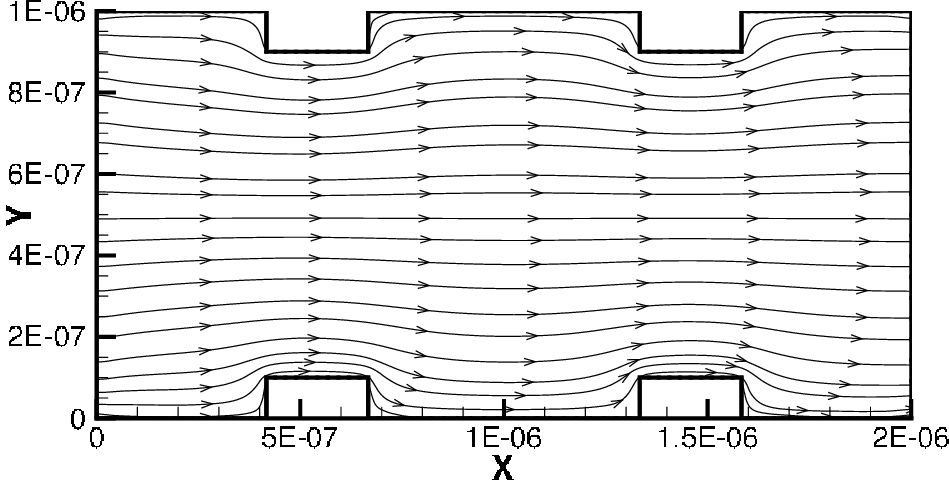}}
	\caption{Geometry model for the microchannel with roughness $h=0.1H$.}
	\label{fig:01h}
\end{figure}

\begin{figure}
	\centering
	\subfigure[Contours of $c_1$]{	\includegraphics[width=0.44\linewidth]{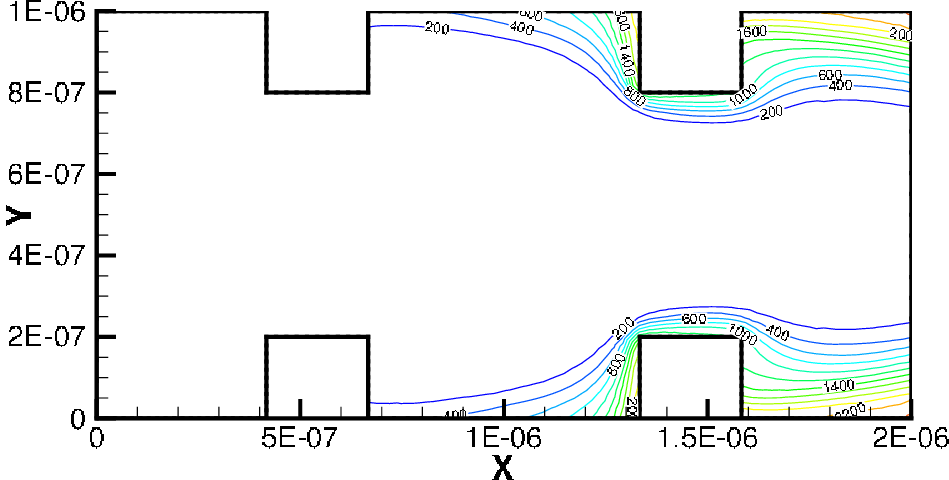}}
	\subfigure[Contours of $c_2$]{	\includegraphics[width=0.44\linewidth]{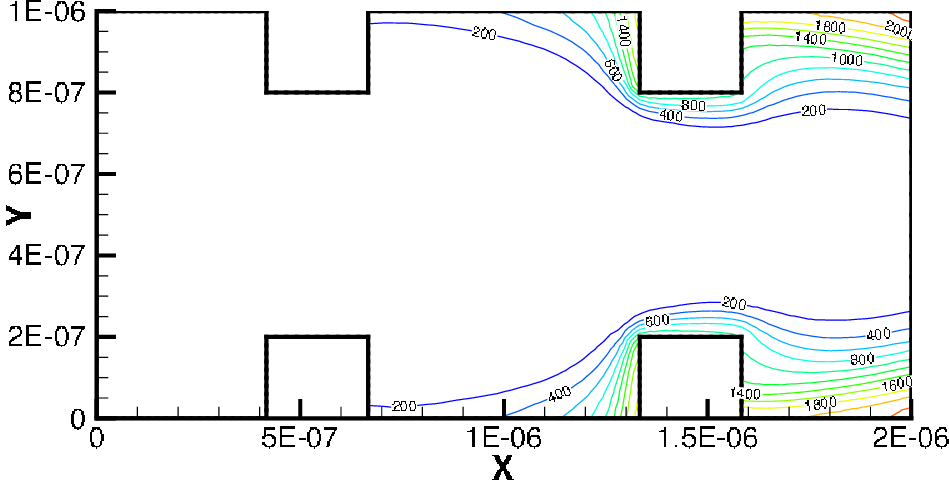}}\\
	\subfigure[Contours of $c_3$]{\includegraphics[width=0.44\linewidth]{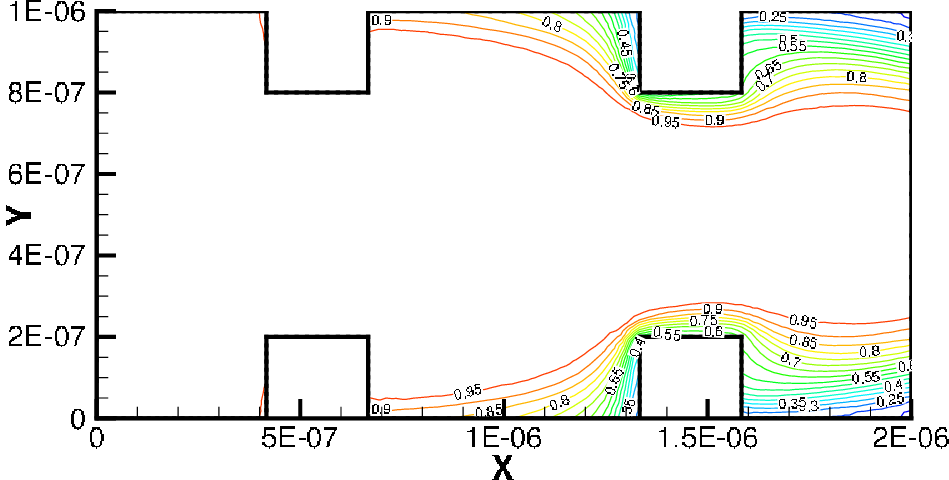}}
	\subfigure[Streamline]{\includegraphics[width=0.44\linewidth]{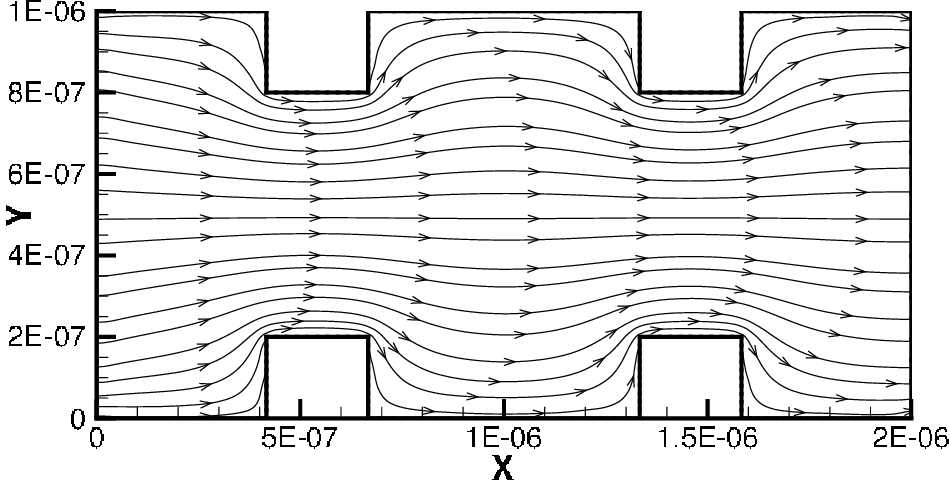}}
	\caption{Numerical results of the roughness in microchannels with $h=0.2H$.}
	\label{fig:02h}7
\end{figure}

\begin{figure}
	\centering
	\subfigure[Contours of $c_1$]{	\includegraphics[width=0.44\linewidth]{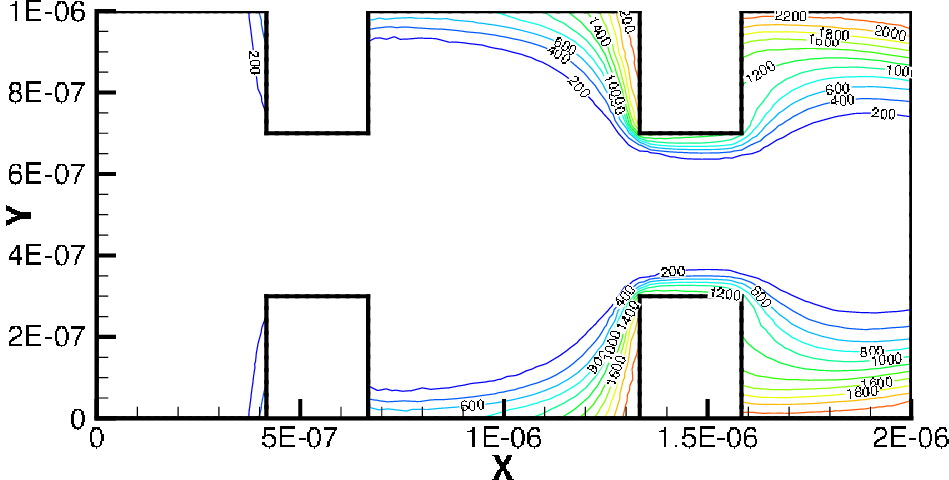}}
	\subfigure[Contours of $c_2$]{	\includegraphics[width=0.44\linewidth]{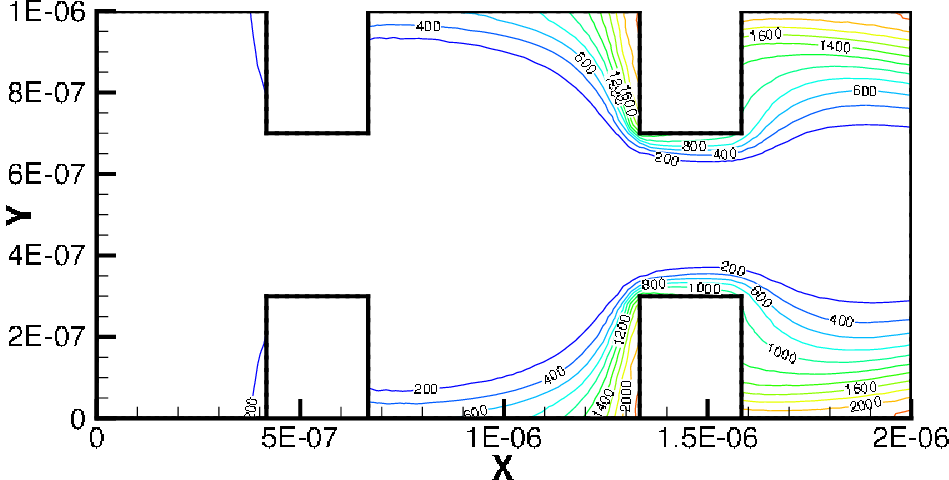}}\\
	\subfigure[Contours of $c_3$]{\includegraphics[width=0.44\linewidth]{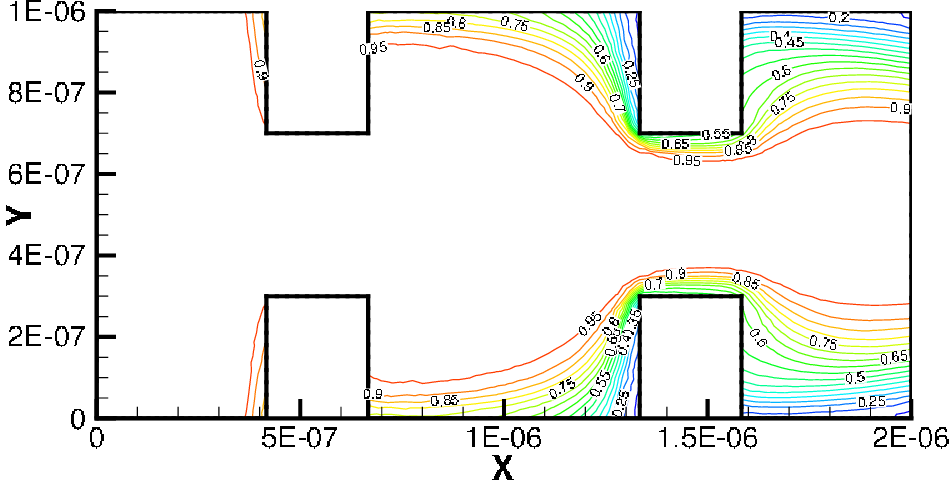}}
	\subfigure[Streamline]{\includegraphics[width=0.44\linewidth]{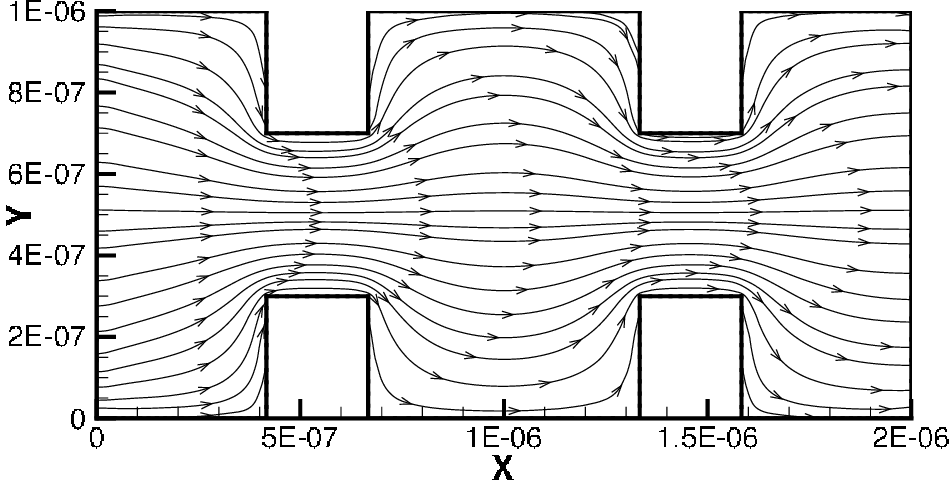}}
	\caption{Numerical results of the roughness in microchannels with $h=0.3H$.}
	\label{fig:03h}
\end{figure}

\begin{figure}
	\centering
	\subfigure[Contours of $c_1$]{	\includegraphics[width=0.44\linewidth]{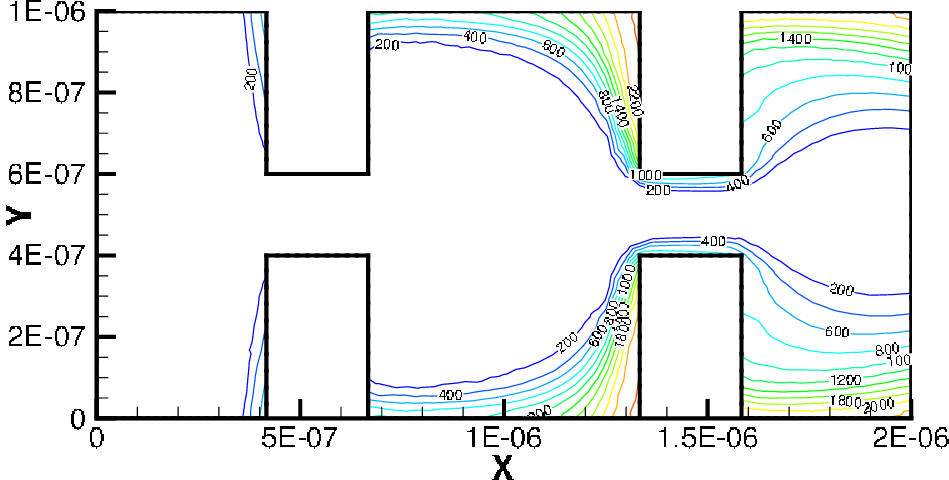}}
	\subfigure[Contours of $c_2$]{	\includegraphics[width=0.44\linewidth]{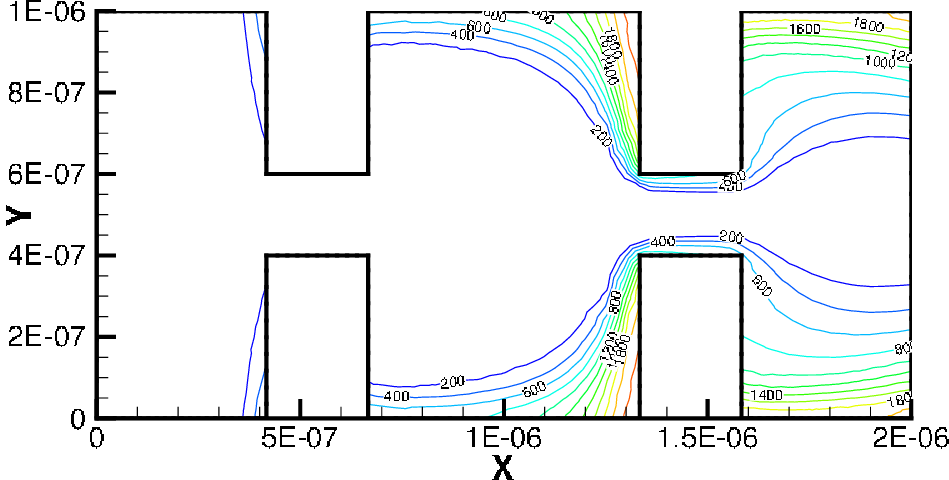}}\\
	\subfigure[Contours of $c_3$]{\includegraphics[width=0.44\linewidth]{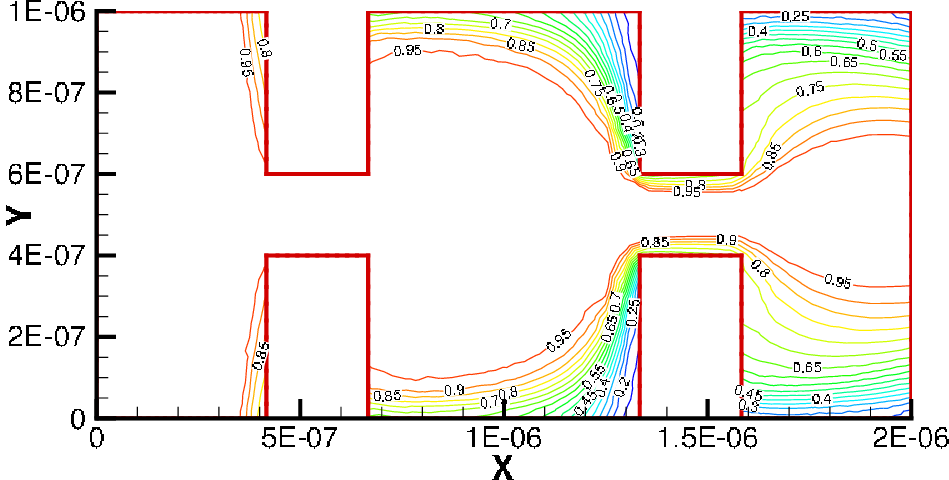}}
	\subfigure[Streamline]{\includegraphics[width=0.44\linewidth]{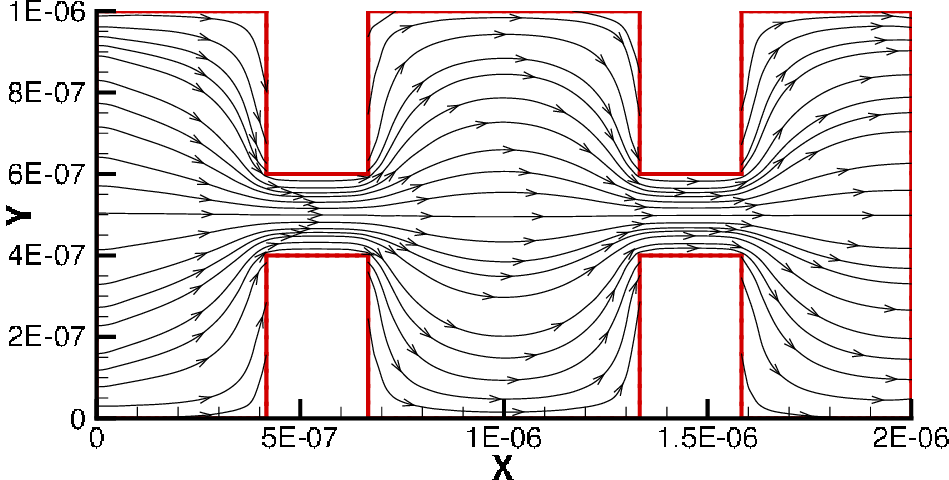}}
	\caption{Numerical results of the roughness in microchannels with $h=0.4H$.}
	\label{fig:04h}
\end{figure}

\begin{figure}
	\centering
	\subfigure[Contours of $u_1$ with $h=0.1h$]{	\includegraphics[width=0.44\linewidth]{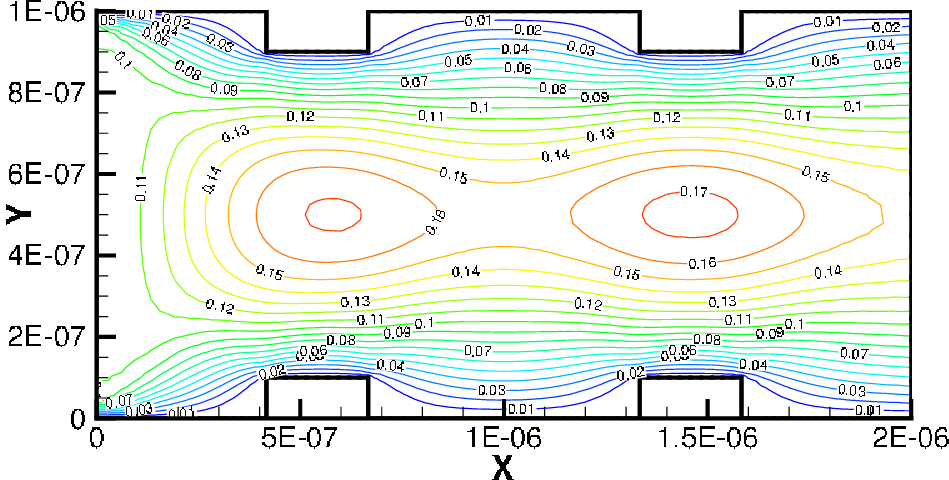}}
	\subfigure[Contours of $u_1$ with $h=0.2h$]{	\includegraphics[width=0.44\linewidth]{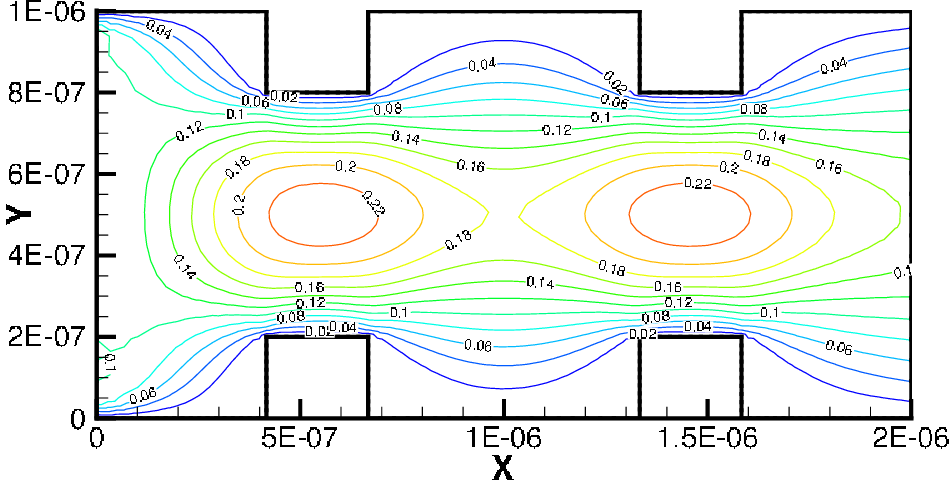}}\\
	\subfigure[Contours of $u_1$ with $h=0.3h$]{\includegraphics[width=0.44\linewidth]{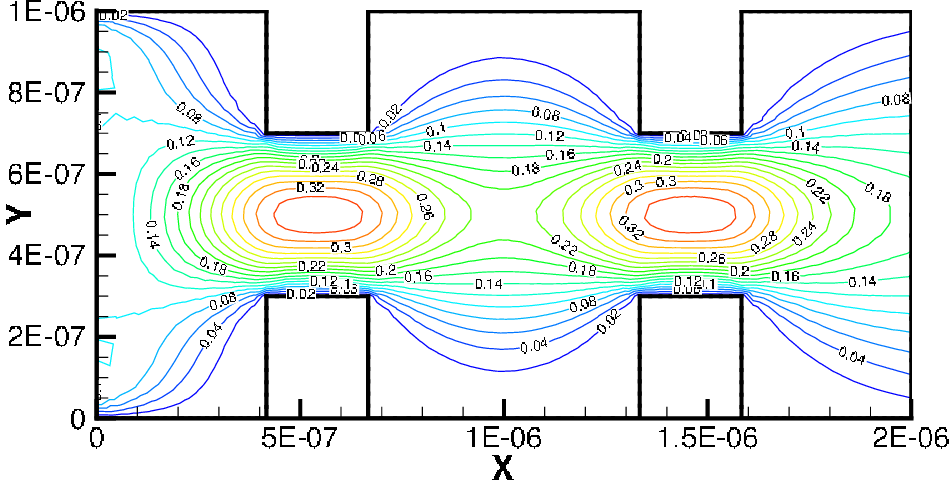}}
	\subfigure[Contours of $u_1$ with $h=0.4h$]{\includegraphics[width=0.44\linewidth]{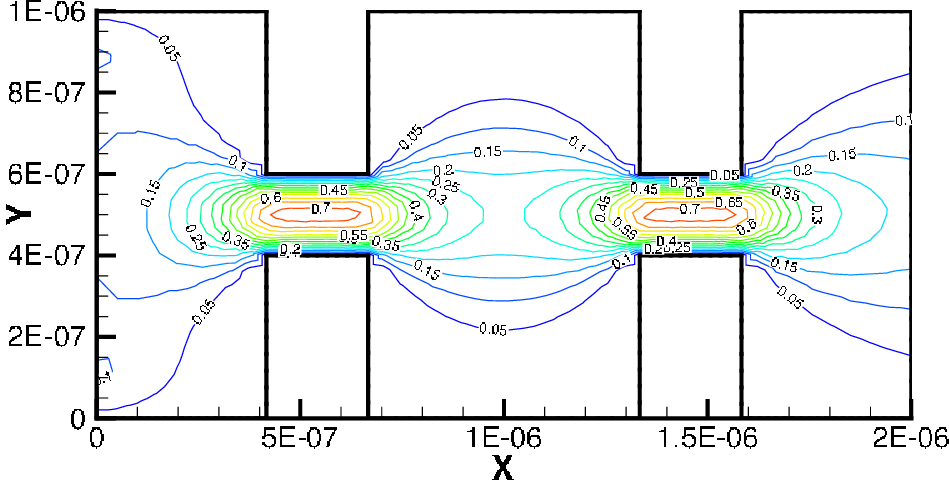}}
	\caption{Counter plots of $u_1$ with different roughness height $h$.}
	\label{fig:U}
\end{figure}

We give the effect of the roughness height  on the electro-osmotic flow in a microchannel.  We choose $h=0.1 H, 0.2H, 0.24H$ and $0.4H$, the numerical results are presented in Figure \ref{fig:01h}, \ref{fig:02h}, \ref{fig:03h} and \ref{fig:04h} at $t=2e-6$, respectively. From Figure \ref{fig:01h}, we can see that the initial condition affects the numerical results strongly. Figure \ref{fig:02h} gives the numerical results with $h=0.2H$. It shows that the roughness height $h$ inflects the numerical results strongly.
Numerical results for $h=0.3H$ and $0.4H$ are given in Figure \ref{fig:03h} and \ref{fig:04h}. In Figure \ref{fig:U}, we give the counter plots of $u_1$ with different roughness height $h$.  As the numerical results of Wang et al \cite{WWC}, it can be seen that the velocity maximum increases with the roughness height.



\end{document}